\documentclass{amsart}
\usepackage{amsmath,amsthm,amssymb}

\theoremstyle{plain}
\newtheorem{theorem}{Theorem}[section]
\newtheorem{lemma}[theorem]{Lemma}

\newtheorem{cor}[theorem]{Corollary}

\newtheorem{prop}[theorem]{Proposition}

\newtheorem*{theoremunnumbered}{Theorem}

\theoremstyle{definition}
\newtheorem{definition}[theorem]{Definition}
\newtheorem{example}[theorem]{Example}
\newtheorem{examples}[theorem]{Examples}

\theoremstyle{remark}
\newtheorem{remark}[theorem]{Remark}
\newtheorem*{remarkunnumbered}{Remark}
\newtheorem*{historicalremark}{Historical remark}

\newtheorem*{notation}{Notation}

\newcommand{\nc}{\newcommand}
\nc{\hb}{\mathbb}
\nc{\M}{\mathcal}
\nc{\mf}{\mathfrak}
\nc{\mbf}{\mathbf}
\nc{\cal}{\mathcal}
\nc{\DMO}{\DeclareMathOperator}
 
\nc{\rom}{\textup}

\nc{\cA}{{\M A}} \nc{\cB}{{\M B}} \nc{\cC}{{\M C}} \nc{\cD}{{\M D}}
\nc{\cE}{{\M E}} \nc{\cF}{{\M F}} \nc{\cG}{{\M G}} \nc{\cH}{{\M H}}
\nc{\cI}{{\M I}} \nc{\cJ}{{\M J}} \nc{\cK}{{\M K}} \nc{\cL}{{\M L}}
\nc{\cM}{{\M M}} \nc{\cN}{{\M N}} \nc{\cO}{{\M O}} \nc{\cP}{{\M P}}
\nc{\cQ}{{\M Q}} \nc{\cR}{{\M R}} \nc{\cS}{{\M S}} \nc{\cT}{{\M T}}
\nc{\cU}{{\M U}} \nc{\cV}{{\M V}} \nc{\cW}{{\M W}} \nc{\cX}{{\M X}}
\nc{\cY}{{\M Y}} \nc{\cZ}{{\M Z}}
        
\nc{\Aa}{{\hb A}} \nc{\Cc}{{\hb C}} \nc{\Gg}{{\hb G}}
\nc{\Nn}{{\hb N}} \nc{\Pp}{{\hb P}}
\nc{\Qq}{{\hb Q}} \nc{\Rr}{{\hb R}} \nc{\Zz}{{\hb Z}}

\DMO*{\trdeg}{trdeg}
\DMO*{\spec}{Spec}
\DMO{\RU}{RU}
\DMO{\RM}{RM}
\DMO{\RC}{RC}
\DMO{\RH}{RH}
\DMO{\RD}{RD}
\DMO{\hgt}{ht}
\DMO{\ord}{ord}
\DMO{\rank}{rank}
\DMO{\ann}{ann}
\DMO{\lc}{lc}
\DMO*{\proj}{Proj}
\DMO{\On}{On}
\DMO{\Ant}{Ant}
\DMO{\Dec}{Dec}
\DMO{\Bad}{Bad}
\DMO{\Incr}{Incr}
\DMO{\Decr}{Decr}
\DMO{\ot}{o}
\DMO{\rk}{rk}
\DMO{\wdt}{wd}
\DMO{\sat}{sat}
\DMO{\reg}{reg}
\DMO{\Greg}{Greg}
\nc{\concat}{\, \widehat{\ }\,}
\DMO{\GL}{GL}

\nc{\ra}{\rightarrow}
\nc{\xra}{\xrightarrow}

\nc{\de}{\delta} \nc{\td}[2]{\trdeg{{#1}/{#2}}}
\nc{\dtd}[2]{\trdeg_{\delta}{({#1}/{#2})}}
\nc{\gen}[2]{ {#1} \langle {#2} \rangle }
\DMO{\dspec}{\ensuremath{\Delta}-Spec}
\DMO{\ddim}{\ensuremath{\Delta}-dim}
\nc{\form}{\Omega}
\nc{\mDCF}{\text{$m$-DCF$_0$}}

\nc{\tuple}[2]{{#1},\ldots,{#2}} \nc{\ptu}[2]{{#1}:\ldots:{#2}}
\nc{\pbrac}[1]{\langle {#1} \rangle}

\nc{\maps}[3]{{#1}\!\colon\!{#2}\ra{#3}}
\nc{\map}[2]{{#1}\ra {#2}} \nc{\res}[2]{{#1} |_{#2}}
\nc{\imbed}{\hookrightarrow}

\nc{\set}[1]{\{#1\}}
\nc{\ptb}[1]{\langle #1 \rangle}
\nc{\abs}[1]{\lvert#1\rvert}

\nc{\bn}{\binom}
\nc{\ini}{\sqsubseteq}
\nc{\propini}{\sqsubset}
\nc{\ext}{\sqsupseteq}
\nc{\KB}{{\operatorname{KB}}}
\nc{\length}{{\operatorname{length}}}
\nc{\lex}{{\operatorname{lex}}}
\nc{\supp}{{\operatorname{supp}\,}}
\nc{\initial}{\operatorname{in}}
\nc{\lm}{\operatorname{lm}}

\nc{\eq}{\text{eq}}
\nc{\tp}{\operatorname{tp}}
\nc{\WKL}{\text{\sf WKL}}

\numberwithin{equation}{section}

\title{Orderings of Monomial Ideals}

\author{Matthias Aschenbrenner}
\address{Department of Mathematics\\
University of California at Berkeley\\
Evans Hall \\
Berkeley, CA 94720 }
\email{maschenb@math.berkeley.edu}

\author{Wai Yan Pong}
\address{Department of Mathematics\\
California State University Dominguez-Hills \\
1000 E. Victoria Street \\
Carson, CA 90747 }
\email{pong@math.csudh.edu}

\thanks{Wai Yan Pong partially supported by NSF Grant
    DMS-0070179.}

\date{April 2003}

\begin{document}

\begin{abstract}
We study the set of monomial ideals in a polynomial ring
as an ordered set, with the ordering given by reverse inclusion. We give a 
short proof of the fact that every antichain of monomial ideals is
finite. Then we investigate ordinal invariants for the complexity of this
ordered set.
In particular, 
we give an interpretation of the height function in terms of the Hilbert-Samuel
polynomial, and we compute upper and lower bounds on the maximal order type.
\end{abstract}

\maketitle

\section*{Introduction}

\noindent
Monomial ideals (that is, ideals generated by monomials) in polynomial or
power series rings play an important role in commutative algebra and algebraic
combinatorics, both from a 
theoretical and a practical perspective. The reason for this is that
more often than not problems about
arbitrary ideals can be reduced to the special case of 
monomial ideals, and hence to questions of a combinatorial
nature. Conversely, monomial ideals may be used to make algebra out of
combinatorics, see, e.g., \cite{Sturmfels}. 
The link between monomial ideals and arbitrary ideals is
provided by the theory of Gr\"ob\-ner bases 
(or standard bases), see, e.g.,
\cite{BeckerWeispfenning}.  

Let $K$ be a field and $R=K[X]=K[X_1,\dots,X_m]$ the ring of polynomials 
in indeterminates $X=\{X_1,\dots,X_m\}$ with coefficients from $K$.
We employ the usual multi-index notation $X^\nu=X^{\nu_1}\cdots X^{\nu_m}$
for monomials, where $\nu=(\nu_1,\dots,\nu_m)$ is an $m$-tuple
of non-negative integers.
Divisibility of monomials in $R$ has the following 
well-known finiteness property:
\begin{quote}
{\it Every sequence $X^{\nu^{(1)}},X^{\nu^{(2)}},\dots,
X^{\nu^{(n)}},\dots$ 
of monomials in $R$ such that
$X^{\nu^{(i)}}$ does not divide $X^{\nu^{(j)}}$, for all $i<j$, is finite.}
\end{quote}
This equally elementary and fundamental
fact, commonly known as ``Dickson's Lemma'', is arguably
``the most frequently rediscovered mathematical
theorem.'' (\cite{BeckerWeispfenning}, p.~184.) 
Among other things, it implies
Hilbert's Basis Theorem, with its numerous consequences.
Recently Diane
Maclagan \cite{Maclagan} proved the following more general result:
\begin{quote}
{\it Every sequence $I^{(1)},I^{(2)},\dots,I^{(n)},\dots$
of monomial ideals in $R$ such that $I^{(i)}\not\supseteq I^{(j)}$
whenever $i<j$ is finite.}
\end{quote}
She also showed how this can be used to give short proofs of
several other finiteness statements like the
existence of a universal Gr\"ob\-ner basis of an ideal in $R$ 
and the finiteness of the number
of atomic fibers of a matrix with non-negative integer entries. 
Galligo's theorem on the existence of generic initial
ideals can also been seen as a consequence of this principle, as
can the upper semi-continuity of fiber dimension (see \cite{maschenb-pong-rd})
and Sit's theorem \cite{Sit} on 
the well-orderedness of the set of Hilbert polynomials under
eventual dominance; for the latter see Section~\ref{HilbertSection} of
the present paper.
It is these remarkable
applications which seem to warrant a further investigation into
combinatorial finiteness phenomena of monomials in $R$. The theory of
Noetherian ordered sets  provides a convenient 
axiomatic framework for this:
Let $(S,\leq)$ be an ordered set, i.e., $S$ is a set and
$\leq$ is a (partial) ordering on
$S$. We call $(S,\leq)$ {\em Noetherian}\/ if every
sequence $s_1,s_2,\dots,s_n,\dots$ in $S$ such that $s_i\not\leq s_j$ for
all $i<j$ is finite.
Dickson's Lemma may then be rephrased as saying that the set
of monomials under divisibility is Noetherian, and
Maclagan's principle just expresses that the set 
$\cM_m$ of monomial ideals in $K[X]$ ordered
by {\sl reverse}\/ inclusion is Noetherian.
Noetherian orderings are usually called ``well-partial-orderings'' or
``well-quasi-orderings'' in the
literature (see, e.g., \cite{Kruskal}). We follow a proposal by
Joris van der Hoeven \cite{vdH:phd} and use the more concise 
(and perhaps more suggestive) term ``Noetherian''.
Noetherian ordered sets play an important role 
in such diverse fields as asymptotic differential algebra \cite{vdH:phd},
Ramsey theory \cite{Kriz-Thomas}, 
theoretical computer science \cite{Dershowitz}, and  
proof theory \cite{Gallier}.

The purpose of this paper is to 
study some aspects of the set of monomial ideals of $K[X]$
from the point of view of combinatorial set theory. In 
Section~\ref{NoetherianOrderedSetsSection}, after reviewing some basic
facts about Noetherian ordered sets, we first give a quick proof of
Maclagan's result. We also indicate a certain generalization, dealing with 
direct products of Noetherian ordered sets (Proposition~\ref{F-Product}), 
which was
stated without proof in \cite{Maclagan} and attributed there to
Farley and Schmidt. 

The complexity of a Noetherian ordered set $(S,\leq)$ 
can be measured in terms of certain ordinal-valued invariants.
We recall their definitions and basic properties in 
Section~\ref{InvariantsSection}. Here is one example:
There always exists a chain in $S$ having maximal possible order type, 
called the {\em height}\/ of $S$; for $x\in S$, the height of the
Noetherian ordered set $S^{\not\geq x}:=\{s\in S:s\not\geq x\}$ is
called the {\em height}\/ of $x$ (in $S$).
From a result of Bonnet and Pouzet
\cite{Bonnet-Pouzet} we deduce that the height of $({\cal M}_m,\supseteq)$
is $\omega^m+1$. In Section~\ref{HilbertSection} we give an interpretation
of the height of a monomial ideal $I$ in terms of the
Hilbert-Samuel polynomial of $M=R/I$.
Recall that for every finitely generated graded $R$-module 
$M=\bigoplus_{s\in\Nn} 
M_s$,  the function which associates with $s\in\Nn$
the dimension of the $K$-vector space $M_s$ agrees, 
for all sufficiently large $s$,
with a polynomial in $\Qq[T]$,
called the Hilbert polynomial of $M$, see \cite{Bruns-Herzog},
Chapter~4. It follows that the function $s\mapsto\dim_K M_{\leq s}$, where
$M_{\leq s}:=\bigoplus_{i=0}^s M_i$, also ultimately agrees with a
polynomial in $\Qq[T]$, which is called the Hilbert-Samuel polynomial
of $M$.
We let $\cI_m$ denote the set of homogeneous
ideals of $R$, considered as a (partially) 
ordered set, with the ordering given 
by reverse inclusion. Given $I\in\cI_m$ we denote the Hilbert-Samuel
polynomial of $R/I$ by $p_I$.  
We totally order $\cS_m=\{p_{I}:I\in\cI_m\}$ by 
eventual dominance: $p_I \leq p_J$ if and only if
$p_{I}(s) \leq p_{J}(s)$ for
all sufficiently large $s$. The map $p\colon\cI_m\to\cS_m$ that maps $I$ to
$p_{I}$ is strictly increasing. 
It is well-known that the map taking each finite $R$-module to its
Hilbert polynomial is the universal additive function on finite 
$R$-modules
which is zero on modules of finite length. (See \cite{Eisenbud},
Section~19.5 for a precise statement.)
The following theorem, proved in
Section~\ref{HilbertSection} below, is in a similar spirit; it
shows that $p$ is universal among
strictly increasing surjections defined on the
ordered set $\cI_m$.
 
\begin{theoremunnumbered}
For every strictly increasing surjection $\varphi\colon\cI_m\to S$, 
where $S$ is
any totally ordered set, there exists a strictly increasing
map $\psi\colon \cS_m\to S$ with $\psi\circ p \leq \varphi$.
\end{theoremunnumbered}

\noindent
Every total ordering extending the ordering $\leq$ of a Noetherian ordered set
$(S,\leq)$ is a well-ordering. This fact gives rise to another invariant of 
$(S,\leq)$: by a theorem of de Jongh and Parikh \cite{dJP}, there
exists a total ordering extending $\leq$ of maximal possible order type,
which we call the {\em type}\/ $\ot(S,\leq)$ of $(S,\leq)$. In 
Section~\ref{OrderingsOfMonomialIdeals} we obtain upper and lower bounds
on $\ot(\cM_m,\supseteq)$: We show that $$\omega^{\omega^{m-1}}+1 \leq 
\ot(\cM_m,\supseteq)
\leq \omega^{\omega^{m+1}}.$$
The proof of the upper bound involves a generalization of a result of
van den Dries and Ehrlich \cite{vdDries-Ehrlich-Erratum},
\cite{vdDries-Ehrlich} on the order type of submonoids of ordered
abelian groups. The lower bound is established by studying a particularly
useful total ordering on ${\cal M}_m$ extending $\supseteq$, inspired by
the Kleene-Brouwer ordering of recursion theory and Kolchin's rankings of
characteristic sets. Both bounds still leave much room for improvement.

We should mention that Dickson's Lemma and Maclagan's principle are only
the first two levels of an infinite hierarchy of finiteness principles:
$\Nn^m$ is ``better-quasi-ordered.'' (This was first shown
by Nash-Williams \cite{NW}.) We refer to \cite{maschenb-hemmecke} 
for the definition of ``better-quasi-ordered'' set and
applications of these more general finiteness properties.

\subsection*{Acknowledgments.}
The first author would like to thank Bernd Sturmfels, whose questions on the
total orderings of monomial ideals inspired this work, and Andreas
Weiermann for an e-mail exchange around the topics of
this paper.

\subsection*{Notations and conventions.}
The cardinality of
a finite set $S$ is denoted by $\abs{S}$.
We let $m,n,\dots$ range over $\Nn:=\{0,1,2,\dots\}$.
For any set $U$, let $U^{\ast}=\bigcup_{n\in\Nn} U^n$ 
denote the set of finite
sequences of elements of $U$. Here $U^0$ consists of the single
element $\varepsilon$ (the empty sequence). 
(So $\emptyset^{\ast}=\{\varepsilon\}$.)
For an element $a=(a_1,\dots,a_n)\in U^{\ast}$ we call the
natural number $n$ the {\it length}\/ of $a$, denoted by
$\length(a)$.
For $a=(a_1,\dots,a_n)$ and $b=(b_1,\dots,b_m)$ in $U^{\ast}$ 
we write $a\ini b$ ($a$ is a {\it truncation}\/ of the sequence $b$)
if $n\leq m$ and $a=(a_1,\dots,a_n)=(b_1,\dots,b_n)$.
By $ab :=
(a_1,\dots,a_n,b_1,\dots,b_m)$ we denote the {\em concatenation}\/ of
the sequences $a=(a_1,\dots,a_n)$ and $b=(b_1,\dots,b_m)$ in $U^{\ast}$.
If, for example, $a=(a_1)$, we shall also write $a_1b$ instead
of $(a_1) \,  b$. With concatenation as monoid operation,
$U^\ast$ is the free monoid generated by $U$ (with identity $\varepsilon$).
We extend concatenation to subsets of
$U^{\ast}$ in the natural way, for example, $aS=
\bigl\{ab:b\in S\bigr\}$ for $a\in U^{\ast}$ and 
$S\subseteq U^{\ast}$.

\section{Noetherian Ordered Sets}\label{NoetherianOrderedSetsSection}

\noindent
In this section we first review the definitions and basic facts about
Noetherian ordered sets. We then give a short proof that the set
of monomial ideals in $K[X]$ is Noetherian, and outline a generalization. 
%We
%finish with the computation of the maximal length of
%strictly decreasing sequences in the lexicographic order on $\Nn^m$,
%which will be used in Section~\ref{HilbertSection}.

\subsection*{Orderings and ordered sets.}
A {\it quasi-ordering}\/ on a set $S$ is a binary relation $\leq$ on $S$
which is reflexive and transitive;
we call $(S,\leq)$ (or simply $S$, if no confusion is possible) a 
{\it quasi-ordered set.}\/ 
If in addition $\leq$ is antisymmetric,
then $\leq$ is called an {\em ordering,}\/ and the pair $(S,\leq)$
is called an {\it ordered set.}\/ 
If $\leq$ is a quasi-ordering on $S$, then so is the inverse relation $\geq$;
likewise for orderings.
If $x$ and $y$ are elements of a 
quasi-ordered set $S$, we write as usual $x<y$ if $x\leq y$ and $y\not\leq x$.
Given a quasi-ordering $\leq$ on a set $S$ and an equivalence
relation $\sim$ on $S$ which is compatible with $\leq$ in the sense
that $x\leq y \Rightarrow x'\leq y'$ for all $x'\sim x$ and $y'\sim y$,
there is  
a unique ordering $\leq_{S/{\sim}}$
on the set $S/{\sim}=\bigl\{x/{\sim}:x\in S\}$ of 
equivalence classes of $\sim$ such that 
$$x/{\sim} \leq_{S/{\sim}} y/{\sim} \qquad\Longleftrightarrow\qquad x\leq y.$$ 
If $\leq$ is an ordering on $S$, then $\leq_{S/\sim}$ is an ordering
on $S/{\sim}$. For any quasi-ordering $\leq$ on $S$,
the equivalence relation on $S$ defined by
$$x\sim y \qquad\Longleftrightarrow\qquad x\leq y \text{ and } 
y\leq x$$ is compatible with $\leq$, and in this case 
$\leq_{S/\sim}$ is an ordering.
Hence by passing from $(S,\leq)$ to $(S/{\sim},\leq_{S/{\sim}})$ if
necessary, we usually can reduce the study of quasi-orderings to the one
of orderings.
In the following, we shall therefore concentrate on {\sl ordered}\/ sets.

\subsection*{Total orderings and directed orderings.}
We say that an ordering on a set $S$ is {\it total}\/ 
if $x \leq y$ or $y \leq x$ for
all $x , y \in S$. 
An ordering $\leq'$ on a set $S$ 
is said to {\it extend}\/ the ordering $\leq$ on $S$
if $x\leq y \Rightarrow x\leq' y$ for all $x,y\in S$. 
Every ordering on a set $S$ can be extended to a total ordering on $S$.
(Szpilrajn's Theorem; the proof uses the Ultrafilter Axiom.)
An ordering on $S$ is {\it directed}\/ if for any $x,y\in S$ there
exists $z\in S$ with $x\leq z$ and $y\leq z$. Any total ordering is directed.

\subsection*{Maps between ordered sets.}
A map $\varphi\colon S\to T$ 
between ordered sets $S$ and $T$ is called
{\it increasing}\/ if  $$x\leq y 
\Rightarrow \varphi(x)\leq\varphi(y) \quad\text{for all $x,y\in S$}$$ and
{\it decreasing}\/ if 
$$x\leq y 
\Rightarrow \varphi(x)\geq\varphi(y) \quad\text{for all $x,y\in S$.}$$
Similarly, we say that  $\varphi\colon S\to T$ is
{\it strictly increasing}\/ if
$$x < y 
\Rightarrow \varphi(x) < \varphi(y) \quad\text{for all $x,y\in S$}$$ and
{\it strictly decreasing}\/ if 
$$x < y 
\Rightarrow \varphi(x) > \varphi(y) \quad\text{for all $x,y\in S$.}$$
We shall write $\Incr(S,T)$ for the set of all increasing maps $S\to T$
and $\Decr(S,T)$ for the set of all decreasing maps $S\to T$.
A map $\psi\colon S\to T$ is a {\em quasi-embedding}\/ 
of $S$ into $T$ if
$$\psi(x)\leq\psi(y) \Rightarrow x\leq y\quad\text{for all $x,y\in S$.}$$
An increasing quasi-embedding $S\to T$
is called an {\em embedding}\/ of $S$ into $T$. 
Finally, a map $S\to T$ is called an
{\em isomorphism}\/ between $S$ and $T$ if it is increasing and
bijective, and its inverse is also
increasing.

\subsection*{Construction of ordered sets.}
Every set $S$ can be equipped with the {\em trivial}\/ ordering, given by
$x\leq y \Longleftrightarrow x=y$.
There are a number of standard constructions for obtaining new (quasi-)
ordered sets from given ones. 
For example, by restricting the ordering,
any subset of an ordered set can be construed as an
ordered set in its own right. Let us explicitly mention some
of the constructions used below. For this, let
$(S,{\leq_S})$ and $(T,{\leq_T})$ be ordered sets. The
disjoint union of the sets $S$ and $T$ is naturally ordered by the relation
$\leq_S \cup \leq_T$; we shall denote this ordered set by $S\amalg T$.
The cartesian product $S\times T$ of $S$ and $T$ can be
made into an ordered set by means of the {\it product ordering:}\/
$$(x,y) \leq (x',y') \qquad :\Longleftrightarrow\qquad x\leq_S x' \text{ and }
y\leq_T y',$$
or the {\it lexicographic ordering:}\/
$$(x,y) \leq_{\text{lex}} (x',y') \qquad :\Longleftrightarrow\qquad x<_S x' \text{ or }
(x=x'\text{ and }y\leq_T y'),$$
for $(x,y),(x',y')\in S\times T$. Taking $S=T$ yields the 
product ordering and the lexicographic ordering 
on $T^2=T\times T$, and by repeating the construction,
on $T^m$ for any $m>0$. More generally,
if $I$ is any set, then
the set $T^I$ of all functions $I\to T$ is ordered by setting
$$f\leq g \qquad:\Longleftrightarrow\qquad 
f(i)\leq_T g(i) \quad\text{for all $i\in I$.}$$
By restriction this yields orderings on the subsets 
$\Incr(S,T)$ and $\Decr(S,T)$ of $T^S$. If the ordering on $S$ is directed, 
we have
(at least) two other ways of defining a quasi-ordering on $T^S$
which extends the product ordering:
\begin{enumerate}
\item using the {\it lexicographic ordering,}\/ defined by 
\begin{align*}
f \leq_{\text{lex}} g \quad :\Longleftrightarrow\quad 
&\text{$f=g$, or there is $y\in S$ with}\\ &\text{$f(x)=g(x)$ for all 
$x<_Sy$ and $f(y)<_Tg(y)$,}
\end{align*} and
\item using the {\it dominance quasi-ordering,}\/ given by
\begin{align*}
f\preceq g\quad :\Longleftrightarrow\quad 
&\text{there is $y\in S$ with $f(x)\leq_T g(x)$ for all $x\geq_S y$.}
\end{align*} 
\end{enumerate}
If both $\leq_S$ and $\leq_T$ are total, then $\leq_{\text{lex}}$ is total.
In general, the dominance quasi-ordering
is neither antisymmetric (i.e., not an ordering on $T^S$)
nor total.

\begin{example}\label{example-Nn}
We consider $\Nn$ as an ordered set under its usual ordering, and we equip
$\Nn^m$ with 
the product ordering.
For $\nu=(\nu_1,\dots,\nu_m)\in\Nn^m$ we put
$\abs{\nu}=\nu_1+\cdots+\nu_m$ 
(the {\em degree}\/ of $\nu$).
Let $X=\{X_1,\dots,X_m\}$ be distinct indeterminates
and $X^\diamond=\{X^\nu: \nu\in\Nn^m\}$ 
the free commutative monoid generated by $X$, where
$X^\nu:=X_1^{\nu_1}\cdots X_n^{\nu_m}$ for $\nu=(\nu_1,\dots,\nu_m)\in\Nn^m$.
We order $X^\diamond$ by divisibility: 
$$X^\nu\leq X^\mu \quad :\Longleftrightarrow \quad
\mu=\nu+\lambda \text{ for some $\lambda\in\Nn^m$.}$$ 
Then $\nu\mapsto X^\nu\colon\Nn^m\to X^\diamond$ is an
isomorphism of ordered sets. The elements of $X^\diamond$ can be seen 
as {\sl monomials}\/ in the polynomial ring $K[X]=K[X_1,\dots,X_m]$,
where $K$ is a  field. Here, the identity element $\varepsilon$ of $X^\diamond$
is identified with the monomial $1$.
\end{example}

\subsection*{Final segments and antichains.}
A {\it final segment}\/ of an ordered set $(S,\leq)$ is a
subset $F \subseteq S$ such that $$x \leq y \wedge x \in F \Rightarrow y
\in F, \qquad \text{for all $x , y \in S$.}$$ 
(Dually, $I \subseteq S$ is called an
{\it initial segment}\/ if $S \backslash I$ is a final segment.)
Given an arbitrary subset $X$ of $S$, we
denote by $$(X) := \bigl\{y \in S : \exists x \in X\ ( x \leq y ) \bigr\}$$
the final segment {\it generated by $X$.}\/
We construe the set $\cF(S)$ of final segments of $S$ as
an ordered set, with the ordering given by {\sl reverse}\/ inclusion.

\begin{example}\label{example-Nn-2}
Under the isomorphism in Example~\ref{example-Nn}, 
final segments of $\Nn^m$ correspond
to {\it ideals}\/ in the commutative
monoid $X^\diamond$, that is, subsets $I\subseteq X^\diamond$ 
such that $vu\in I$ for all $u\in I$ and $v\in X^\diamond$.
Considering the elements of $X^\diamond$ as monomials in a polynomial ring
$K[X]$ over a field $K$, the ordered set $\cF(\Nn^m)$ becomes isomorphic to the set
of {\sl monomial ideals}\/ of $K[X]$ (that is, ideals of $K[X]$ which are
generated by monomials), ordered by reverse inclusion.
%Henceforth we will also call the elements of $\cF(\Nn^m)$ monomial ideals.
\end{example}

\noindent
We write $x||y$ if $x,y\in S$ are
{\em incomparable,}\/ that is, if $x\not\leq y$ and $y\not\leq x$.
An {\it antichain}\/ of $S$ is a subset $A \subseteq S$
such that any two distinct elements $x$ and $y$ of $A$ are incomparable. 
(For example, 
a generating set of a final segment $F$ of $S$ is a minimal generating
set for $F$ if and only if it is an antichain.)
A subset $C$ of $S$ is called a {\it chain}\/ if 
the restriction of the ordering of $S$ to $C$ is total, that is, if
for all $x,y\in C$ we have $x\leq y$ or $y\leq x$.

\subsection*{Noetherian orderings.}
An ordered set $S$ is {\it well-founded}\/ if there is no
infinite strictly decreasing sequence $x_0 > x_1 > \cdots $ in $S$. 
We say that an
ordered set $S$ is {\it Noetherian}\/ if it is well-founded
and every antichain of $S$ is finite. For example, every
finite ordered set is Noetherian.
Since every antichain of a
totally ordered set consists of at most one
element, a totally ordered set $S$ is Noetherian if and only if it is
well-founded; in this case $S$ is called {\it well-ordered.}\/
For every well-ordered set $S$ there exists a unique ordinal number, called
the {\it order type}\/ $\ot(S)$ of $S$, which is isomorphic to $S$.

%\noindent
%Note that if $S$ is well-founded and $F$ is a final segment of $S$, then 
%the antichain $F_{\min}$ of minimal elements of $F$ is the smallest generating set for
%$F$, that is, $(F_{\min})=F$ and $F_{\min}\subseteq X$ for every $X\subseteq F$ with $(X)=F$.

An infinite sequence 
$x_0,x_{1}, \ldots$ in $S$ is {\it good}\/ if $x_i \leq x_j$ for some $i <
 j$, and {\it bad,}\/ otherwise. (For instance, if $\{x_0,x_1,\dots\}$ is an
antichain, then $x_0,x_1,\dots$ is bad.)
The following characterization of Noetherian orderings is folklore; we
omit the proof. (For the details, see, e.g., \cite{maschenb-hemmecke}.)

\begin{prop}\label{Folk}
The following are equivalent, for an ordered set $S$: 
\begin{enumerate}
\item $S$ is Noetherian.
\item Every infinite sequence $x_0 , x_1 , \ldots$ in $S$ contains an
increasing subsequence.
\item Every infinite sequence $x_0 , x_1 , \ldots$ in $S$ is good.
\item Any final segment of $S$ is finitely generated.
\item $\bigl(\mathcal{F} ( S ),\supseteq\bigr)$ is well-founded \textup{(}i.e., the ascending chain
condition with respect to inclusion holds for final segments of $S$\textup{).}
%\item $\bigl(\mathcal{I} ( S ),\subseteq\bigr)$ is well-founded \textup{(}i.e., the descending chain
%condition with respect to inclusion holds for initial segments of $S$\textup{).}
\item Every total ordering on $S$ which extends $\leq$ is a well-ordering. \qed
\end{enumerate}
\end{prop}

\noindent
The proposition immediately implies:
\begin{examples}\label{Folk-Examples}
Let $S$ and $T$ be ordered sets.
\begin{enumerate}
\item If there exists an increasing surjection $S\to T$, and
$S$ is Noetherian, then so is $T$.
In particular: If $S$ is Noetherian, then 
any ordering on $S$ which extends the given ordering is Noetherian;
if $S$ is Noetherian and 
$\sim$ is an equivalence relation on $S$ which is compatible with
the ordering of $S$, then $S/{\sim}$ is Noetherian.
\item If there exists a quasi-embedding $S\to T$, and $T$
is Noetherian, then $S$ is Noetherian.
In particular, if $T$ is Noetherian, then any subset of $T$ with 
the induced ordering is Noetherian.
\item  If $S$ and $T$ are Noetherian and
$U$ is an ordered set which contains both ordered sets $S$ and $T$,
then $S\cup T$ is Noetherian. 
In particular, it follows that  $S\amalg T$ is Noetherian.
\item If $S$ and $T$ are Noetherian, then so is
$S\times T$ with the product ordering.
Inductively, it follows that 
if the ordered set $S$ is Noetherian, 
then so is $S^m$ equipped with the product ordering, for
every $m$.
In particular, for each $m$, the ordered set $\Nn^m$ is
Noetherian (``Dickson's Lemma''). 
\end{enumerate}
\end{examples} 

\noindent
For future use we also remark:

\begin{lemma}\label{Folk-Lemma}
Let $\varphi\colon S\to T$ be a map between ordered sets
$S$ and $T$, with $S$ Noetherian.
\begin{enumerate}
\item If $\varphi$ is strictly increasing, then $\varphi$ has finite fibers.
\item If $\varphi$ is decreasing, and $T$ is well-founded,
then the image of $\varphi$ is finite.
\end{enumerate}
\end{lemma}

\subsection*{Noetherianity of the set of monomial ideals.}

By Proposition~\ref{Folk}, if $S$ is Noetherian, then the ordered set
$\cF(S)$ of final segments of $S$ is well-founded.
In general, it is {\sl not}\/ true that if $S$ is Noetherian, then 
$\cF(S)$ is Noetherian. A counterexample 
was found by Rado \cite{Rado}. (Rado's example is indeed ``generic''
in the sense that a Noetherian ordered set $S$ contains an isomorphic copy 
of this example if and only if $\cF(S)$ is non-Noetherian; 
see, e.g., \cite{maschenb-hemmecke}.) 
We will now give a short proof of the fact that
the ordered set $\cF(\Nn^m)$ of monomial ideals is Noetherian.
Here is a key observation:

\begin{lemma}\label{Decr-Lemma}
$\cF(S\times T)\cong\Decr\bigl(S,\cF(T)\bigr)$, for ordered sets $S$ and $T$.
\end{lemma}
\begin{proof}
For a final segment $F\in\cF(S\times T)$ let $\varphi_F\colon S\to\cF(T)$ be defined by
$\varphi_F(x)=\bigl\{y\in T:(x,y)\in 
F\bigr\}$ for $x\in S$.
It is straightforward to verify that $\varphi_F$ is decreasing, and
$F\mapsto\varphi_F$ is an isomorphism $\cF(S\times T)\to
\Decr\bigl(S,\cF(T)\bigr)$.
\end{proof}
\noindent
In particular, we have $\cF(\Nn^{m})\cong
\Decr\bigl(\Nn,\cF(\Nn^{m-1})\bigr)$ for $m>0$.
This fact allows us to analyze $\cF(\Nn^m)$ by induction on $m$; it
also makes it necessary to take a closer look at
decreasing maps $\Nn\to\cF(\Nn^{m-1})$.
More generally, for any ordered set $S$, let us use
$S^{(\geq)}$ to denote the set
$\Decr(\Nn,S)$ of all infinite
decreasing sequences $s=(s_0,s_1,\dots)$ of elements 
$s_0\geq s_1\geq\cdots$ of $S$, ordered component-wise (that is, by restriction
of the product ordering on $S^\Nn$ to $S^{(\geq)}$). 

\begin{prop}\label{S-Prop}
If $S$ is Noetherian, then so
is $S^{(\geq)}$. 
\end{prop}
\begin{proof}
The proof is inspired by Nash-Williams' proof
\cite{NW2} of Higman's Lemma (see Lemma~\ref{Higman} below):
Suppose for a contradiction that $S$ is Noetherian and 
$s^{(0)},s^{(1)},\dots$ is a bad sequence in $S^{(\geq)}$; 
we write $s^{(i)}=\bigl(s_0^{(i)},s_1^{(i)},\dots\bigr)$.
Every sequence $s=(s_0,s_1,\dots)$ in $S^{(\geq)}$
becomes eventually stationary; we let
$j(s)$ denote the smallest index $j\in\Nn$ such that $s_j=s_{j+1}=\cdots$.
We may assume that the bad sequence is chosen in such a way that
for every $i$,
$j\bigl(s^{(i)}\bigr)$ is {\it minimal} among the $j(s)$, where $s$ ranges over all
elements of $S^{(\geq)}$
with the property that $s^{(0)},s^{(1)},\dots,s^{(i-1)},s$ can be continued to 
a bad sequence in $S^{(\geq)}$. 
We may further assume that there is
an index $i_0$ such that
$j\bigl(s^{(i)}\bigr)>0$ for all $i\geq i_0$.
Now consider the sequence $s_0^{(i_0)},s_0^{(i_0+1)},\dots$ in $S$. Since
$S$ is Noetherian, there exists an infinite sequence $i_0\leq i_1<i_2<\cdots$
of indices such that $s_0^{(i_1)}\leq s_0^{(i_2)}\leq\cdots$. Put
$t^{(i_k)}:=\bigl(s_1^{(i_k)},s_2^{(i_k)},\dots\bigr)$ for all $k>0$. 
It is now easily seen
that then $s^{(0)},\dots,s^{(i_1-1)},t^{(i_1)},t^{(i_2)},\dots$ is a bad
sequence in $S^{(\geq)}$. But
$j\bigl(t^{(i_k)}\bigr)=j\bigl(s^{(i_k)}\bigr)-1$, contradicting 
the minimality property of our original bad sequence.
\end{proof}
\noindent
The ordered set $\cF(\Nn)$ is clearly well-ordered (of order type $\omega+1$),
hence Noetherian. This is the base case for an induction on $m$, which yields,
using Proposition~\ref{S-Prop}:

\begin{cor}\label{Noetherian-Cor}
The set of monomial ideals in $K[X_1,\dots,X_m]$, ordered by reverse
inclusion, is Noetherian, for any $m$. \qed
\end{cor}
\noindent
See \cite{Maclagan} and \cite{maschenb-hemmecke} for other proofs of this
result. The proof in \cite{Maclagan} uses primary decomposition of 
monomial ideals in $K[X]$; the
proof in \cite{maschenb-hemmecke} is based on Ramsey's Theorem.

\subsection*{Higman's lemma.}
In the following, we will often make use of 
a fundamental fact due to Higman \cite{Higman}.
Let $S$ be an ordered set.
We define an ordering on the set $S^{\ast}$ of finite sequences of
elements of $S$ as follows:
$$
(x_1,\dots,x_m) \leq^\ast (y_1,\dots,y_n) \quad :\Longleftrightarrow \quad
\begin{cases}
&\text{\parbox{160pt}{there exists a strictly increasing function $\varphi\colon
\{1,\dots,m\}\to\{1,\dots,n\}$ such that $x_i\leq y_{\varphi(i)}$ for
all $1\leq i\leq m$.}}\end{cases}
$$
\begin{lemma}\label{Higman}
\rom{(Higman.)}
If $S$ is Noetherian, then 
the ordering $\leq^\ast$ on $S^{\ast}$ is Noetherian.
\end{lemma}
\noindent
The equivalence relation $\sim$ on $S^*$ defined by
$$
(x_1,\dots,x_m) \sim (y_1,\dots,y_n) \quad :\Longleftrightarrow \quad
\begin{cases}
&\text{\parbox{160pt}{$m=n$ \&
there exists a permutation $\sigma$ of 
$\{1,\dots,m\}$ such that $x_i=y_{\sigma(i)}$ for
all $1\leq i\leq m$}}\end{cases}
$$
is compatible with $\leq^\ast$, and hence induces an ordering on 
$S^\diamond:=S^*/{\sim}$, which we denote by $\leq^\diamond$. 
By Higman's Lemma, if $(S,\leq)$ is Noetherian, then so
is $(S^\diamond,\leq^\diamond)$.
For $x=(x_1,\dots,x_m)\in S^*$ we denote by $[x]=[x_1,\dots,x_m]\in S^\diamond$ 
the equivalence class of $x$, and we put $\abs{w}=m$ for $w=[x_1,\dots,x_m]\in
S^\diamond$. 
We may think of the elements of $S^*$ as
{\em non-commutative words}\/ in the alphabet $S$, and of the elements of
$S^\diamond$ as {\em commutative words}\/ in $S$.
Note that $S^\diamond$, with concatenation of commutative words, is
the free commutative monoid generated by $S$.

\begin{remark}\label{Sdiamond}
We identify $S$ with a subset of $S^\diamond$ in a natural way.
Let us call a total ordering $\leq$ of $S^\diamond$ 
extending the ordering on $S$ a {\it term ordering}\/
if $\varepsilon\leq v$ and 
$v\leq w \Rightarrow sv\leq sw$, for all $v,w\in S^\diamond$ and
$s\in S$. Then, for $v,w\in S^\diamond$:
$$v\leq^\diamond w \quad\Longleftrightarrow\quad
v\leq w \text{ for all term orderings $\leq$ of $S^\diamond$.}$$
This follows, e.g., from the Artin-Schreier theory of formally real fields 
applied to the quotient field 
of the monoid
ring $\Qq[S^\diamond]$; see, e.g., \cite{BCR}, Proposition~1.1.10. 
%(Probably an analogous characterization can be given of
%$\leq^\ast$.) 
In the case $S=X=\{X_1,\dots,X_m\}$ 
ordered such that $X_1<\cdots<X_m$, monomial ideals of $K[X]$
whose corresponding final segment $E\in\cF(\Nn^m)$ is also
a final segment with respect to $\leq^\diamond$ are called
{\em strongly stable.}
\end{remark}

\begin{remark}\label{multiset}
If $(S,\leq)$ is Noetherian, then so is any ordering on $S^\diamond$ which
extends $\leq^\diamond$. An important example is
the {\em multiset ordering}\/ on $S^\diamond$, defined by
$$s \preceq\!\!\!\preceq t \qquad :\Longleftrightarrow\qquad
\begin{cases}
&\text{\parbox{180pt}{
$s=t$, or for each $i\in\{1,\dots,n\}$ there exists $j\in\{1,\dots,l\}$
with $s_{m+i}\leq t_{m+j}$,}}
\end{cases}
$$
where we write
$s=[s_1,\dots,s_{m+n}]$, $t=[t_1,\dots,t_{m+l}]$
with $s_i=t_i$ for $1\leq i\leq m$ and
$\{s_{m+1},\dots,s_{m+n}\}\cap\{t_{m+1},\dots,t_{m+l}\}=\emptyset$.
To show that $s\leq^\diamond t\Rightarrow s\preceq\!\!\!\preceq t$
we proceed by induction on $m$, the case $m=0$ being trivial. Suppose $m>0$
and $s\leq^\diamond t$. Then there exists an injective function
$\varphi\colon \{1,\dots,m+n\}\to\{1,\dots,m+l\}$
such that $s_k \leq t_{\varphi(k)}$ for all $1\leq k\leq m+n$.
Given $i$ with $1\leq i\leq n$, we have to show that $s_{m+i}\leq t_{m+j}$ for some 
$j$ with $1\leq j\leq l$.
Injectivity of $\varphi$ implies that there exists $r\in\Nn$ such that
$\varphi(m+i),\varphi^2(m+i)=\varphi\bigl(\varphi(m+i)\bigr),\dots,\varphi^{r}(m+i)\in \{1,\dots,m\}$ and
$\varphi^{r+1}(m+i)>m$. Hence $j=m-\varphi^{r+1}(m+i)$ works. 
\end{remark}

\subsection*{A generalization.}
The ideas used to establish Noetherianity of $\cF(\Nn^m)$ above can be
generalized somewhat to give a proof of the following fact:

\begin{prop}\label{F-Product}
Let $S$ and $T$ be ordered sets. If $\cF(S)$ and $\cF(T)$ are 
Noetherian, then $\cF(S\times T)$ is Noetherian.
\end{prop}

\noindent
Recall that for ordered sets $S$ and $T$, we use $\Decr(S,T)$ to denote the
set of all decreasing maps $S\to T$, ordered point-wise: $\varphi\leq\psi$
if and only if $\varphi(x)\leq\psi(x)$ for all $x\in S$.

\begin{lemma}\label{Decr-Lemma-2}
Let $S\neq\emptyset$ be an ordered set. The following are equivalent:
\begin{enumerate}
\item For every Noetherian ordered set $T$, $\Decr(S,T)$ is Noetherian.
\item For some Noetherian ordered set $T$ with 
$\abs{T}>1$, $\Decr(S,T)$ is Noetherian.
\item $\Decr(S,2)$ is Noetherian. \rom{(}Here $2=\{0,1\}$ with $0<1$.\rom{)}
\item $\cF(S)$ is Noetherian.
\end{enumerate}
\end{lemma}
\begin{proof}
The implication ``(1)~$\Rightarrow$~(2)'' is trivial. For
``(2)~$\Rightarrow$~(3)'', let $T$ be a Noetherian ordered set with
more than $1$ element, such that  $\Decr(S,T)$ is Noetherian.
If $T$ is an antichain, then so is $\Decr(S,T)$. Hence
$\Decr(S,T)$, and thus $S$, are finite. Therefore $\Decr(S,2)$ is finite,
hence Noetherian.
If $T$ is not an antichain, then there exists a  quasi-embedding
$\Decr(S,2)\to\Decr(S,T)$, showing that $\Decr(S,2)$ is Noetherian.
For ``(3)~$\Rightarrow$~(4)'',
note that for every $F\in\cF(S)$, the function
$\varphi_F\colon S\to 2=\{0,1\}$ given by
$$\varphi_F(x) = \begin{cases}
0 &\text{if $x\in F$} \\
1 &\text{if $x\notin F$}
\end{cases}$$
(the characteristic function of $S\setminus F$) is decreasing, and
$\varphi_F \leq\varphi_G$ if and only $F\supseteq G$, for $F,G\in\cF(S)$.
For ``(4)~$\Rightarrow$~(1)'', suppose that $\cF=\cF(S)$ is Noetherian, and
let $T$ be a Noetherian ordered set. Then the image of every decreasing
map $S\to T$ is finite. (Lemma~\ref{Folk-Lemma},~(2).)
For $\varphi\in\Decr(S,T)$ and $y\in T$, the inverse
image $\varphi^{-1}(T_y)$ of the initial segment 
$T_y=T\setminus (y)=\{z\in T:z\not\geq y\}$ of $T$ is a final segment of
$S$.
We define a map $$\Psi\colon\Decr(S,T)\to (T\times\cF)^{\diamond}$$ 
as follows:
Given $\varphi\in\Decr(S,T)$ let $y_1,\dots,y_k\in T$ be the distinct elements
of $\varphi(S)$; hence $\varphi^{-1}(T_{y_i}) \supseteq\varphi^{-1}(T_{y_j})$
if $y_i>y_j$. We put
$$\Psi(\varphi) = \bigl[ \bigl(y_1,\varphi^{-1}(T_{y_1})\bigr),\dots,
\bigl(y_k,\varphi^{-1}(T_{y_k})\bigr)\bigr].$$ 
One checks easily that $\Psi$ is a
quasi-embedding, where $(T\times\cF)^{\diamond}$ is equipped with the
ordering $\leq^\diamond$. Since
$(T\times\cF)^{\diamond}$ is Noetherian, so is 
$\Decr(S,T)$, as desired.
\end{proof}
\begin{remark}
In \cite{Brookfield}, Lemma~2.12 it is shown that if $S$ is Noetherian and
$T$ is well-founded, then $\Decr(S,T)$ is well-founded.
\end{remark}

\noindent
Proposition~\ref{F-Product} now follows from
Lemmas~\ref{Decr-Lemma} and \ref{Decr-Lemma-2}.

\section{Invariants of Noetherian Ordered Sets}\label{InvariantsSection}

\noindent
Here, we introduce certain ordinal numbers associated
to Noetherian ordered sets, and establish (or recall) some fundamental
facts about them. After some preliminaries concerning ordinal arithmetic
we discuss the {\it height}\/ of a Noetherian, or more generally
well-founded, ordered set. We then define the {\it type}\/ and {\it width}\/ 
of a Noetherian ordered set $S$ in terms of the heights of certain well-founded
trees associated to $S$, and we state some of the basic relations between
them. We relate another invariant (the {\it minimal order type}\/ of $S$) 
with the height, and we compute the height of 
a certain modification of one of the trees associated with $S$.
We finish by computing these invariants for the ordered set $S=\Nn^m$.

\subsection*{Natural sum and product of ordinals.}

We denote the class of all ordinal numbers by $\On$. We identify each
ordinal with the set of its predecessors; thus $\alpha<\beta$ is
synonymous with $\alpha\in\beta$, for $\alpha,\beta\in\On$. The smallest
infinite ordinal is denoted by $\omega$. 
Any non-zero ordinal $\alpha$ can be expressed
in the form
$$\alpha=\omega^{\gamma_1}a_1 + \omega^{\gamma_2}a_2 + \cdots + \omega^{\gamma_n}a_n,$$
where $\gamma_1 > \gamma_2 > \cdots > \gamma_n$ are ordinals and
$a_1,\dots,a_n\in\Nn$. If we require in addition that the $a_i$ are positive,
then this representation of $\alpha$ is unique and called the
{\it Cantor normal form}\/ of $\alpha$.
The (Hessenberg) {\it natural
sum $\alpha\oplus\beta$}\/ of two ordinals
\begin{equation}\label{Ordinal-1}
\alpha=\omega^{\gamma_1}a_1 + \omega^{\gamma_2}a_2 + \cdots + \omega^{\gamma_n}a_n
\end{equation}
and
\begin{equation}\label{Ordinal-2}
\beta=\omega^{\gamma_1}b_1 + \omega^{\gamma_2}b_2 + \cdots + \omega^{\gamma_n}b_n
\end{equation}
(where $a_i,b_j\in\Nn$) is defined by 
$$\alpha\oplus\beta = \omega^{\gamma_1}(a_1+b_1)+\omega^{\gamma_2}(a_2+b_2)
+\cdots+\omega^{\gamma_n}(a_n+b_n).$$
In particular, we have $0\oplus\alpha=\alpha\oplus 0=\alpha$ for all
$\alpha\in\On$. The operation $\oplus$ on $\On$ is associative, commutative
and strictly increasing when  one of the arguments is
fixed: $\alpha<\beta \Rightarrow \alpha\oplus\gamma <\beta\oplus\gamma$, for
all $\alpha,\beta,\gamma\in\On$. It follows that $\oplus$ is 
cancellative (i.e., $\alpha\oplus\gamma=\beta\oplus\gamma \Rightarrow
\alpha=\beta$).
%By transfinite induction, we 
%define the natural sum $\bigoplus_{\gamma<\delta} \alpha_\gamma$ of
%a sequence $\{\alpha_\gamma\}_{\gamma<\delta}$ in $\On$ indexed by all
%ordinals less than some ordinal $\delta$: If $\delta=0$, we put
%$\{\alpha_\gamma\}_{\gamma<\delta}=0$, in the successor case we define
%$\{\alpha_\gamma\}_{\gamma<\delta+1}=\{\alpha_\gamma\}_{\gamma<\delta}\oplus
%\alpha_{\delta}$, and if $\delta$ is a limit ordinal, we put
%$\{\alpha_\gamma\}_{\gamma<\delta}=\sup_{\delta'<\delta} 
%\{\alpha_\gamma\}_{\gamma<\delta'}$.
The {\it natural product}\/ of ordinals $\alpha$ and
$\beta$ written as in \eqref{Ordinal-1} and \eqref{Ordinal-2} above,
respectively, is given by
$$\alpha\otimes\beta = \bigoplus_{i,j} \omega^{\gamma_i\oplus\gamma_j}a_ib_j.$$
The natural product, too, is associative, commutative, and
strictly increasing in both arguments (hence cancellative).
The distributive law for $\oplus$ and $\otimes$
holds: $\alpha\otimes(\beta\oplus\gamma)=
(\alpha\otimes\beta)\oplus(\alpha\otimes\gamma)$. We refer to
\cite{Bachmann} for more information about the natural operations on $\On$.
Below, we will make use of the identity
\begin{equation}\label{NatSum}
\alpha\oplus\beta = \sup\bigl\{\alpha'\oplus\beta+1, \alpha\oplus\beta'+1:
\alpha'<\alpha,\beta'<\beta\bigr\}
\end{equation}
for ordinals $\alpha,\beta$.

\subsection*{Height functions.}
Let $S$ be a well-founded ordered set.
For a proof of
the following lemma see, e.g., \cite{Fraisse}, \S 2.7.
By convention, $\sup\emptyset := 0\in\On$.

\begin{lemma}\label{hgt-char}
The following are equivalent, for a map $h\colon S\to\On$:
\begin{enumerate}
\item $h$ is strictly increasing, and if $h'\colon S\to\On$ is
a strictly increasing map, then $h(x)\leq h'(x)$ for all $x\in S$.
\item $h$ is strictly increasing, and $h(S)$ is an 
initial segment of $\On$.
\item
$h(x)=\sup\bigl\{h(y)+1 : y<x\bigr\}$ for all $x\in S$.
\end{enumerate}
\end{lemma}

\noindent
There exists a unique function $h=\hgt_S\colon S\to\On$ satisfying
the equivalent conditions of the lemma. If $S$ is clear from the context,
we shall just write $\hgt$ for $\hgt_S$.
The ordinal $\hgt(x)$ is called the {\it height}\/ of $x$ in $S$, and
the image
$$\hgt(S)=\sup\bigl\{\hgt(x)+1:x\in S\bigr\} \in\On$$ of $S$ under $\hgt$
is called the 
{\it height}\/ of the well-founded ordered set $S$. 
The height of $S$ is the
smallest ordinal $\alpha$ such that there exists a strictly increasing
function $S\to\alpha$. Equivalently:
\begin{equation}\label{height-sup}
\hgt(S) = \sup\bigl\{ \ot(C) : \text{$C\subseteq S$ chain}\bigr\}.
\end{equation}
In particular, if $S$ is well-ordered, then the height
$\hgt(S)$ of $S$ agrees with the order type $\ot(S)$ of $S$, and 
the height function $\hgt\colon S\to\hgt(S)$ is the unique isomorphism 
$S\to\ot(S)$. 

\begin{lemma}\label{hgt-properties}
Let $S$ and $T$ be non-empty well-founded ordered sets. 
\begin{enumerate}
\item If there exists a strictly increasing map $S\to T$, 
then $\hgt(S)\leq\hgt(T)$.
\item $\hgt(S\amalg T)=\max\bigl(\hgt(S),\hgt(T)\bigr)$.
\item $\hgt_{S\times T}(s,t)=\hgt_S(s)\oplus\hgt_T(t)$ for
all $s\in S$, $t\in T$.
\item $\max\bigl(\hgt(S),\hgt(T)\bigr)\leq
\hgt(S\times T)<\hgt(S)\oplus\hgt(T)$.
\end{enumerate}
\end{lemma}
\begin{proof}
Part (1) follows immediately from (1) in the previous lemma. 
Part (2) is obvious. For a proof
of (3) and (4) see \cite{Fraisse}, 4.8.3. 
\end{proof}
\noindent
The following lemma will be used in Section~\ref{HilbertSection} below.
Let us call a map $\varphi\colon S\to T$ between ordered sets $S$ and $T$
{\it non-decreasing}\/ if $x<y\Rightarrow \varphi(x)\not\geq\varphi(y)$,
for all $x,y\in S$. (So if $T$ is totally ordered, then non-decreasing is
equivalent to strictly increasing.)

\begin{lemma}\label{Aux-Height-Lemma}
Let $S$ be a well-founded ordered set. The following are equivalent,
for a strictly increasing surjection $h\colon S\to T$, where $T$ is a
well-ordered set:
\begin{enumerate}
\item $\hgt_S=\hgt_T\circ h$.
\item For every ordered set $T'$ and strictly increasing map
$h'\colon S\to T'$ there exists a non-decreasing map
$\psi\colon T\to T'$ such that $\psi\circ h\leq h'$.
\item For every totally ordered set $T'$ and strictly increasing map
$h'\colon S\to T'$ there exists a strictly increasing map
$\psi\colon T\to T'$ such that $\psi\circ h\leq h'$.
\end{enumerate}
Given strictly increasing surjections $h\colon S\to T$ and
$h'\colon S\to T'$ satisfying these conditions, with totally ordered sets
$T$ and $T'$, there exists an
isomorphism $\varphi\colon T\to T'$ such that $h'=\varphi\circ h$.
\end{lemma}
\begin{proof}
For (1)~$\Rightarrow$~(2), suppose $\hgt_S=\hgt_T\circ h$, and let
$T'$ be an ordered set, $h'\colon S\to T'$ strictly increasing.
Then $\hgt_{T'}\circ h'$ is strictly increasing, hence
$\hgt_T\circ h =\hgt_S\leq \hgt_{T'}\circ h'$, and 
$\hgt(S)=\hgt_T(h(S))$ is an initial
segment of the range of $\hgt_{T'}\circ h'$. We let $\iota$ denote the
natural inclusion $\hgt(S)\imbed (\hgt_{T'}\circ h')(S)$.
For $y\in\hgt(T')$ let $\hgt_{T'}^{-1}(y)\in T'$ denote the smallest
$x\in T'$ such that $\hgt_{T'}(x)=y$. Then
$\hgt_{T'}\bigl(\hgt_{T'}^{-1}(y)\bigr)=y$ and 
$\hgt_{T'}^{-1}\bigl(\hgt_{T'}(x)\bigr)\leq x$ for all $x\in T',y\in\hgt(T')$.
The map 
$\psi=\hgt_{T'}^{-1}\circ\iota\circ\hgt_T$ is non-decreasing and satisfies
$\psi\bigl(h(s)\bigr) \leq h'(s)$ for all $s\in S$ as required.

The implication (2)~$\Rightarrow$~(3) is trivial.
Suppose that $h$ satisfies (3). Then for every strictly increasing
function $h'\colon S\to\On$ there exists a strictly increasing map
$\psi\colon T\to \On$ such that $\psi\circ h\leq h'$. Hence
$\psi'=\psi\circ\hgt_T^{-1}$ is a strictly increasing embedding of the
ordinal $\hgt(T)$ into $\On$ with $\psi'\circ (\hgt_T\circ h) \leq h'$.
Since $\psi'(\alpha)\geq\alpha$ for all ordinals $\alpha<\hgt(T)$, it follows
that $\hgt_T\circ h\leq h'$. Hence $\hgt_T\circ h=\hgt_S$.

For the second part, let $T$ and $T'$ be totally ordered sets, and let
$h\colon S\to T$ and
$h'\colon S\to T'$ be strictly increasing surjections
satisfying these equivalent conditions. Then 
$\hgt_T\circ h=\hgt_S=\hgt_{T'}\circ h'$, hence 
$h'=\varphi\circ h$ for the isomorphism
$\varphi=\hgt_{T'}^{-1}\circ\hgt_{T}\colon T\to T'$.
\end{proof}

\subsection*{Trees.}

A {\it tree}\/ on a set $U$
is a non-empty final segment $T$ of $(U^{\ast},\sqsupseteq)$.
(Recall that $a\sqsupseteq b \Longleftrightarrow$ $b$ is a truncation of $a$.)
The empty sequence $\varepsilon$ is the largest element of a tree
$T$ on $U$,  called the {\it root}\/ of $T$.
The elements of $T$ are called the {\it nodes}\/ of the tree $T$, and
the minimal elements of $T$ are called the {\it leafs}\/ of $T$.
Given a node $a=(a_1,\dots,a_n)$ of $T$, we denote by $T_a$ the tree 
$$T_a := \bigl\{ b\in U^*:a  b\in T\bigr\},$$ called the
{\it subtree}\/ of $T$ with root at $a$. 
Let $S$ be a tree on $U$ and $T$ a tree on $V$.
A map $\varphi\colon S\to T$  is called
{\it length-preserving}\/ 
if $\length(\varphi(s))=\length(s)$ for all $s\in S$.
Any increasing length-preserving map $S\to T$ is strictly increasing,
and the image of $S$ is a tree on $V$. Moreover:

\begin{lemma}
For a map $\varphi\colon S\to T$, the following are equivalent:
\begin{enumerate}
\item $\varphi$ is increasing and length-preserving.
\item $\varphi(\varepsilon)=\varepsilon$, and
for all $s\in S$, $a\in U$ there exists $b\in V$ with
$\varphi(s a) = \varphi(s)  b$.
\end{enumerate}
\end{lemma}
%\begin{proof}
%Let $\varphi$ be a tree embedding. 
%Since $\varepsilon \in \varphi(S)$ (because $\varphi(S)$ is a tree)
%and $\varphi$ is strictly increasing, we get $\varphi(\varepsilon)=\varepsilon$.
%Let $s\in S$, $a\in U$. Then $s \propini s\wu
%Conversely, suppose that $\varphi$ satisfies (2). Then clearly $\varphi$ is
%strictly increasing. In order to show that $\varphi(S)$ is a tree on $V$,
%let $s\in S$, $t\in V^{\ast}$
%with $t\ini\varphi(s)$; we have to show $t\in\varphi(S)$. 
%We may assume $s\neq\varepsilon$. Moreover, we may assume that $s$ has minimal
%length with $t\ini\varphi(s)$.
%Write $s=s' a$ for some $a\in U$ and some $s'\in S$.
%Then $\varphi(s)=\varphi(s') b$ for some $b\in V$; so if
%$t\neq\varphi(s)$, then $t\ini\varphi(s')$, contradicting the minimality of
%$s$. Hence $t=\varphi(s)\in\varphi(T)$.
%\end{proof}
\noindent
Given a map $\varphi\colon
U\to V$ we obtain an increasing length-preserving map $U^{\ast}\to V^{\ast}$,
also denoted by $\varphi$, by setting $$\varphi(a_1,\dots,a_n):=
\bigl(\varphi(a_1),\dots,\varphi(a_n)\bigr)\quad\text{for 
$a_1,\dots,a_n\in U$.}$$

\subsection*{Rank of a well-founded tree.}

Let $T$ be a well-founded tree on a set $U$. 
Then $\hgt(T)=\hgt(\varepsilon)+1$. (Recall that $\varepsilon$ is the root of $T$). 
We call the ordinal
$\hgt(\varepsilon)$ the {\it rank}\/ of the tree $T$, denoted by
$\rk(T)$. Hence $\rk(T)=0$ if and only if $T=\{\varepsilon\}$,
and $\hgt(a)=\rk(T_a)$ for any $a\in T$.
Also note that
\begin{equation}\label{Thgt}
\rk(T) = \sup\bigl\{ \hgt\bigl((x)\bigr)+1 : x\in U, (x)\in T\bigr\}.
\end{equation}
Every tree $T'$ on $U$ with $T'\subseteq T$ is well-founded with 
$\rk(T')\leq\rk(T)$.
More generally:

\begin{lemma}\label{Kechris-Lemma}
Let $S$, $T$ be trees on $U$ and $V$, respectively. 
If $T$ is well-founded, then the following are equivalent:
\begin{enumerate}
\item $S$ is well-founded with $\rk(S)\leq\rk(T)$.
\item There exists a length-preserving increasing map $S\to T$.
\item There exists a strictly increasing map $S\to T$.
\end{enumerate}
\end{lemma}
\begin{proof}
For a proof of ``(1)~$\Rightarrow$~(2)'' 
see, e.g., \cite{Kechris}, (2.9). 
The implication ``(2)~$\Rightarrow$~(3)'' is trivial, and
``(3)~$\Rightarrow$~(1)'' 
follows from part~(1) of Lemma~\ref{hgt-properties}.
\end{proof}

\subsection*{Invariants of Noetherian ordered sets.}
Every Noetherian ordered set $S$ is well-founded, hence has a certain
height $\hgt(S)$.
Following \cite{Kriz-Thomas}, 
we now introduce two other ordinal-valued invariants
associated to every Noetherian ordered set $S$, 
called the {\it type}\/ and the {\it width}\/ of $S$. Together, they 
measure the complexity of $S$.
First, for an ordered set $S$ 
we define the following trees on $S$:
\begin{align*}
\Dec(S) &:= \bigl\{(s_1,\dots,s_n)\in S^{\ast}:s_i>s_j
\text{ for all $1\leq i<j\leq n$}\bigr\}, \\
\Ant(S) &:= \bigl\{(s_1,\dots,s_n)\in S^{\ast}:s_i\,\,||\,\, s_j
\text{ for all $1\leq i<j\leq n$}\bigr\},\\
\Bad(S) &:= \bigl\{(s_1,\dots,s_n)\in S^{\ast}: s_i\not\leq s_j
\text{ for all $1\leq i<j\leq n$}\bigr\}.
\end{align*}
We call $\Dec(S)$ the {\it tree of decreasing sequences of $S$,}
$\Ant(S)$ the {\it tree of anti\-chains of $S$,} and
$\Bad(S)$ the {\it tree of bad sequences of $S$.}\/
Note that $S$ is well-founded if and only if $\Dec(S)$ is well-founded,
and $\hgt(S)=\rk\bigl(\Dec(S)\bigr)$. The tree $\Ant(S)$ is 
well-founded if and only if
every antichain of $S$ is finite, 
and $S$ is Noetherian if and only if $\Bad(S)$ is well-founded.
For any quasi-embedding 
$\varphi\colon S\to T$ of Noetherian ordered sets we have 
$\varphi\bigl(\Ant(S)\bigr)
\subseteq\Ant(T)$ and $\varphi\bigl(\Bad(S)\bigr)
\subseteq\Bad(T)$.

\begin{definition}\label{Def-wdt-ot}
Let $S$ be a Noetherian ordered set. Then
\begin{enumerate}
\item $\wdt(S):= \rk\bigl(\Ant(S)\bigr)$ is called the {\it width}\/ of $S$, and
\item $\ot(S) := \rk\bigl(\Bad(S)\bigr)$ is called the {\it type}\/ of $S$.
\end{enumerate}
\end{definition}
\noindent
If the ordering on $S$ is total (i.e., $S$ is well-ordered) then
$\Dec(S)=\Bad(S)$, hence $\rk\bigl(\Bad(S)\bigr)$ is the order type 
of $S$, justifying our choice of notation. Note that
\begin{equation}\label{wdhgt}
\hgt_{\Ant(S)}(s) = \wdt(S^{||s}) \qquad\text{for $s\in\Ant(S)$},
\end{equation}
where $S^{||s}:=\bigl\{y\in S: \text{$y||s_i$ for all $i$}\bigr\}$
for $s=(s_1,\dots,s_n)\in S^{\ast}$, and
\begin{equation}\label{ohgt}
\hgt_{\Bad(S)}(s) = \ot(S^{\not\geq s}) \qquad\text{for $s\in\Bad(S)$},
\end{equation}
where $S^{\not\geq s}:=\bigl\{y\in S: \text{$y\not\geq s_i$ for all $i$}\bigr\}$
for $s=(s_1,\dots,s_n)\in S^{\ast}$.

\subsection*{Characterization of height and type.}
The height and the type of a Noetherian ordered set allow important
reinterpretations.
Recall that any total ordering extending a Noetherian ordering
is a well-ordering,
cf.~Proposition~\ref{Folk}.

\begin{theorem}\label{dJP}
Let $S$ be a Noetherian ordered set.
\begin{enumerate}
\item There exists a total ordering on $S$ extending the given ordering
of maximal possible order type; this order type equals $\ot(S)$.
\rom{(De Jongh-Parikh \cite{dJP}.)}
\item There exists a chain of $S$ of maximal possible order type; this
order type equals $\hgt(S)$.
\rom{(Wolk \cite{Wolk}.)}
\end{enumerate}
\end{theorem}
\noindent
Because of 
(1), the type of a Noetherian ordered set $S$ is sometimes also called
the {\it maximal order type}\/ of $S$.
The width of $S$ is finite if and only if 
there is some $n$ such that 
every antichain in $S$ has size $\leq n$; the smallest such $n$ is $\wdt(S)$. 
In this case, $S$ is a union of
$\wdt(S)$ many chains (Dilworth \cite{Dilworth}).
In general (i.e., for infinite width),
no characterization of the width similar to (1) or (2) in the theorem 
above seems to be known. (See \cite{Kriz-Thomas}, p.~77.) We refer to
\cite{Abraham}, \cite{Perles} for
discussions of Dilworth's Theorem in the case of infinite width.)

\subsection*{Basic facts about type and width.}
We record some basic properties:

\begin{prop}\label{oProp}
Let $S$ and $T$ be Noetherian ordered sets.
\begin{enumerate}
\item %If $\ot(S)$ is a limit ordinal, then 
$\ot(S)=\sup\bigl\{\ot(S^{\not\geq x})+1 : x\in S \bigr\}$.
%where $$S^{\not\geq x} := S\setminus (x) = \{y\in S : x\not\leq y\}\qquad
%\text{for $x\in S$.}$$
\item If there exists a quasi-embedding $S\to T$, then $\ot(S)\leq \ot(T)$. 
\rom{(}In particular if $S\subseteq T$, then $\ot(S)\leq\ot(T)$.\rom{)}
\item If there exists an increasing surjection $S\to T$, then
$\ot(S)\geq\ot(T)$.
\item $\ot(S\amalg T)=\ot(S) \oplus \ot(T)$ and
$\ot(S\times T)=\ot(S)\otimes \ot(T)$.
\end{enumerate}
\end{prop}
\begin{proof}
Part (1) follows from the identities \eqref{Thgt} and \eqref{ohgt},
and part (2) from Lemma~\ref{hgt-properties},~(1) and  
the remarks preceding Definition~\ref{Def-wdt-ot}.
For (3) suppose that $\varphi\colon S\to T$
is an increasing surjection. For every $t\in T$ choose $\psi(t)\in S$
with $\varphi\bigl(\psi(t)\bigr)=t$. If $(t_1,\dots,t_n)\in\Bad(T)$ then
$\bigl(\psi(t_1),\dots,\psi(t_n)\bigr)\in\Bad(S)$, so $\psi$ induces an
increasing and length-preserving map $\Bad(T)\to\Bad(S)$. Hence
$\ot(T)=\rk\bigl(\Bad(T)\bigr) \leq \rk\bigl(\Bad(S)\bigr)=\ot(S)$ by
Lemma~\ref{hgt-properties},~(1).
For a proof of (4) see \cite{dJP} or \cite{Kriz-Thomas}.
\end{proof}
\noindent
The computation of $\ot(S)$, for a concrete given Noetherian
ordered set $S$, is often quite hard; see, e.g., \cite{Schmidt}.
In Section~\ref{OrderingsOfMonomialIdeals} we will bound the maximal order
type of the Noetherian ordered set of monomial ideals.

\begin{prop}\label{wdprop}
Let $S$ and $T$ be Noetherian ordered sets.
\begin{enumerate}
\item %If $\ot(S)$ is a limit ordinal, then 
$\wdt(S)=\sup\bigl\{\wdt(S^{||x})+1 : x\in S \bigr\}$.
%where $$S^{||x} :=  \{y\in S : x || y\}\qquad\text{for $x\in S$.}$$
\item If there exists a quasi-embedding $S\to T$, then $\wdt(S)\leq \wdt(T)$. 
\rom{(}In particular if $S\subseteq T$, then $\wdt(S)\leq\wdt(T)$.\rom{)}
\item If there exists an increasing surjection $S\to T$, then
$\wdt(S)\geq\wdt(T)$.
\item $\wdt(S\amalg T)=\wdt(S)\oplus\wdt(T)$.
\end{enumerate}
\end{prop}
\begin{proof}
Part~(1) follows from \eqref{Thgt} and \eqref{wdhgt}, and part~(2) 
again from  Lemma~\ref{hgt-properties},~(1). 
The proof of (3) is similar to the proof of Proposition~\ref{oProp},~(3).
Part~(4) is
shown by induction on $\wdt(S)\oplus\wdt(T)$: Since
$(S\amalg T)^{||x}=S^{||x}\amalg T$ if $x\in S$ and
$(S\amalg T)^{||x}=S\amalg T^{||x}$ if $x\in T$, we get from part (1)
$$\wdt(S\amalg T)=\sup\bigl\{\wdt\bigl(S^{||x}\amalg T\bigr)+1,
\wdt\bigl(S\amalg T^{||y}\bigr)+1 : x\in S, y\in T\bigr\},$$
and so by inductive hypothesis and \eqref{NatSum}
we obtain $\wdt(S\amalg T)=\wdt(S)\oplus\wdt(T)$ as desired.
\end{proof}
\noindent
A formula for the width of $\alpha\times\beta$ ordered component-wise, 
where $\alpha$ and $\beta$ are
ordinals, can be found in \cite{Abraham}, and a formula for the width of
$S\times T$ ordered lexicographically, for Noetherian ordered sets $S$ and $T$,
in \cite{Abraham-Bonnet}.

%Note that by (2) and (3) in Proposition~\ref{oProp}, if $S$ is a Noetherian
%ordered set and $S_0,S_1$ are disjoint subsets of $S$ with $S=S_0\cup S_1$,
%then $\ot(S)\leq\ot(S_0)\oplus\ot(S_1)$. More generally:
%
%\begin{lemma}
%Suppose that $S$ is a Noetherian ordered set which can be written as the 
%disjoint union
%$$S=\coprod_{\alpha<\lambda} S_\alpha$$
%of a family $\{S_\alpha\}_{\alpha<\lambda}$ of subsets of $S$
%\rom{(}indexed by all ordinals less than some ordinal $\lambda$\rom{)}. Then
%$$\ot(S) \leq \bigoplus_{\alpha<\lambda} \ot(S_\alpha),$$
%with equality if any two elements $s\in S_\alpha$ and $t\in S_\beta$ with 
%$\alpha<\beta<\lambda$ are incomparable.
%\end{lemma}
%\begin{proof}
%By transfinite induction on $\lambda$. The case where $\lambda$ is a
%successor ordinal is covered by the remarks preceding the lemma.
%The limit stage follows using \dots
%\end{proof}
%\noindent
%By Proposition~\ref{oProp},~(3), this implies that $S$ is a Noetherian
%ordered set which can be written as a disjoint union
%$S=\coprod_{\alpha<\lambda} S_\alpha$
%of a family $\{S_\alpha\}_{\alpha<\lambda}$ of subsets of $S$
%(not necessarily pairwise incomparable), then
%$o(S) \leq \bigoplus_{\alpha<\lambda} o(S_\alpha)$.

\subsection*{Connections between the invariants.}
Height, width and type are related by the following fundamental inequality:

\begin{prop}\rom{(Height-Width Theorem, \cite{Kriz-Thomas}.)}\label{Height-Width}
Let $S$ be a Noetherian ordered set. Then
$$\ot(S) \leq \hgt(S)\otimes\wdt(S).$$
\end{prop}
\noindent
This generalizes the well-known fact that a finite ordered set with
at least $rs+1$ elements contains a chain with $r+1$ elements or an
antichain with $s+1$ elements.
\begin{proof}
Since we will need a similar idea below (Lemma~\ref{Ant-Prime}), 
we sketch the proof of 
Proposition~\ref{Height-Width}. Let $g=\hgt_{\Ant(S)}$ be the
height function of the tree of antichains of $S$, and define
$h\colon\Bad(S)\setminus\{\varepsilon\}\to\wdt(S)$ by 
\begin{multline*}
h(s_1,\dots,s_n) := \\\min\bigl\{ g(s_{i_1},\dots,s_{i_m}):
1\leq i_1<\cdots<i_m=n, \hgt(s_{i_1})\leq\cdots\leq\hgt(s_{i_m})\bigr\}
\end{multline*}
for $(s_1,\dots,s_n)\in\Bad(S)$, $n\geq 1$.
It is easy to see that $f(\varepsilon):=\varepsilon$ and
\begin{multline*}
f(s_1,\dots,s_n) := \\\bigl(\bigl(\hgt(s_1),h(s_1)\bigr),
\bigl(\hgt(s_2),h(s_1,s_2)\bigr),\dots,
\bigl(\hgt(s_n),h(s_1,\dots,s_n)\bigr)\bigr)
\end{multline*}
defines a strictly increasing map $f\colon\Bad(S)\to
\Bad\bigl(\hgt(S)\times\wdt(S)\bigr)$. Hence
$$\ot(S)\leq
\ot\bigl(\hgt(S)\times\wdt(S)\bigr)=\hgt(S)\otimes\wdt(S)$$ by Lemma
\ref{Kechris-Lemma} and Proposition~\ref{oProp}~(4).
\end{proof}
\noindent
The following proposition connects the type of a Noetherian ordered set
$S$ with the height of the well-founded
ordered set $\bigl(\cF(S),{\supseteq}\bigr)$ of its final segments:

\begin{prop}\rom{(Bonnet-Pouzet \cite{Bonnet-Pouzet}.)}\label{Bonnet-Pouzet}
For every Noetherian ordered set $S$,
$$\hgt\bigl(\cF(S)\bigr)=\ot(S)+1.$$
\end{prop}

\noindent
Let us outline the main idea of the proof of this fact. 
First we note that Proposition~\ref{Bonnet-Pouzet}
is a consequence of the following lemma,
the characterization \eqref{height-sup} of the height,
Theorem~\ref{dJP}~(1), and the
fact that $\hgt\bigl(\cF(\alpha)\bigr)=\alpha+1$ for any ordinal $\alpha$.

\begin{lemma}
Let $(S,\leq)$ be an ordered set.
The assignment ${\leq'} \mapsto \cF(S,{\leq'})$ is a one-to-one
correspondence between the total orderings $\leq'$ on $S$ extending $\leq$
and the maximal chains of the ordered set $\bigl(\cF(S,\leq),\supseteq\bigr)$.
\end{lemma}

\noindent
It is easy to verify that the map 
in the lemma is well-defined, and it is clearly one-to-one.
Now let $\cC$ be a maximal chain of $\cF(S,\leq)$. Define a binary relation
$\leq_{\cC}$ on $S$ by $x\leq_{\cC} y :\Longleftrightarrow$ every $F\in\cC$
which contains $x$ also contains $y$. 
The main work consists in establishing that for
any two distinct elements $x\neq y$
of $S$ for
which there exists $F\in\cF(S,\leq)$ with $x\in F$, $y\notin F$,
there exists $G\in\cC$ with $x\in G$, $y\notin G$.
(For the details see \cite{Bonnet-Pouzet}.) 
From this it is straightforward to check that
$\leq_{\cC}$ is a total ordering on $S$ extending $\leq$, and
$\cF(S,\leq_{\cC})=\cC$.

\subsection*{Height and minimal order type.}
Let $(S,\leq)$ be a Noetherian ordered set. 
Then there is a smallest ordinal $\alpha$ such that $\leq$ has an
extension to a well-ordering on $S$ of order type $\alpha$. We call
$\alpha$ the {\it minimal order type}\/ of $S$, denoted by
$\ot^*(S)$. We show here that this ordinal agrees with the height of
$S$ if $\hgt(S)$ is a limit ordinal, and differs from the height of
$S$ at most by a finite ordinal otherwise. (This was observed in
\cite{Schmidt}, p. 8--10.) In the first case 
we also show how to obtain an
extension $\leq^*$ of $\leq$ to a well-ordering of $S$ of order
type $\ot^*(S)$.
This will all be based on the following observation:

\begin{lemma} 
The height function $\hgt\colon S\to\hgt(S)$ has finite fibers.
\end{lemma}
\begin{proof}
By Lemma~\ref{Folk-Lemma},~(1) 
and the fact that $\hgt$ is strictly increasing.
\end{proof}
\noindent
Let now $\leq'$ be any extension of $\leq$ to a well-ordering on $S$.
Define a binary relation $\leq^*$ on $S$ by
\begin{equation}\label{MinimalOT}
x\leq^* y \qquad :\Longleftrightarrow\qquad \bigl(\hgt(x),x\bigr) 
\leq_{\text{lex}} \bigl(\hgt(y),y\bigr).
\end{equation}
Here $\leq_{\text{lex}}$ denotes the lexicographic product of
the ordering of $\hgt(S)$ and $\leq'$, that is,
$$(\alpha,x) \leq_{\text{lex}} (\beta,y) \qquad\Longleftrightarrow\qquad \alpha<\beta \text{ or }
(\alpha=\beta\text{ and } x\leq' y),$$
for all $\alpha,\beta<\hgt(S)$ and $x,y\in S$. 
It is straightforward to check that
$\leq^*$ is an extension of $\leq$ to a well-ordering of $S$.
We denote the height function of 
$(S,\leq^*)$ by $\hgt^*\colon S\to\hgt^*(S)$.

\begin{lemma}
For all $x\in S$, $\hgt^*(x) < \hgt(x)+\omega$.
\end{lemma}
\begin{proof}
By transfinite induction on $\alpha=\hgt(x)$. If $\alpha=0$, then $x$ is
one of the finitely many minimal elements of $(S,\leq)$. Hence there
are only finitely many $y\in S$ with $y\leq^* x$, hence $\hgt^*(x)<\omega$.
For the successor case, suppose that $\hgt(x)=\alpha+1$, and choose
$y<x$ with $\hgt(y)=\alpha$. There are only finitely many elements
$y=y_0<^* y_1<^*\cdots<^*y_m=x$ of $S$ which lie between $y$ and $x$ in
the ordering $\leq^*$. So
$\hgt^*(x)=\hgt^*(y)+m$, and by induction we get $\hgt^*(y)<\hgt(y)+\omega=
\alpha+\omega$. Hence $\hgt^*(x)<(\alpha+1)+\omega=\hgt(x)+\omega$.
Now suppose that $\alpha$ is a limit ordinal.
Let $x_0<^*x_1<^*\cdots<^*x_n=x$ be the elements of height $\alpha$
which are $\leq^*x$; so $\hgt^*(x)=\hgt^*(x_0)+n$.
We have $\hgt(y)<\hgt(x)$ for all $y\in S$ with $y<^* x_0$, hence
$\hgt^*(y)<\hgt(y)+\omega\leq\hgt(x)$ by inductive hypothesis and since
$\hgt(x)=\alpha$ is a limit ordinal. Therefore
$$\hgt^*(x_0)=\sup\bigl\{\hgt^*(y)+1:y<^*x_0\bigr\}\leq\hgt(x),$$
and hence $\hgt^*(x) < \hgt(x)+\omega$ as required.
\end{proof}

\begin{cor}\label{oStarCor}\

\begin{enumerate}
\item If $\hgt(S,\leq)$ is a successor ordinal, then
$$\hgt(S,\leq) \leq \ot^*(S,{\leq})\leq\ot(S,\leq^*) < \hgt(S,\leq) + \omega.$$
\item If $\hgt(S,\leq)$ is a limit ordinal, then 
$$\ot^*(S,{\leq})=\ot(S,{\leq^*})=\hgt(S,\leq).$$ \qed
\end{enumerate}
\end{cor}
\noindent
In Section~\ref{OrderingsOfMonomialIdeals} we will apply this corollary in
the following situation: Suppose that $S$ is a Noetherian ordered set
with a largest element $s_0$ whose height $\hgt(s_0)$ is a limit ordinal.
Then the Noetherian ordered set $S_0=S\setminus\{s_0\}$ has height
$\hgt(s_0)$. By part (1) of the corollary, it follows that
$\ot^*(S_0)=\hgt(s_0)$ and hence $S$ has minimal order type 
$\ot^*(S)=\hgt(s_0)+1=\hgt(S)$.

\subsection*{Total orderings of monomials.}

As an illustration for the material in this section, we now compute
the invariants $\hgt$, $\ot$, $\wdt$, and $\hgt^*$ for the Noetherian
ordered set $\Nn^m$, and hence for the set of monomials in the
polynomial ring $K[X_1,\dots,X_m]$ over a field $K$, ordered
by divisibility (see Example~\ref{example-Nn}).
It is convenient to consider, slightly more generally, 
Noetherian ordered sets of the form
$\Nn^m\times S$, where $S$ is a finite non-empty ordered set.

\begin{lemma}
Let $S$ be a finite non-empty ordered set and $m>0$. Then
$$\hgt(\Nn^m\times S)=\omega,\ \ot(\Nn^m\times S)=\omega^m\abs{S},\  
\hgt^*(\Nn^m\times S)=\omega.$$
\end{lemma}
\begin{proof}
The function $(\nu,s)\mapsto\abs{\nu}+\hgt_S(s)\colon \Nn^m\times S\to\Nn$ is
strictly increasing. Hence $\hgt\bigl((\nu,s)\bigr)=\abs{\nu}+\hgt_S(s)$ 
for all $(\nu,s)\in\Nn^m\times S$, and
$\hgt(\Nn^m\times S)=\omega$. 
By Corollary~\ref{oStarCor},~(2) this yields $\hgt^*(\Nn^m\times S)=\omega$.
By Proposition~\ref{oProp},~(4) we get $\ot(\Nn^m\times S)=\ot(\Nn^m)\otimes
\ot(S)=\omega^m\abs{S}$.
\end{proof}

\noindent
The lexicographic ordering $\leq_{\text{lex}}$ 
on $\Nn^m$ is an example for a
total ordering of $\Nn^m$ extending the product ordering $\leq$ and having 
maximal order type. Given any total ordering $\leq'$ on $\Nn^m$
extending $\leq$, we obtain a total ordering $\leq^*$ on $\Nn^m$
of minimal order type $\omega$ extending $\leq$,
as shown in \eqref{MinimalOT}:
$$\nu \leq^* \mu \quad :\Longleftrightarrow\quad
\text{$\abs{\nu}<\abs{\mu}$, or $\abs{\nu}=\abs{\mu}$ and $\nu\leq'\mu$.}$$
For $\leq'=\leq_{\text{lex}}$ (the lexicographic ordering of $\Nn^m$) the
ordering obtained in this way is commonly called
the {\it degree-lexicographic ordering}\/ of $\Nn^m$.
Orderings of the form $\leq^*$ are called {\em degree-compatible.} 
In applications, one is usually interested in total orderings $\leq'$
extending $\leq$ which  are {\it semigroup orderings,}\/ that is,
which satisfy the condition
$$\nu \leq' \mu \Rightarrow \nu+\lambda \leq' \mu+\lambda
\qquad\text{for all $\nu,\mu,\lambda\in\Nn^m$.}$$
A total semigroup ordering on $\Nn^m$ extending $\leq$ is called a
{\it term ordering.}\/ 
Via the usual identification of $\Nn^m$ with $X^\diamond$,
term orderings on $\Nn^m$ hence correspond to term orderings on
$X^\diamond$ (as defined in Remark~\ref{Sdiamond})
where $X=\{X_1,\dots,X_m\}$ carries the trivial ordering.
The lexicographic and 
degree-lexicographic orderings of $\Nn^m$ are term-orderings.
A complete description of all term orderings on $\Nn^m$
is available (see \cite{Robbiano} or \cite{Weispfenning-Orders}): 
For any such ordering
$\leq'$ there exists an invertible $m\times m$-matrix $A$ with real
coefficients such that 
\begin{equation}\label{Termorder}
\nu\leq'\mu \quad \Longleftrightarrow\quad A\nu \leq_{\text{lex}}
A\mu \qquad\text{for all $s,t\in S$, $\nu,\mu\in\Nn^m$,}
\end{equation} 
where $\leq_{\text{lex}}$ denotes the lexicographic ordering on 
$\Rr^m$. Conversely, any matrix $A\in\operatorname{GL}(m,\Rr)$ 
satisfying $Ae_i\geq_{\text{lex}} 0$ for all $i=1,\dots,m$ 
(where $e_i$ denotes the $i$-th unit vector in $\Rr^m$) gives
rise to a term ordering $\leq'$ on $\Nn^m$, via \eqref{Termorder}.
In particular, the order types of term orders on $\Nn^m$ are
the ordinals of the form $\omega^k$ with $1\leq k\leq m$.
There are only $m!$ many different term
orderings of maximal order type $\omega^m$ on $\Nn^m$ (obtained by choosing
permutation matrices for $A$), and for each $1\leq k<m$ 
there are continuum many
term orderings on $\Nn^m$ with order type $\omega^k$. 
(See \cite{Martin-Scott}.)

More generally, a 
{\em ranking} of $\Nn^m\times S$, where $S$ is a finite non-empty
set, is a total ordering of $\Nn^m\times S$ which extends the product ordering
on $\Nn^m\times S$ (where $S$ is equipped with the trivial ordering) and
satisfies
$$(\nu,s)\leq (\mu,t) \Rightarrow (\nu+\lambda,s)\leq (\mu+\lambda,t)
\qquad\text{for all $\nu,\mu,\lambda\in\Nn^m$.}$$
Rankings play a role in algorithmic differential algebra (e.g.,
in the theory of Riquier-Janet bases) similar to the role played by 
term orderings in ordinary 
algorithmic algebra (in the theory of Gr\"ob\-ner bases), see \cite{Kolchin},
\cite{Rust}: the elements $(\nu,i)$ of 
$\Nn^m\times\{1,\dots,n\}$ correspond to the derivatives
$\partial^{|\nu|}y^i/\partial X_1^{\nu_1}\cdots\partial X_m^{\nu_m}$, where
$y^1,\dots,y^n$ are differential indeterminates over a differential ring
with $m$ commuting derivations $\partial/\partial X_1,\dots,
\partial/\partial X_m$.
Rankings also naturally arise when Gr\"ob\-ner basis theory is generalized to
finitely generated free modules over $K[X]$, see \cite{Eisenbud}, Chapter~15.
We refer to 
\cite{Rust} for a (rather involved) classification of rankings which extends 
the one of term orderings described above. It would be interesting to
determine the possible order types of rankings from this classification.

We now turn to the width of $\Nn^m\times S$.
By the Height-Width 
Theorem~\ref{Height-Width} we have $\wdt(\Nn^m\times S)\geq\omega^{m-1}
\abs{S}$, since
$\hgt(\Nn^m\times S)=\omega$ and $\ot(\Nn^m\times S)=\omega^{m}\abs{S}$. 
We will show:

\begin{prop}\label{wd-Nnm}
$\wdt(\Nn^m\times S)=\omega^{m-1}\abs{S}$, for all $m>0$ and all
finite non-empty ordered sets $S$.
\end{prop}
\noindent
In the proof, we will use the following lemma.

\begin{lemma}
Let $T_1,\dots,T_m$ be well-ordered sets, $T=T_1\times\cdots\times T_m$. 
Then, for any $a=(a_1,\dots,a_m)\in T$:
$$\wdt(T^{||a}) = \bigoplus_{\varepsilon} \wdt(T_1^{\varepsilon_1 a_1}\times
\cdots\times T_m^{\varepsilon_n a_m}),$$
where the sum runs over all $\varepsilon=(\varepsilon_1,\dots,\varepsilon_m)\in
\{ {\leq}, {>} \}^m$ such that for some $i,j$, we have
$\varepsilon_i={\leq}$ and $\varepsilon_j={>}$.
\end{lemma}
\begin{proof}
By \eqref{wdhgt},  the fact that
$$T^{||a} = \coprod_{\varepsilon} T_1^{\varepsilon_1 a_1}\times
\cdots\times T_m^{\varepsilon_m a_m},$$ and
Proposition~\ref{wdprop},~(4).
\end{proof}
\noindent
In order to prove Proposition~\ref{wd-Nnm}, it suffices to show
$\wdt(\Nn^m\times S)\leq \omega^{m-1}\abs{S}$, for all $m>0$ and
finite ordered sets $S\neq\emptyset$. We proceed by induction
on $m$. Note first that if $M\neq\emptyset$ is an ordered set, then
we have a natural quasi-embedding of $M\times S$ into $M\amalg \cdots\amalg M$
($|S|$ many times); hence if $M$ is Noetherian, then $\wdt(M\times S)\leq
\wdt(M)\abs{S}$. Taking $M=\Nn$ this yields $\wdt(\Nn\times S)\leq\abs{S}$,
and hence the base case $m=1$ of our induction. Now suppose $m>1$. 
With $T=\Nn^m$, we have
$$\wdt(T) = \sup\bigl\{ \wdt(T^{||a})+1 : a\in T\bigr\}.$$
Let $\varepsilon=(\varepsilon_1,\dots,\varepsilon_m)\in
\{ {\leq}, {>} \}^m$ be such that for some $i,j$, we have
$\varepsilon_i={\leq}$ and $\varepsilon_j={>}$, 
and let $a=(a_1,\dots,a_m)\in\Nn^m$. For $b\in\Nn$, 
if $\varepsilon={\leq}$, then $\Nn^{\varepsilon b}$ is finite, and
if $\varepsilon={>}$, then $\Nn^{\varepsilon b}\cong\Nn$.
It follows that $\Nn^{\varepsilon_1 a_1}\times
\cdots\times \Nn^{\varepsilon_n a_m} \cong \Nn^k\times U$ for some $1\leq k<m$
and some non-empty finite ordered set $U$. By induction
$\wdt(\Nn^{\varepsilon_1 a_1}\times
\cdots\times \Nn^{\varepsilon_n a_m})\leq \omega^{k-1}\abs{U}$. 
By the lemma, this yields
$\wdt\bigl((\Nn^m)^{||a}\bigr)<\omega^{m-1}$ and thus 
$\wdt(\Nn^m\times S)\leq\wdt(\Nn^m)\abs{S}\leq\omega^{m-1}\abs{S}$ as desired. \qed

\subsection*{A variant of the tree of antichains.}
Let $S$ be a Noetherian ordered set.
In Section~\ref{OrderingsOfMonomialIdeals}, 
we will need a variant of the tree of antichains $\Ant(S)$ of $S$.
For this, we fix a total ordering $\leq'$ extending the ordering $\leq$ of $S$.
We define a tree 
$$\Ant_{\leq'}(S) := \bigl\{(s_1,\dots,s_n)\in S^{\ast}:s_i\,\,||\,\, s_j
\text{ and } s_i<' s_j
\text{ for all $1\leq i<j\leq n$}\bigr\}$$
on $S$, an ordered subset of the tree $\Ant(S)$ of antichains of $S$.
We clearly
have $\rk(\Ant_{\leq'}(S)\bigr)\leq\rk\bigl(\Ant(S)\bigr)=\wdt(S)$, and
we conjecture that in general, the reverse inequality is also true. Here, we 
confine ourselves to showing:

\begin{lemma}\label{Ant-Prime}
$\rk\bigl(\Ant_{\leq'}(\Nn^m\times S)\bigr)=\omega^{m-1}|S|$
for any finite non-empty ordered set $S$ and total ordering
$\leq'$ on $\Nn^m\times S$ of order type $\omega$ extending the product
ordering.
\end{lemma}
\noindent
This follows immediately from 
Proposition~\ref{wd-Nnm} above and the following fact:

\begin{lemma}
$\ot(S) \leq \ot(S,\leq')\otimes\rk\bigl(\Ant_{\leq'}(S)\bigr)$,
for any Noetherian ordered set $(S,\leq)$ and any total ordering $\leq'$
extending $\leq$.
\end{lemma}
\begin{proof}
Put $\alpha=\rk\bigl(\Ant_{\leq'}(S)\bigr)$.
Let $g=\hgt_{\Ant_{\leq'}}(S)$ 
be the height function of $\Ant_{\leq'}(S)$. Define
$h\colon\Bad(S)\setminus\{\varepsilon\}\to\alpha$ by
\begin{multline*}
h(s_1,\dots,s_n) := \\\min\bigl\{ g(s_{i_1},\dots,s_{i_m}):
1\leq i_1<\cdots<i_m=n, s_{i_1} <' \cdots <' s_{i_m}\bigr\},
\end{multline*}
for $(s_1,\dots,s_n)\in\Bad(S)$, $n\geq 1$.
Then $f(\varepsilon):=\varepsilon$ and
$$
f(s_1,\dots,s_n) := \bigl(\bigl(s_1,h(s_1)\bigr),
\bigl(s_1,h(s_1,s_2)\bigr),\dots,
\bigl(s_n,h(s_1,\dots,s_n)\bigr)\bigr)
$$
defines a strictly increasing map $$f\colon\Bad(S)\to
\Bad\bigl((S,\leq')\times\alpha\bigr).$$ Hence
$\ot(S)\leq
\ot\bigl((S,\leq')\times\alpha\bigr)=
\ot(S,\leq')\otimes\alpha$.
\end{proof}

\section{The Ordered Set of Hilbert Polynomials}\label{HilbertSection}

\noindent
In this section we discuss the sets of Hilbert and Hilbert-Samuel
polynomials of finitely generated graded $K$-algebras (where $K$ is a field)
as ordered sets,
with the ordering given by the relation of eventual dominance.
Macaulay's Theorem on the possible Hilbert functions of such $K$-algebras
will play an important role. We begin by recalling this theorem and some
of its consequences, in particular a description of all Hilbert-Samuel
polynomials of finitely generated graded $K$-algebras. This description will
be used to give an interpretation of the height function on $\cF(\Nn^m)$ in
terms of the coefficients of the Hilbert-Samuel polynomial. We give two
applications concerning increasing chains of ideals in polynomial rings.

\subsection*{Integer-valued polynomials.}
Recall that a polynomial $f(T)\in\Qq[T]$ (in a single variable $T$) is called
an {\it integer-valued polynomial}\/ if $f(s)$ is an integer for all $s\in\Nn$. 
For example,
$$\binom{T}{j}=\frac{T(T-1)\cdots (T-j+1)}{j!}\in\Qq[T] \qquad
(j\in\Nn)$$
is integer-valued. The polynomials $\binom{T+i}{i}$ (for $i\in\Nn$)
form a basis for the $\Zz$-submodule of $\Qq[T]$ consisting of the
integer-valued polynomials. In other words, every non-zero integer-valued 
polynomial
$f(T)\in\Qq[T]$ can be uniquely written in the form
$$f(T)=b_d\binom{T+d}{d}+b_{d-1}\binom{T+d-1}{d-1}+\cdots+b_0\binom{T+0}{0}$$
with $b_0,\dots,b_d\in\Zz$, $b_d\neq 0$.
%By the {\em leading coefficient}\/ of $f(T)$ we
%shall mean the non-zero integer $b_d$. (It is $d!$ times the coefficient 
%of the $T^d$ term in $f(T)$.)
We totally order the integer-valued polynomials by {\em
dominance}\/: if $f(T)=\sum_{i=0}^d b_i\binom{T+i}{i}$ and
$g(T)=\sum_{j=0}^d c_j\binom{T+j}{j}$ with $b_i,c_j\in\Zz$, then
\begin{align*}
f(T) \preceq g(T) \quad :\Longleftrightarrow&\quad  f(s) \leq g(s) 
\text{ for all $s \gg 0$} \\
\Longleftrightarrow&\quad (b_d,b_{d-1},\dots,b_0) \leq_{\text{lex}} (c_d,c_{d-1},\dots,c_0) \text{ in $\Zz^{d+1}$.}
\end{align*} 
With this ordering, the ring of integer-valued polynomials becomes an
ordered integral domain.

\subsection*{Hilbert polynomials of homogeneous ideals.}

In this section, $K$ denotes a field.
Let $I$ be a homogeneous ideal of a polynomial ring
$K[X]=K[X_1,\dots,X_m]$ over $K$; that is, the ideal $I$ is generated
by homogeneous elements of positive degree.
Then as a $K$-vector space,
$R=K[X]/I$ has a direct-sum decomposition $R=R_0\oplus R_1\oplus\cdots$
given by
$$R_s := \bigl\{f+I\in R : \text{$f\in K[X]$ has total degree $s$}\bigr\}.$$
This decomposition makes $R$ into a {\it graded $K$-algebra}: we have $R_0=K$
and $R_s\cdot R_t\subseteq R_{s+t}$ for all $s,t$. Each component $R_s$ is
a finite-dimensional vector space over $K$. The function $H_{I}\colon\Nn\to\Nn$
defined by $H_{I}(s)=\dim_K R_s$ is called the {\it Hilbert function}\/ of $I$.
There exists an integer-valued polynomial 
$P_I$ of degree $< m$ (the {\em Hilbert polynomial}\/ of $I$)
such that $s\mapsto P_I(s)$ agrees with
$s\mapsto H_I(s)$  for sufficiently large $s$. The degree of $P_I$ is
one less than
the Krull dimension of the ring $R$. (See, e.g., \cite{Eisenbud}, 
Corollary~13.7.) If $I=(1)$ is the unit ideal in $K[X]$, 
we put $H_I(s)=0$ for all $s$ and $P_I=0$. As usual the degree of the
zero polynomial is $\deg 0:=-1$.

For a final segment $E$ of $\Nn^m$,
we call $H_E:=H_{I_E}$ and $P_E:=P_{I_E}$ the {\it Hilbert function}\/
and {\it Hilbert polynomial}\/ of $E$, respectively, 
where $I_{E}\subseteq\Qq[X]$ is the monomial ideal corresponding to $E$.
Given any final segment $E$ of $\Nn^m$, 
let us write $V_E:=\Nn^m\setminus E$ for the complement of $E$ in $\Nn^m$ (an initial 
segment of $\Nn^m$). We then have $H_E(n)=\abs{V_{E,n}}$ for all $n$, where
$Z_n := \bigl\{z\in Z : \abs{z}=n\bigr\}$
for $Z\subseteq\Nn^m$.

\subsection*{Macaulay's Theorem.}

A classical theorem of Macaulay characterizes exactly those functions
$f\colon\Nn\to\Nn$ which arise as Hilbert functions of homogeneous
ideals $I\subseteq K[X]$. 
Before we state Macaulay's Theorem, we have to introduce some more notation.
(As a general reference for this material, we recommend \cite{Bruns-Herzog},
Chapter~4.)
Given an integer $d\geq 1$, 
every positive integer $a$ can be written uniquely in the form
$$a=\binom{a_d}{d}+\binom{a_{d-1}}{d-1}+\cdots+\binom{a_1}{1}$$
where $a_d>a_{d-1}>\cdots >a_1\geq 0$. 
This sum is called the
{\it $d$-th Macaulay representation}\/ of $a$, 
and $(a_d,\dots,a_1)$ are called
the {\it $d$-th Macaulay coefficients}\/ of $a$.
We have $a\leq b$ if and only if $(a_d,\dots,a_1) 
\leq_{\text{lex}} (b_d,\dots,b_1)$. We define
$$a^{\langle d\rangle}=\binom{a_d+1}{d+1}+\binom{a_{d-1}+1}{(d-1)+1}+\cdots+
\binom{a_1+1}{1+1},$$
and $0^{\langle d\rangle}:=0$. 
We have
(for a proof see \cite{Bruns-Herzog}, p.~162):

\begin{theorem}\rom{(Macaulay, \cite{Macaulay})} \label{Macaulay}
Let $f\colon\Nn\to\Nn$. The following are equivalent:
\begin{enumerate}
\item There exists a homogeneous ideal $I\subseteq K[X]$ with $H_I=f$.
\item We have $f(1)=m$, and if $M_n$
denotes the set of the first
$f(n)$ elements of $\Nn^{m}$ of degree $n$, in the 
lexicographic ordering, then
$M=\bigcup_{n\in\Nn} M_n$ is an initial segment of $\Nn^{m}$.
\item $f(0)=1$, $f(1)=m$, 
and $f(n+1)\leq f(n)^{\langle n\rangle}$ for all $n\geq 1$.
\end{enumerate}
\end{theorem}
\noindent
A final segment $E$ of $\Nn^m$ is  called a
{\it lex-segment}\/ of $\Nn^m$ if for every $n$ the set
$E_n=\bigl\{e\in E:\abs{e}=n\}$ of elements of $E$ having degree $n$ is
a  final segment of $(\Nn^m)_n$ under the lexicographic ordering. 
(This terminology is used slightly differently
in \cite{Bruns-Herzog}.) 
If $M$ is as in statement (2) of the theorem, then clearly $\Nn^m\setminus M$ 
is a lex-segment of $\Nn^m$
with Hilbert function $f$.

The zero ideal of $K[X]=K[X_1,\dots,X_m]$ has Hilbert polynomial $\binom{T+m-1}{m-1}$.
The following characterization of Hilbert polynomials of non-zero ideals
is well-known:

\begin{cor}\label{Macaulay-Corollary}
A polynomial $P(T)\in\Qq[T]$ is a Hilbert polynomial of
some non-zero homogeneous  ideal of $K[X]=K[X_1,\dots,X_m]$ if and only if
\begin{equation}\label{Hilbert-Poly}
P(T) = \binom{T+a_1}{a_1} + \binom{T+a_2-1}{a_2} + \cdots
+\binom{T+a_s-(s-1)}{a_s}
\end{equation}
for certain integers $m-1 > a_1\geq a_2\geq \cdots \geq a_s\geq 0$, with
$s\geq 1$.
\end{cor}
\begin{proof}
Let $I\subseteq K[X]$, $I\neq (0)$,
be a homogeneous ideal with Hilbert
function $f=H_I$ and Hilbert polynomial $P$.
Macaulay's Theorem implies (see \cite{Bruns-Herzog},
Corollary~4.2.14) that
there exists an integer $n_0\in\Nn$ such that
$f(n+1)=f(n)^{\langle n\rangle}$ for all $n\geq n_0$.
We have $f(n)=P(n)$ for all $n\geq n_0$: Let
$$f(n_0)=\binom{c_{n_0}}{n_0} + \binom{c_{n_0-1}}{n_0-1}+\cdots
+\binom{c_1}{1}$$
be the $n_0$-th Macaulay representation of $f(n_0)$, and $j\geq 1$
minimal with $c_{j}\geq j$, so $c_{n_0}>c_{n_0-1}>\cdots>c_j\geq j\geq 1$.
Then for $n=n_0+k$ with $k\geq 0$:
\begin{align*}
f(n)&=\binom{c_{n_0}+k}{n_0+k} + \binom{c_{n_0-1}+k}{n_0-1+k}+\cdots
+\binom{c_j+k}{j+k}\\
%&=\binom{c_{n_0}+k}{c_{n_0}-n_0} + \binom{c_{n_0-1}+k}{c_{n_0-1}-n_0+1}+\cdots
%+\binom{c_j+k}{c_j-j}\\
&=\binom{n+a_1}{a_1}+\binom{n+a_2-1}{a_2}+\cdots+\binom{n+a_s-(s-1)}{a_s},
\end{align*}
where $a_{n_0-i+1}=c_{i}-i$ for $i=j,\dots,n_0$, and $s=n_0-j+1>0$.
We have $a_1=\deg P = \dim I -1 < m-1$, and hence $m-1 >
a_1\geq \cdots\geq a_s\geq 0$.
Conversely, suppose $P(T)$ is an integer-valued polynomial in the form given
in the corollary; we may assume $P\neq 0$.  
For $n\geq s$
$$P(n) = \binom{n+a_1}{n} + \binom{n+a_2-1}{n-1} + \cdots
+\binom{n+a_s-(s-1)}{n-(s-1)},$$
with $n+a_1>n+a_2-1>\cdots>n+a_s-(s-1)>0$ is the $n$-th Macaulay representation of
$P(n)$. Hence $P(n+1)=P(n)^{\langle n\rangle}$ for all $n\geq s$.
The $n$-th Macaulay coefficients of $\binom{n+m-1}{n}$ are
$(n+m-1,0,\dots,0)$. Since $s+a_1<s+m-1$ it follows that $P(n)<
\binom{n+m-1}{n}$ for all $n\geq s$.
Define $f\colon\Nn\to\Nn$ by $f(n)=\binom{n+m-1}{n}$ for $n< s$
and $f(n)=P(n)$ for $n\geq s$.
Then $f(n+1)\leq f(n)^{\langle n\rangle}$ for all $n$.
Moreover, if $s>1$, then $f(1)=m$, and if $s=1$, then $f(1)=a_1+1<m$.
By
Theorem~\ref{Macaulay} it follows that there exists a homogeneous ideal
$I\subseteq K[X]$, $I\neq (0)$, with $f=H_I$, and hence $P=P_I$.
\end{proof}

\noindent
The integers $a_1,\dots,a_s$ describing a Hilbert polynomial $P$ as in
the corollary
are uniquely determined by $P$. For a homogeneous ideal $I$ of $K[X]$ let
$$n_0(I):=\min\bigl\{n_0\in\Nn : \text{$H_I(n+1)=H_I(n)^{\langle n\rangle}$
for all $n\geq n_0$}\bigr\}.$$
If $E$ is a lex-segment, then $n_0(I_E)$ agrees with 
the largest degree of a minimal generator of $E$, see 
\cite{Bruns-Herzog}, Corollary~4.2.9.
Note that given a Hilbert polynomial
$P$ of a non-zero homogeneous ideal as in \eqref{Hilbert-Poly},
the integer
$$\varphi(P):= \min\bigl\{ n_0(I) : 
\text{$I\subseteq K[X]$ homogeneous ideal with $P_I=P$} \bigr\}$$
coincides with $s$. (By the proof of Corollary~\ref{Macaulay-Corollary}.) 
We put $\varphi(I):=\varphi(P_I)$ for any non-zero homogeneous ideal
$I$ of $K[X]$.
We have $H_I(n)=P_I(n)$ for all $n\geq \varphi(I)$.
We also note:

\begin{cor}\label{Macaulay-Corollary-2}
Let $P(T)\in\Qq[T]$ and $\deg P<m-1$.
If $P$ is the Hilbert polynomial of some non-empty final segment of $\Nn^n$
\rom{(}for some $n$\rom{)}, then $P$ is the Hilbert polynomial of
some non-empty final segment of $\Nn^m$. \qed
\end{cor}

\subsection*{The ordered set of Hilbert polynomials.}
Let us write
$$\cH_m:=\bigl\{H_E : E\in\cF(\Nn^m) \bigr\}$$
for the set of Hilbert functions of final segments of
$\Nn^m$, and put $\cH:=\bigcup_m \cH_m$.
We consider $\cH$ as an ordered set via the product ordering:
$$H_E\leq H_F \qquad:\Longleftrightarrow\qquad
\text{$H_E(s)\leq H_F(s)$ for all $s$.}$$
We have a strictly increasing surjection
$$\bigl(\cF(\Nn^m),\supseteq\bigr)\to \cH_m\colon E\mapsto H_E.$$
Hence by Corollary~\ref{Noetherian-Cor}:

\begin{cor}\label{Hilbert-Wellorder}
The ordered set $\cH_m$ is Noetherian. \qed
\end{cor}

\begin{remark}
In fact, the ordered set $\cH$ is also Noetherian. This can be shown
using Nash-Williams' theory of ``better-quasi-orderings'', see
\cite{maschenb-hemmecke}.
\end{remark}
\noindent
We write
$$\cP_m:=\bigl\{P_E : E\in\cF(\Nn^m)\bigr\}$$
for the set of Hilbert polynomials (of final segments of
$\Nn^m$), and $\cP:=\bigcup_m \cP_m$.
We totally order
$\cP$ via the dominance ordering $\preceq$.
%$$P_E \preceq P_F \qquad:\Longleftrightarrow\qquad
%\text{$P_E(s)\leq P_F(s)$ for all sufficiently large $s$.}$$
Clearly $H_E\leq H_{F} \Rightarrow P_E\preceq P_{F}$, so
$$\cH_m\to \cP_m\colon H_E\mapsto P_E$$
is an increasing surjection. 
A variant of the following fundamental fact 
has first been proved by Sit \cite{Sit} using different
methods:

\begin{cor}\label{HilbertPoly-Wellorder}
The
dominance ordering on the set $\cP$ of Hilbert polynomials is a well-ordering.
\end{cor}
\begin{proof}
By Corollary~\ref{Hilbert-Wellorder} and the preceding remarks,
$\cP_m$ is well-ordered, for every $m$. Moreover,
the leading coefficients of polynomials $P,Q\in\cP$ are positive, so if
$\deg P < \deg Q$, then $P\prec Q$. This implies that for every
decreasing sequence $P_0\succeq P_1 \succeq\cdots$ in $\cP$ 
there exists some $m$ such that $P_i\in\cP_m$ for all $i\gg 0$, and hence
$P_i=P_{i+1}$ for all $i\gg 0$.
This shows that $\cP$ is well-ordered.
\end{proof}

%Proposition~\ref{KeyProp} yields:
%
%\begin{cor}\label{KeyProp-Cor}
%For any function $f\colon\Nn\to\Nn$ there exists a finite set
%$\cM_{\succeq f}$ of final segments of $\Nn^m$ such that for
%every final segment $E$ of $\Nn^m$:
%$$P_E(s) \geq f(s) \text{\ for all $s\gg 0$}
%\quad\Longleftrightarrow\quad \text{$E\subseteq F$ for
%some $F\in\cM_{\succeq f}$.}$$
%\qed
%\end{cor}
%\noindent

\noindent
The following will be used in \cite{maschenb-pong-rd}:

\begin{cor}\label{Add-Cor}
If $P,Q\in\cP$ then $P+Q\in\cP_M$ where $M=\max\{\deg P,\deg Q\}+2$ 
and $P\cdot Q \in \cP_N$ where $N=(\deg P+2)(\deg Q+2)$. 
\rom{(}In particular, $\cP$ is a sub-semiring of the 
ring of all integer-valued polynomials.\rom{)}
\end{cor}
\begin{proof}
By Corollary~\ref{Macaulay-Corollary-2} we have
$P=P_I$ and $Q=P_J$ for non-zero monomial ideals $I\subseteq \Qq[X]$ and
$J\subseteq \Qq[Y]$,
where $X=\{X_1,\dots,X_m\}$ and $Y=\{Y_1,\dots,Y_n\}$
are disjoint sets of distinct indeterminates and
$m=\deg P+2$, $n=\deg Q+2$.
Consider the homomorphism
of graded $\Qq$-algebras
$$\Qq[X,Y] \longrightarrow \Qq[X]/I \oplus \Qq[Y]/J$$
defined by $$X_i\mapsto x_i=X_i+I\text{ and } 
Y_j\mapsto y_j=Y_j+J \text{ for all $i,j$.}$$
It is easy to see that its kernel $K_0$ 
is generated by $I\cup J\cup \{X_iY_j:
1\leq i\leq m, 1\leq j\neq n\}$ 
(hence is a monomial ideal of $\Qq[X,Y]$), and
$$\dim_\Qq \bigl(\Qq[X,Y]/K_0\bigr)_s=
\dim_\Qq \bigl(\Qq[X]/I\bigr)_s + \dim_\Qq
\bigl(\Qq[Y]/J\bigr)_s$$
for all $s$ except possibly $0$. In particular,
$P_{K_0}=P_I+P_J\in\cP_M$. As to $P\cdot Q$, it is well-known that
$P\cdot Q=P_S$, where $S\subseteq\Qq[Z_1,\dots,Z_N]$ is the homogeneous
ideal corresponding to the image of $V(I)\times V(J)\subseteq \Pp^{m-1}
\times\Pp^{n-1}$ in $\Pp^{N-1}$ under the Segre embedding.
\end{proof}

\begin{remark}
We write $\cP_m^* := \cP_m\setminus\{P_\emptyset\}$, where
$P_\emptyset=\binom{T+m-1}{m-1}$. By the proof of the corollary
$\cP_m^*$ is closed under addition. 
\end{remark}

\subsection*{Hilbert-Samuel polynomials.}
We now associate another integer-valued polynomial to a homogeneous
ideal. Given a homogeneous ideal $I$ of $K[X]$ and
$R=K[X]/I$ as in the beginning of this section, 
where $K$ is a field, let $$h_I(s)=\dim_K (R_0\oplus\cdots\oplus R_s)\qquad
\text{for
$s\in\Nn$.}$$ We call the function $h_I\colon\Nn\to\Nn$ the 
{\em Hilbert-Samuel function}\/ of the ideal $I$. We put $h_{(1)}(s)=0$
for all $s$. If $E\in\cF(\Nn^m)$ we put $h_E:=h_{I_E}$, where $I_E$ is the
monomial ideal in $\Qq[X]$ corresponding to $E$. 
With the notation
$Z_{\leq s} := \bigl\{ z\in\Nn^m : \abs{z}\leq s\bigr\}$ for $Z\subseteq\Nn^m$ and
$s\in\Nn$ we then have
$h_E(s)=\abs{V_{E,{\leq s}}}$ for all $s$.

\begin{lemma}\label{HilbertLemma}
Given a homogeneous ideal $I$ of $K[X]$
there exists an integer-valued polynomial $p_I$ of degree $\leq m$ such that
$p_I(s)=h_I(s)$ for all $s\gg 0$ in  $\Nn$.
\end{lemma}
\noindent
Indeed, the function $h_I$ is nothing but
the Hilbert function of the homogeneous ideal $IS$ of the polynomial ring
$S:=K[X_0,X_1,\dots,X_m]$.
We call the polynomial $p_I$ the {\em Hilbert-Samuel polynomial}\/ of 
the homogeneous ideal $I$. We put $p_E:=p_{I_E}$ for $E\in\cF(\Nn^m)$.

\begin{lemma}\label{Macaulay-Corollary-Lemma}
A polynomial $p(T)\in\Qq[T]$ is the Hilbert-Samuel polynomial of a
non-zero homogeneous ideal of $K[X_1,\dots,X_m]$ if and only if $p$ is
the Hilbert polynomial of some non-zero homogeneous ideal of
$K[X_0,X_1,\dots,X_m]$.
\end{lemma}
\begin{proof}
The ``only if'' part follows from the preceding discussion. 
Conversely, suppose there exists a homogeneous 
ideal $I\neq (0)$ of $K[X_0,\dots,X_m]$ such that 
$p=P_{I}$. We may assume that the maximal ideal $(X_0,\dots,X_m)$ is not
an associated prime of $I$. (Otherwise, $p=P_{I}=0$.)
Then, for a generic linear form $h$ of
$K[X_0,\dots,X_m]$,
multiplication by $h$ on $R=K[X_0,\dots,X_m]/I$ is injective (see, e.g., 
\cite{Bruns-Herzog},
Proposition~1.5.12). So we have a short exact
sequence of graded $K$-algebras and degree $0$ maps:
$$0\longrightarrow R(-1) \overset{h}{\longrightarrow} R \longrightarrow
S\longrightarrow 0,$$
where $R(-1)=\bigoplus_{s\geq 0} R(-1)_s$ with $R(-1)_0=0$,
$R(-1)_s=R_{s-1}$ for $s\geq 1$,  and $S=K[X_0,\dots,X_m]/J$ with $J=I+(h)$.
Hence
$$H_{I}(n)-H_{I}(n-1) = H_{J}(n) \qquad\text{for $n\geq 1$}$$
and so $H_{I}(n)=\sum_{i=0}^n H_{J}(i)$ for all $n$. 
Note that $H_{J}(1)=H_{I}(1)-1\leq m$, so
by Macaulay's Theorem there exists a non-zero homogeneous ideal $J'$ of
$K[X_1,\dots,X_m]$ with $H_J=H_{J'}$. Hence $H_I=h_J=h_{J'}$ and thus
$p=P_{I}=p_{J'}$.
\end{proof}

\begin{remark}
From the  characterization of Hilbert polynomials which was established in
Corollary~\ref{Macaulay-Corollary}, together with the
previous lemma, we obtain an analogous one for 
Hilbert-Samuel polynomials.
%(A related, but more complicated, description is given in
%\cite{ddp}, Chapter~2.) 
For the empty 
final segment $\emptyset\subseteq\Nn^m$
we have $p_{\emptyset}(T)=\binom{T+m}{m}$.
If $p(T)\in\Qq[T]$ has degree $m-1$ and is of the form 
$p=p_F$ for some non-empty final segment $F$ of $\Nn^n$
\rom{(}for some $n$\rom{)}, then there
exists a non-empty final segment $E$ of $\Nn^m$ such that $p=p_E$. 
\end{remark}
\noindent
Somewhat more generally, we can also 
define the {\it Hilbert function $H_E$}\/ and
{\it Hilbert polynomial $P_E$} of an $n$-tuple $E=(E_1,\dots,E_n)$ of 
final segments of $\Nn^m$ by setting
$$H_{E}=H_{E_1}+\cdots+H_{E_n}, \qquad
P_E=P_{E_1}+\cdots+P_{E_n}.$$
Similarly we define the {\em Hilbert-Samuel function}\/ $h_E$ and the
{\em Hilbert-Samuel polynomial}\/ $p_E$ of $E$.
(We will use these constructions in \cite{maschenb-pong-rd}.)
Given $n$-tuples
$E=(E_1,\dots,E_n)$ and $F=(F_1,\dots,F_n)$ of final segments of $\Nn^m$, 
we will write
$E\supseteq F$ if $E_i\supseteq F_i$ for all
$i=1,\dots,n$; that is, $\supseteq$ denotes the product order on 
${\cal F}(\Nn^m)^n$. 
%Then $E\subseteq F$ implies $H_E\geq H_F$ and
%$P_E\succeq P_F$, and $h_E\geq h_F$,
%$p_E\succeq p_F$. 
The map that assigns to $E$ 
its Hilbert function has finite fibers. In fact:

\begin{lemma}
The maps $$E\mapsto H_E, \quad E\mapsto h_E, \quad E\mapsto p_E,$$
where $E=(E_1,\dots,E_n)\in \cF(\Nn^m)^n$, are strictly increasing
and hence have finite fibers.
\end{lemma}
\begin{proof}
Let
$E$ and $E'$ be $n$-tuples of final segments of $\Nn^m$ with $E\supset E'$.
Then clearly $H_E\leq H_{E'}$ (thus $h_E\leq h_{E'}$ and $p_E\leq p_{E'}$) and
$H_E\neq H_{E'}$. 
Say $H_E(s_0) < H_{E'}(s_0)$ for some $s_0$; then
$p_E(s) < p_{E'}(s)$ for all $s\geq s_0$ sufficiently large.
The rest now follows from
Lemma~\ref{Folk-Lemma},~(1) and the Noetherianity of $\cF(\Nn^m)^n$.
\end{proof}
\noindent
It might be worth pointing out that although in general there are
infinitely many final segments of $\Nn^m$ with a given Hilbert polynomial $P$
(for example, every non-empty final segment of $\Nn$ has Hilbert polynomial
$0$), it is not difficult to see that
there always exists a {\it smallest}\/ 
Hilbert function $H_E$ with $P_E=P$. 

\subsection*{The ordered set of Hilbert-Samuel polynomials}
In the following, we will write
$$\cS_m:=\bigl\{p_E : E\in\cF(\Nn^m) \bigr\}$$
for the set of Hilbert-Samuel polynomials of final segments of $\Nn^m$, 
and $\cS:=\bigcup_m \cS_m$. We also let $\cS_{m}^*:=\cS_m\setminus
\{p_\emptyset\}$, with $p_\emptyset=\binom{T+m}{m}$.
By Lemma~\ref{Macaulay-Corollary-Lemma}, we have
$\cS_m^*=\cP_{m+1}^*$ for all $m$. By
Corollaries~\ref{HilbertPoly-Wellorder} and \ref{Add-Cor}:

\begin{cor}
The set $\cS$ is a well-ordered sub-semiring of the ordered ring 
of integer-valued polynomials. \qed
\end{cor}
\noindent
Corollary~\ref{Macaulay-Corollary} of Macaulay's Theorem
makes it possible to describe the unique 
isomorphism between $\cS_m$ and its order type $\ot(\cS_m)$
in a rather explicit way:

\begin{definition}\label{Defpsi}
Every Hilbert-Samuel polynomial $p(T)$ of a non-empty final segment of $\Nn^m$ 
can be written uniquely in the form
$$p(T) = \binom{T+a_1}{a_1} + \binom{T+a_2-1}{a_2} + \cdots
+\binom{T+a_s-(s-1)}{a_s}$$
for certain integers $m > a_1\geq a_2\geq \cdots \geq a_s\geq 0$ and
$s\geq 0$. We put
$$\mbf{c}_p := (c_{m-1},c_{m-2},\dots,c_0)\in\Nn^m,$$
where $c_i$ denotes 
the number of occurrences of $i\in\{0,\dots,m-1\}$ among the coefficients
$a_1,a_2,\dots,a_s$. Note that $\abs{\mbf{c}_p}=s$.
 We define an ordinal
$$\psi_p := \omega^{m-1}c_{m-1} + \omega^{m-2}c_{m-2} + \cdots + c_0.$$
For the Hilbert-Samuel polynomial $p(T)=\binom{T+m}{m}$ of the empty 
subset of $\Nn^m$ 
we set $\psi_p := \omega^m$.
\end{definition}

\noindent
The following observation is now easy. 
(This gives another proof of the well-orderedness of $\cS_m$.)

\begin{cor}\label{PsiCor}
The map
$$p \mapsto \psi_p \colon \cS_m \to \omega^{m}+1$$
is an isomorphism of ordered sets. \qed
\end{cor}

\noindent
In \cite{maschenb-pong-rd} we will use the previous corollary to define a
new model-theoretic rank for definable sets in differentially closed fields
of characteristic zero (via their Kolchin polynomials, see \cite{Kolchin}
or \cite{ddp}).

\subsection*{Computation of $\psi_p$}
Here is how the $\mbf{c}_p$ can be computed recursively, following
\cite{ddp}. (In \cite{ddp}, $c_{0},\dots,c_{m-1}$ are called the
{\em minimizing coefficients}\/ of $p$.) Let 
$$p(T)=b_d\binom{T+d}{d}+b_{d-1}\binom{T+d-1}{d-1}+\cdots+b_0\binom{T+0}{0}$$
with $b_0,\dots,b_d\in\Zz$, $b_d\neq 0$, be an integer-valued polynomial
of degree $d$. We define a sequence $\tilde{\mbf{c}}_p\in\Zz^m$,
where $m=d+1$, by
induction on $d$ as follows: If $p=0$ or $d=0$, so $p(T)=b_0$ is constant, 
we put
$\tilde{\mbf{c}}_p:=(b_0)$. If $d>1$, we consider the integer-valued
polynomial
$$q(T) := p(T+b_d) - \binom{T+d+1+b_d}{d+1} + \binom{T+d+1}{d+1}.$$
Note that $e:=\deg q <d$, so $\tilde{\mbf{c}}_q=(\tilde{c}_{q,e},\dots
\tilde{c}_{q,0})\in\Zz^{e+1}$ has been defined already. We
let $\tilde{\mbf{c}}_p := (b_d,0,\dots,0,\tilde{c}_{q,e},\dots
\tilde{c}_{q,0}) \in \Zz^m$.

\begin{lemma}
$p(T)\in\cS_m^*$ if and only if $\tilde{\mbf{c}}_p \geq 0$, and in this
case $\tilde{\mbf{c}}_p=\mbf{c}_p$.
\end{lemma}
\begin{proof}
We proceed by induction on $d$. The case $d=0$ is trivial.
Suppose $d>0$, and assume
first that $p(T)\in\cS_m^*$, say $p(T)=p_I(T)$ for some
non-zero monomial ideal $I$ of $R=\Qq[X_1,\dots,X_m]$. 
Since $p(s)>0$ for $s\gg 0$, we clearly have $b_d>0$. For $i=1,\dots,m$
we let  $\nu_i\in\Nn$ be the smallest natural number such that
$X_i^{\nu_i}X^\mu\in I$ for some $\mu\in\Nn^m$ with $\mu_i=0$. 
Multiplication by
$X^\nu$, where $\nu=(\nu_1,\dots,\nu_m)$, induces a short exact sequence
$$0 \longrightarrow \bigl(R/(I:X^\nu)\bigr)(-\abs{\nu}) \overset{X^\nu}{\longrightarrow}
R/I \longrightarrow R/(X^\nu) \longrightarrow 0.$$
Hence for all $s$: $$H_{(I:X^\nu)}(s)=
H_{I}(s+\abs{\nu})-H_{(X^\nu)}(s+\abs{\nu}).$$
Using the short exact sequences
$$0 \longrightarrow \bigl(R/(I:X_i^{\nu_i})\bigr)(-\abs{\nu_i}) \overset{X_i^{\nu_i}}{\longrightarrow}
R/I \longrightarrow R/(X_i^{\nu_i}) \longrightarrow 0$$
for $i=1,\dots,m$ it is easy to see that $\abs{\nu}=b_d$.
It follows that 
$$h_{(I:X^\nu)}(s)=p_I(s+b_d)-\binom{s+d+1+b_d}{d+1} + \binom{s+d+1}{d+1}$$
for all $s\gg 0$, and therefore
$q=p_{(I:X^\nu)}$. By induction we get
that $\tilde{\mbf{c}}_q\geq 0$ and thus $\tilde{\mbf{c}}_p\geq 0$.
Conversely, suppose that $\tilde{\mbf{c}}_p\geq 0$.
By induction we may write $q=p_J$ for some monomial ideal $J$
of $R'=K[X_1,\dots,X_{m-1}]$. We put $I=(X_m^{b_d+1},X_m^{b_d}J)$, a
monomial ideal of $R=K[X_1,\dots,X_m]$. Then, using the short exact sequence
$$0 \longrightarrow \bigl(R/(X_m,J)\bigr)(-b_d) \overset{X_m^{b_d}}{\longrightarrow}
R/I \longrightarrow R/(X_m^{b_d}) \longrightarrow 0$$
and the fact that $R/(X_m,J) \cong R'/J$
we obtain $p=p_I$ as required. The identity $\tilde{\mbf{c}}_p=\mbf{c}_p$
follows from Corollary~\ref{PsiCor},
Lemma~\ref{Aux-Height-Lemma}~(1) and 
the observation 
that $(\cS_m^*,{\preceq})\to(\Nn^m,{\leq_\lex})\colon p\mapsto\tilde{\mbf{c}}_p$
is strictly increasing and surjective.
\end{proof}

%\begin{example}
%For the Hilbert-Samuel polynomial $p=p_E$ of $E=\Nn^m$, we have 
%$\mbf{c}_p=(0,0,\dots,0)$ and hence
%$\psi_p=0$.
%\end{example}

\begin{example}\label{OrdinaryExample}
Let $p(T)=a(T+1)+b \in\Zz[T]$ with $a,b\in\Zz$, $a\neq 0$. Then $p(T)$ 
is the Hilbert-Samuel polynomial
of a non-empty final segment of $\Nn^2$ if and only if $a>0$ and
$b+\binom{a}{2}\geq 0$. In this case,
the sequence $\mbf{c}_p$ is given by
$\left(a,b+\binom{a}{2}\right)$, 
so
$\psi_p = \omega a + b+\binom{a}{2}$.
(Writing $p(T)=dT+1-g$ this yields the well-known inequality
$g\leq \binom{d-1}{2}$ relating degree and genus
of a projective curve.)
\end{example}

\subsection*{Application: length of increasing chains of ideals}
The results of this section, in particular Corollary~\ref{PsiCor}, can
be used to study increasing chains of ideals in polynomial rings. We give
two applications. First let us prove the theorem stated in the introduction.
We denote the set of homogeneous ideals of $K[X]=
K[X_1,\dots,X_m]$, ordered by reverse inclusion, by $\cI_m$. Since
$K[X]$ is Noetherian, $\cI_m$ is well-founded.
We write $p\colon\cI_m\to\cS_m$ for
the map $I\mapsto p_I$.

\begin{lemma}\label{Div-Lemma}
$\hgt(\cI_m)=\hgt\bigl(\cF(\Nn^m)\bigr)=\omega^m+1$.
\end{lemma}
\begin{proof}
The first equality holds since there exists a strictly increasing surjection 
$\cI_m\to\cF(\Nn^m)$. 
This is a well-known consequence of the division 
algorithm in $K[X]$ (see, e.g., \cite{Eisenbud}, Chapter 15):
Choose a
term ordering $\leq$ on $X^\diamond$; given a non-zero polynomial $f\in K[X]$
let $\lm(f)$ be the {\em leading monomial}\/ of $f$, that is, the
largest monomial in the ordering $\leq$ which occurs in $f$ with a non-zero
coefficient. Given an ideal $I$ of $K[X]$ we denote by $\lm(I)$ the monomial
ideal generated by the $\lm(f)$, where $0\neq f\in I$.
Now suppose that $I\supset J$ are ideals in $K[X]$. Then $\lm(I)\supset\lm(J)$:
Choose $f\in I\setminus J$ such that $\lm(f)$ is minimal in the ordering
$\leq$; we claim that $\lm(f)\in\lm(I)\setminus\lm(J)$.
Otherwise $\lm(f)=\lm(g)$ for some $0\neq g\in J$, and we can write
$f=qg+r$ for some $q,r\in K[X]$,
$r\neq 0$, with $\lm(r)<\lm(f)$. Since $r\in I\setminus J$, this is a 
contradiction. Hence
the map which associates to a homogeneous ideal $I$ the
monomial ideal $\lm(I)$ is
strictly increasing.  The second equality
follows from Proposition~\ref{Bonnet-Pouzet}.
\end{proof}
\noindent
We now show:

\begin{theorem}
For every strictly increasing surjection $\varphi\colon\cI_m\to S$,
where $S$ is an ordered set, there
exists a non-decreasing map $\psi\colon\cS_m\to S$ such that
$\psi\circ p\leq\varphi$.
\end{theorem}
\begin{proof}
The map $I\mapsto\psi_{p_I}\colon\cI_m\to\omega^m+1$ is strictly
increasing and surjective. Hence $\psi_{p_I}=\hgt_{\cI_m}(I)$ 
for all $I\in\cI_m$, by the last lemma. 
In fact, $\hgt_{\cI_m}=\hgt_{\cS_m}\circ p$.
The claim now follows from Lemma~\ref{Aux-Height-Lemma}~(2).
\end{proof}
\noindent
(By the second part of Lemma~\ref{Aux-Height-Lemma}, the Hilbert-Samuel
polynomial $p\colon\cI_m\to\cS_m$
is characterized up to isomorphism 
by the property expressed in the theorem, in the
category of strictly increasing surjections $\cI_m\to S$, where $S$ is a totally
ordered set.)

\begin{remark}
Let $\cG_m$ denote the set of isomorphism classes of finitely generated
graded $R$-modules, where $R=K[X]$. 
We define a binary relation $\leq$ on $\cG_m$ by
$M\leq N \Longleftrightarrow$ there exists a surjective homomorphism of
graded $R$-modules $N\to M$. Since every surjective endomorphism of a 
finitely generated $R$-module is an isomorphism (see \cite{Eisenbud}), 
it follows that
$\leq$ is an ordering on $\cG_m$. By Noetherianity of $R$, $\leq$ is
well-founded. We ask: Does the theorem above remain true when $\cI_m$ is
replaced by $\cG_m$ and $p$ by the map which assigns to every $M\in\cG_m$
its Hilbert-Samuel polynomial?
\end{remark}

\noindent
By the theorem above, every strictly increasing chain of non-zero 
homogeneous ideals in $K[X]$
gives rise to a strictly
decreasing sequence in the lexicographically ordered set
$\omega^m$. What can be said about the length of such sequences?
For this, let us fix an increasing function $f\colon\Nn\to\Nn$, and
consider finite sequences
$$\nu_0 >_\lex \nu_1 >_\lex \cdots >_\lex \nu_{\ell-1}$$
of $m$-tuples $\nu_i\in\Nn^m$, strictly decreasing with respect to the
lexicographic ordering on $\Nn^m$, with the property that $|\nu_i|\leq f(i)$
for all $i$. For the purpose of this section, let us call such a sequence
an {\em $f$-bounded sequence}\/ in $\Nn^m$.
By K\"onig's Lemma (e.g., 
\cite{Kechris}, p.~20) applied to the tree whose nodes
are the $f$-bounded sequences it follows that there exists an $f$-bounded 
sequence with maximal length $\ell=\ell(m,f)$. 
It is not difficult to compute an explicit formula for $\ell(m,f)$:

\begin{lemma}\label{Length-Lemma}
We have $\ell(1,f)=f(0)+1$ and
$$\ell(m,f)=1+\ell(m-1,f_1)+\cdots+\ell(m-1,f_{f(0)}) \qquad\text{for $m>1$,}$$
with $f_i\colon\Nn\to\Nn$
defined by $$f_i(j)=f\bigl(j+1+\ell(m-1,f_1)+\cdots+\ell(m-1,f_{i-1})\bigr)-
f(0)+i$$
for $i,j\in\Nn$, $i\geq 1$.
\end{lemma}
\begin{proof}
By induction on $m$. The case $m=1$ is trivial.
Suppose that $m>1$, and let
$$\nu_0 >_\lex \nu_1 >_\lex \cdots >_\lex \nu_{\ell-1}$$
be an $f$-bounded sequence in $\Nn^m$ of maximal length $\ell=\ell(m,f)$.
We must have $\nu_0=\bigl(f(0),0,\dots,0\bigr)$; otherwise
(since $f$ is increasing)
$$\bigl(f(0),0,\dots,0\bigr)>_\lex\nu_0 >_\lex \nu_1 >_\lex \cdots >_\lex 
\nu_{\ell-1}$$
would be a longer $f$-bounded sequence. For a similar reason,
$\nu_1$ must have the form $\nu_1=\bigl(f(0)-1\bigr) \, \mu_0$ for some
$\mu_0\in\Nn^{m-1}$. It follows that
$$\nu_1=\bigl(f(0)-1\bigr) \, \mu_{0},\dots,
\nu_{\ell_1}=\bigl(f(0)-1\bigr) \, \mu_{\ell_1-1}$$ 
for some $f_1$-bounded sequence
$\mu_0>_\lex\mu_1>_\lex\cdots>_\lex\mu_{\ell_1-1}$ in $\Nn^{m-1}$
of maximal length $\ell_1=\ell(m-1,f_1)$. The next terms in the sequence
must then have the form 
$\nu_{\ell_1+i}=\bigl(f(0)-2\bigr) \, \lambda_{i-1}$ for some
$f_2$-bounded sequence
$\lambda_0>_\lex\lambda_1>_\lex\cdots>_\lex\lambda_{\ell_2-1}$ in $\Nn^{m-1}$
of maximal length $\ell_2=\ell(m-1,f_2)$, and so on. This leads to
the displayed formula for $\ell(m,f)$.
\end{proof}

\noindent
We can use this to show the following statement about 
uniform bounds for the length of ascending chains of homogeneous ideals.
Recall that for any homogeneous ideal $I$ of $K[X]=K[X_1,\dots,X_m]$, 
we denote by $\varphi(I)$ the smallest natural number $n_0$ such that
for any homogeneous ideal $J$ of $K[X]$ with Hilbert polynomial $P_J=P_I$,
we have $H_J(n+1)=H_J(n)^{\langle n\rangle}$ for all $n\geq n_0$.
(Cf.~remarks following Corollary~\ref{Macaulay-Corollary}.)

\begin{prop}\label{Greg-Prop}
Let $f\colon\Nn\to\Nn$ be any function
and $m\geq 1$. There exists a natural number $t_m(f)$
depending only on $m$ and $f$, and
primitive recursive in $f$, such that for any field $K$ and any
strictly increasing chain
$$I_0 \subset I_1 \subset \cdots \subset I_{t-1}$$
of non-zero homogeneous ideals in $K[X_1,\dots,X_m]$
such that $\varphi(I_i)\leq f(i)$ for all $i$,
we have $t\leq t_m(f)$.
\end{prop}
\noindent
Here as usual, a function $F\colon\Nn^r\to\Nn$ (for $r\in\Nn$) 
is called {\em primitive recursive}\/ in a given collection
$F_1,\dots,F_k$ of functions $F_i\colon\Nn^{r_i}\to\Nn$, $i=1,\dots,k$,
if it can be obtained from $F_1,\dots,F_k$ as well as the
the constant function $0$, the successor function
$x\mapsto x+1$, coordinate permutations $(x_1,\dots,x_n)\mapsto
(x_{\sigma(1)},\dots,x_{\sigma(n)})$, and the projections
$(x_1,\dots,x_{n+1})\mapsto (x_1,\dots,x_n)\colon \Nn^{n+1}\to\Nn^n$,
by finitely many applications of the following rules
(substitution and induction, respectively):
\begin{enumerate}
\item if $F\colon\Nn^r\to\Nn$ and $G_1,\dots,G_r\colon\Nn^s\to\Nn$ are
primitive recursive in $F_1,\dots,F_k$, then so is
$H=F(G_1,\dots,G_r)\colon\Nn^s\to\Nn$;
\item if $F\colon\Nn^r\to\Nn$ and $G\colon\Nn^{r+2}\to\Nn$ are primitive
recursive in $F_1,\dots,F_k$, then so is the function
$H\colon\Nn^{r+1}\to\Nn$ defined by
$$H(x,y)=\begin{cases}
F(x)&\text{if $y=0$} \\
G(x,y-1,H(x,y-1))& \text{if $y>0$,} \end{cases}$$
for $x\in\Nn^r$ and $y\in\Nn$.
\end{enumerate}
If $k=0$, we obtain the plain {\em primitive recursive}\/ functions.
Standard number-theoretic functions like addition $(x,y)\mapsto x+y$,
multiplication $(x,y)\mapsto x\cdot y$ or exponentiation $(x,y)\mapsto x^y$
are primitive recursive.
The class of primitive recursive functions forms a proper subclass of
all recursive, or computable, functions. A prominent example of a
recursive but not primitive recursive function is the Ackermann function.
(See, e.g., \cite{Moreno-Socias}, Section~2 
for the definition of the Ackermann function.)

Before we prove Proposition~\ref{Greg-Prop} we show the following lemma.
Given $Q(T)\in\Qq[T]$ we put $\Delta Q(T)=Q(T)-Q(T-1)\in\Qq[T]$.
Note that if $\Delta Q=P$ and $Q(0)=P(0)$, 
then $Q(n)=\sum_{i=0}^n P(i)$ for all $n\in\Nn$.

\begin{lemma}
Let
$$P(T) = \binom{T+a_1}{a_1} + \binom{T+a_2-1}{a_2} + \cdots
+\binom{T+a_s-(s-1)}{a_s}
$$
with integers $a_1\geq a_2\geq \cdots \geq a_s\geq 0$ and
$s\geq 1$. Then
$$
Q(T) = \binom{T+(a_1+1)}{a_1+1} + \binom{T+(a_2+1)-1}{a_2+1} + \cdots
+\binom{T+(a_s+1)-(s-1)}{a_s+1}
$$
is the unique polynomial in $\Qq[T]$ such that $Q(0)=P(0)$ and $\Delta Q=P$.
\end{lemma}
\begin{proof}
Clearly $Q(T)$ satisfies $Q(0)=1=P(0)$. The
well-known identity $\binom{a}{b}-\binom{a-1}{b}=\binom{a-1}{b-1}$ for
$a\geq b\geq 0$ implies $\Delta Q=P$ as required.
\end{proof}

\begin{proof}[Proof \rom{(Proposition~\ref{Greg-Prop})}]
First replacing $f$ by the function $g\colon\Nn\to\Nn$ defined by
by $i\mapsto \max\{f(0),\dots,f(i)\}$,
if necessary, we may assume that $f$ is increasing.
If $I$ is a non-zero homogeneous ideal in $K[X]$
with Hilbert-Samuel polynomial $p=p_I$ and $s=\varphi(I)$, then for
$n\geq s$ we have
$$p(n)=\sum_{i=0}^n H_I(i) = Q(n) + k$$
where $Q(n)=\sum_{i=0}^n P_I(i)$ and $k=
\sum_{i=0}^{s-1}\bigl(H_I(i)-P_I(i)\bigr)$. By the previous lemma it follows
that $\abs{\mbf{c}_p}=s+k$, and since $j\leq\binom{s-1+m}{m}$, we have
$\abs{\mbf{c}_p}\leq s+\binom{s-1+m}{m}=:h_m(s)$. Hence
the function $t_m(f):=\ell\bigl(m,h_m\circ f\bigr)$, with $\ell$ as defined in 
Lemma~\ref{Length-Lemma}, bounds the length of every
strictly increasing chain of ideals as in Proposition~\ref{Greg-Prop}.
It is a tedious but straightforward exercise, left to the reader, 
to verify that  $t_m(f)$ is primitive recursive in $f$, for given $m$.
\end{proof}
\noindent
The proposition above yields the following theorem of
Moreno Soc{\'\i}as \cite{Moreno-Socias}:

\begin{cor}\label{MS}
Let $f\colon\Nn\to\Nn$ be any function
and $m\geq 1$. There exists a natural number $t_m(f)$
which is primitive recursive in $f$ such that for any field $K$ and any
strictly increasing chain
\begin{equation}\label{MS-Eq}
I_0 \subset I_1 \subset \cdots \subset I_{t-1}
\end{equation}
of ideals in $K[X]=K[X_1,\dots,X_m]$
such that $I_i$ is generated by polynomials of degree at most $f(i)$,
for every $i$,
we have $t\leq t_m(f)$.
\end{cor}
\begin{proof}
First we show that we may restrict ourselves to
chains \eqref{MS-Eq} where each ideal $I_i$ is monomial. 
To see this choose a term ordering $\leq$ on
$X^\diamond$ which is degree-compatible. Then, as in the proof of the
first equality in Lemma~\ref{Div-Lemma}, one shows that
if $I\subset J$ are ideals and $J$ is generated by polynomials of
degree $\leq d$, then there exists a monomial $X^\nu\in\lm(J)\setminus\lm(I)$
of degree $\abs{\nu}\leq d$.
As in \cite{Moreno-Socias}, Section~4 one further reduces to the
case where every monomial ideal $I_i$ is of the form
$I_i=I_{E_i}$ for a lex-segment $E_i\subseteq\Nn^m$.
By the remarks following Corollary~\ref{Macaulay-Corollary} we then have 
$\varphi(I_i)\leq f(i)$ for all $i$. Hence
$t_m(f)$ as defined in the previous proposition works.
\end{proof}
\noindent
In \cite{Moreno-Socias} (Corollary~7.5)
it is also shown that $t_m(f)$ is not primitive
recursive in $m$, even for an affine function $f(i)=p+iq$ 
($p,q\in\Nn$). In fact, $m\mapsto t_m(f)$ grows like the Ackermann 
function, and hence extremely rapidly.

Moreno Soc\'\i{}as' result \ref{MS} 
may be interpreted as a quantitative variant of Dickson's Lemma (and thus,
of the
Hilbert Basis Theorem). We finish this section with outlining the proof
of a similar finitary formulation of Maclagan's principle.
This fact can be seen to provide primitive recursive complexity estimates
for algorithms whose
termination has been shown using the Noetherianity of $\cF(\Nn^m)$.
The proof is based on ideas of Harvey Friedman
\cite{Friedman}.
It also gives a different argument for Proposition~\ref{MS} in the
case where $f$ is affine,
by reducing to the case of ascending chains of monomial ideals, as in 
the argument in the beginning of the proof of \ref{MS}. 

\begin{prop}\label{PrimRecProp}
Let $p,q\in\Nn$, $m\geq 1$. There exists a natural number $r_m(p,q)$,
which is primitive recursive in $p$ and $q$, such that for any bad sequence
$$F_0, F_1, \dots, F_{r-1}$$
of final segments of $\Nn^m$, with $F_i$
generated by elements of degree at most $p+iq$,
we have $r\leq r_m(p,q)$.
\end{prop}
\begin{proof}[Sketch of the proof]
We fix $m\geq 1$.
Let $T$ be the first-order theory of the structure
$\mbf{N}=(\Nn,\leq)$ in the language $\cL_0$ consisting of
the binary relation 
symbol $\leq$ and a constant symbol for every element of $\Nn$.
Every model of $T$ is an ordered set
containing an isomorphic copy of $\mbf{N}$ as an initial segment, which we
identify with $\mbf{N}$. For $r\in\Nn$ we let $\cL_r$ be the language $\cL$ 
augmented by $m$-ary predicate symbols $F_0,\dots,F_{r-1}$.
For $p,q,r\in\Nn$ let $T_{p,q,r}$ be the
union of the $\cL_r$-theory $T$ together with sentences that express
that $F_0,\dots,F_{r-1}$ form a bad sequence of final segments, and each
$F_i$ is generated by elements of $\Nn^m$ of degree $\leq p+iq$.
Given $p,q\in\Nn$, the set $T_{p,q}=\bigcup_r T_{p,q,r}$
of sentences in the language $\cL=\bigcup_r \cL_r$ is inconsistent, 
by Noetherianity of $\cF(\Nn^m)$. The Completeness Theorem
of first-order logic implies that for some $r$, $T_{p,q,r}$ is inconsistent.
Clearly any such $r$ bounds the length of a bad sequence
in $\cF(\Nn^m)$ with the $i$-th element in the sequence
generated in degrees $\leq p+iq$.
In order to show that $r$ can be found primitive recursively in
$p,q$, we use some facts about  the so-called second
principal system of reverse mathematics $\WKL_0$, see \cite{Simpson-Book}:
First, the proof of the Completeness Theorem for countable languages
can be carried out in $\WKL_0$ (\cite{Simpson-Book}, Section~IV.3).
It is a routine exercise to verify that the
Noetherianity of $\cF(\Nn^m)$ is also provable in $\WKL_0$. 
(For example, the proof given in the next section can be easily formalized.)
The inconsistency of $T_{p,q,r}$ can be
expressed by an existential formula $\varphi(p,q,r)$ in the language of
arithmetic, and the $\forall\exists$-sentence
$\forall p\forall q\exists r\varphi(p,q,r)$ is provable in $\WKL_0$. 
By a theorem of Friedman and Harrington
(\cite{Simpson-Book}, IX.3) there exists a primitive recursive function
$\Nn^2\to\Nn\colon (p,q)\mapsto r_m(p,q)$ such that 
$\WKL_0 \vdash \forall p\forall q\varphi\bigl(p,q,r_m(p,q)\bigr)$.
This $r_m(p,q)$ has the required properties.
\end{proof}

\begin{remark}\label{PrimRecRemark}
The precise form of bounding function $i\mapsto p+iq$ used in Proposition~\ref{PrimRecProp}
is not essential: 
Let $g\colon\Nn^{k+1}\to\Nn$ be a primitive recursive function, $k\in\Nn$.
Then there exists a primitive recursive function $r_{m,g}
\colon\Nn^k\to\Nn$ such that that for any $p_1,\dots,p_k\in\Nn$ and any
bad sequence
$F_0, F_1, \dots, F_{r-1}$
of final segments of $\Nn^m$, with $F_i$
generated by elements of degree at most $g(p_1,\dots,p_k,i)$,
we have $r\leq r_{m,g}(p_1,\dots,p_k)$.
\end{remark}

\noindent
It is well-known that given an ideal $I=(f_1,\dots,f_n)$ in a 
polynomial ring $K[X]=K[X_1,\dots,X_m]$ over a field $K$, with $\deg f_i
\leq d$ for all $d$, the ideal $\lm(I)$ of leading monomials of elements
of $I$, with respect to a degree-compatible 
term ordering, can be generated by monomials of whose degree is
bounded by $d^{2^m}$
(see, e.g., \cite{Dube}). 
Since this bound is primitive recursive in $d$, Remark~\ref{PrimRecRemark}
(applied to $g\colon\Nn^3\to\Nn$ 
given by $g(p,q,i)=(p+iq)^{2^m}$)
implies:

\begin{cor}
Let $p,q\in\Nn$, $m\geq 1$. There exists a natural number $s_m(p,q)$,
which is primitive recursive in $p$ and $q$, such that for any 
field $K$ and any sequence
$$I_0, I_1, \dots, I_{s-1}$$
of ideals of $K[X_1,\dots,X_m]$ with 
$s\geq s_m(p,q)$ and with each $I_i$
generated by elements of degree at most $p+iq$, there exists
$0\leq i<j<s$ such that $\lm(I_i) \supseteq \lm(I_j)$
\rom{(}and hence in particular $H_{I_i}\leq H_{I_j}$\rom{)}. \qed
\end{cor}

\section{Total Orderings of Monomial Ideals}\label{OrderingsOfMonomialIdeals}

\noindent
In this section, we study the ordered 
set $\cF(\Nn^m)$ of final segments  of $\Nn^m$,
with the ordering given by the superset relation
(see Example~\ref{example-Nn-2}). 
We give an upper bound on $\ot\bigl(\cF(\Nn^m)\bigr)$
and we explicitly describe several ways of extending the ordering on
$\cF(\Nn^m)$ to a well-ordering. Finally, 
we compute the order type of one particularly
useful ordering, called the Kleene-Brouwer ordering of $\cF(\Nn^m)$.

\subsection*{Bounding the type of $\cF(\Nn^m)$.}
Our computation of an upper bound for the type of $\cF(\Nn^m)$ is
based on the following idea (which, incidentally, gives yet another proof of
the Noetherianity of $\cF(\Nn^m)$). Recall that an ideal in a commutative
ring is called {\em irreducible}\/ if it cannot be written as the intersection
of two strictly larger ideals. (For example, prime ideals are irreducible.)
By Noetherianity of $K[X]$, every ideal $I$ of $K[X]$ (where $K$ is a field)
can be written as an
intersection of irreducible ideals. Such a representation
$I=J_1\cap\cdots\cap J_r$ ($r\in\Nn$, $r>0$)
of $I$ as an intersection of irreducible ideals
$J_1,\dots,J_r$, however, is not
necessarily unique, even if we require it to be {\em irredundant,} that
is, $J_i \not\subseteq J_j$ for $i\neq j$. However, an irredundant
decomposition $I=J_1\cap\cdots\cap J_r$ of a {\it monomial}\/ 
ideal $I$ is unique,
and in this case the irreducible components
$J_i$ are monomial ideals as well. %(See \cite{Miller}
%for an elegant proof of this uniqueness.) 
It is easy to see that 
every irreducible monomial ideal
is of the form ${\frak m}^\nu := (X_i^{\nu_i}:\nu_i>0)$ for some
$\nu=(\nu_1,\dots,\nu_m)\in\Nn^m$. Note that 
\begin{multline*}
[\nu^{(1)},\dots,\nu^{(r)}] \leq^\diamond
[\mu^{(1)},\dots,\mu^{(s)}]\text{ in $(\Nn^m)^\diamond$} \quad\Rightarrow\quad\\
{\frak m}^{\nu^{(1)}}\cap\cdots\cap {\frak m}^{\nu^{(r)}} \supseteq
{\frak m}^{\mu^{(1)}}\cap\cdots\cap {\frak m}^{\mu^{(s)}},
\end{multline*} for all
$\nu^{(i)},\mu^{(j)}\in\Nn^m$ such that $\supp\bigl(\nu^{(i)}\bigr)\supseteq
\supp\bigl(\mu^{(j)}\bigr)$ for all $i,j$.  
Here $\supp\nu=\{i:\nu_i>0\}$ denotes the
{\em support}\/ of $\nu=(\nu_1,\dots,\nu_m)\in\Nn^m$. (See Section~\ref{NoetherianOrderedSetsSection} for
the definition of $\leq^\diamond$.)
Given a vector $\nu=(\nu_1,\dots,\nu_m)\in\Nn^m$ let us write
$\langle\nu\rangle:=(\nu_{i_1},\dots,\nu_{i_k})\in\Nn^k$,
where $1\leq i_1<\cdots<i_k\leq m$ are the elements of the support of
$\nu$ listed in increasing order.
Given a final segment $E\in\cF(\Nn^m)$ and $\sigma\subseteq\{1,\dots,m\}$
we denote by $\varphi(E,\sigma)$ 
the commutative word
$\bigl[ \langle\nu^{(1)}\rangle,\dots,\langle\nu^{(r)}\rangle\bigr]\in (\Nn^{\abs{\sigma}})^\diamond$,
where ${\frak m}^{\nu^{(1)}},\dots,{\frak m}^{\nu^{(r)}}$ are the irreducible
components of $I_E$
with $\supp\nu^{(i)}=\sigma$. Here $I_E$ is the monomial ideal
of $\Qq[X_1,\dots,X_m]$ corresponding to $E$.
Combining the various $\varphi(E,\sigma)$ we
obtain a quasi-embedding
\begin{equation}\label{QE-m}
\cF(\Nn^m) \to \prod_{\sigma\subseteq\{1,\dots,m\}} \bigl(\Nn^{\abs{\sigma}}\bigr)^\diamond\colon E\mapsto \bigl(\varphi(E,\sigma):\sigma\subseteq\{1,\dots,m\}\bigr).
\end{equation}
By Proposition~\ref{oProp} it follows that
$$\ot\bigl(\cF(\Nn^m)\bigr) \leq \bigotimes_{\sigma\subseteq\{1,\dots,m\}}
\ot\left(\bigl(\Nn^{\abs{\sigma}}\bigr)^\diamond\right).$$
In order to continue our majorization, we need to bound the type of 
the Noetherian ordered sets $\bigl(\Nn^{\abs{\sigma}}\bigr)^\diamond$.
Recall that $\varepsilon_0$ is the supremum of the sequence of ordinals
$\omega,\omega^\omega,\omega^{\omega^\omega},\dots$; in other words,
$\varepsilon_0$ is the smallest solution of the equation $\omega^x=x$ in
ordinals. Given an ordinal $\alpha$ we define 
$$\alpha'=\begin{cases}
\alpha&\text{if $\alpha<\varepsilon_0$,} \\
\omega\otimes\alpha&\text{if $\alpha\geq\varepsilon_0$.}
\end{cases}$$
We will show:

\begin{lemma}\label{vdDE-Lemma}
Let $S$ be a Noetherian ordered set of type $\alpha=\ot(S)$. Then
$S^\diamond$ is Noetherian of type $\ot(S^\diamond)\leq\omega^{\alpha'}$.
\end{lemma}

\noindent
In \cite{Schmidt} one finds that $\ot(S^\ast)=\omega^{\omega^{\alpha^*}}$, 
for every
Noetherian ordered set of type $\alpha=\ot(S)$. Here 
$$\alpha^* = \begin{cases}
\alpha-1 &\text{if $0<\alpha<\omega$,} \\
\alpha+1 &\text{if $\alpha=\varepsilon+n$ for some $n<\omega$ and some
$\varepsilon$ with $\varepsilon=\omega^\varepsilon$,} \\
\alpha &\text{otherwise.}
\end{cases}$$
This yields
the cruder upper bound $\ot(S^\diamond) \leq \omega^{\omega^{\alpha^*}}$. 
The lemma above was inspired by the following consequence of it, a quantitative
version of a well-known result of B. H. Neumann.
(See the erratum \cite{vdDries-Ehrlich-Erratum} to
\cite{vdDries-Ehrlich}.)

\begin{cor} \rom{(van den Dries-Ehrlich.)}
Let $\Gamma$ be an ordered abelian group and
$S\subseteq\Gamma^{\geq 0}$ well-ordered of order type $\alpha=\ot(S)$.
Then the monoid $[S]$ generated by $S$ in $\Gamma$ is well-ordered of
order type $\leq \omega^{\alpha'}$.
\end{cor}
\begin{proof}
Since $S^\diamond$ is the free commutative monoid generated by $S$, we have
a natural surjective monoid homomorphism $S^\diamond \to [S]$.
This homomorphism is increasing when $S^\diamond$ is equipped with the
ordering $\leq^\diamond$ 
and $[S]$ with the well-ordering induced from $\Gamma$.
The claim now follows from the last lemma and Proposition~\ref{oProp},~(3).
\end{proof}

\noindent
Lemma~\ref{vdDE-Lemma} together with \eqref{QE-m} yields the following
upper bound on the type of $\cF(\Nn^m)$:
\begin{equation}\label{UpperBound}
\ot\bigl(\cF(\Nn^m)\bigr) \leq \bigotimes_{\sigma\subseteq\{1,\dots,m\}}
\omega^{\omega^{\abs{\sigma}}} = \omega^{\bigoplus_\sigma \omega^{\abs{\sigma}}}= \omega^{(\omega + 1)^{\otimes m}},
\end{equation}
where $\alpha^{\otimes m}=\alpha\otimes\alpha\otimes\cdots\otimes\alpha$
($m$ times) for $\alpha\in\On$. (If the Cantor normal form of $\alpha$ has
leading term $\omega^{\gamma}$, then the leading term of $\alpha^{\otimes m}$
is $\omega^{\gamma m}$. This implies the bound on
$\ot\bigl(\cF(\Nn^m)\bigr)$ given in the introduction.)

\begin{proof}[Proof of Lemma~\ref{vdDE-Lemma}]
We proceed by transfinite induction on $\alpha$. The case $\alpha=0$ is
trivial ($S=\emptyset$), so let $\alpha>0$. We distinguish two cases.
First suppose that $\alpha$ is not additively indecomposable, that is,
$\alpha=\alpha_1\oplus\alpha_2$ for some ordinals $\alpha_1,\alpha_2<\alpha$.
Hence $S$ is a disjoint union $S=S_1\cup S_2$ with $\ot(S_1)\leq\alpha_1$
and $\ot(S_2)\leq\alpha_2$. (Here each $S_i$ is equipped with the
restriction of the ordering of $S$ to $S_i$.) We have a
bijective increasing map $S_1^\diamond \times S_2^\diamond \to S^\diamond$,
so $$\ot(S^\diamond)\leq \ot\bigl(S_1^\diamond \times S_2^\diamond\bigr)=
\ot(S_1^\diamond)\otimes\ot(S_2^\diamond)\leq
\omega^{\alpha_1'}\otimes\omega^{\alpha_2'}=\omega^{\alpha'}$$ 
by Proposition~\ref{oProp} and using the
induction hypothesis. Now suppose that $\alpha$ is additively indecomposable.
It is well-known that then $\alpha$ has the form $\alpha=\omega^\beta$ for
some $\beta>0$. By Proposition~\ref{oProp},~(1) it suffices to show
that $\ot\bigl((S^\diamond)^{\not\geq w}\bigr)<\omega^{\alpha'}$ for
all $w\in S^\diamond$. We show this by induction on $\abs{w}$. For
$\abs{w}=0$ there is nothing to show, since then $(S^\diamond)^{\not\geq w}=
\emptyset$. 
Suppose $\abs{w}>0$, say $w=[s_0,\dots,s_{m-1}]$ with
$s_0,\dots,s_{m-1}\in S$. There exists a 
quasi-embedding
$$\psi\colon (S^\diamond)^{\not\geq w}\to
(S^{\not\geq s_0})^\diamond\amalg
\left(S\times (S^\diamond)^{\not\geq w'}\right),$$ where $w'=[s_1,\dots,s_{m-1}]$. 
In order to see this, let $v=[t_0,\dots,t_{n-1}]\in S^\diamond$ with
$v\not\geq w$. Then either $t_i\not\geq s_0$ for
all $i$; or $t_i\geq s_0$ for some $i$, so
after reordering the $t$'s we may assume 
$t_0\geq s_0$, and $v'=[t_1,\dots,t_{n-1}]\not\geq
[s_1,\dots,s_{n-1}]$. In the first case, we
put $\psi(v)=v\in (S^{\not\geq s_0})^\diamond$, and in the second case, we
put $\psi(v)=(t_0,v')\in S\times (S^\diamond)^{\not\geq w'}$.  It is
easy to check that $\psi$ is a quasi-embedding. Hence, by 
Proposition~\ref{oProp},
$$\ot\left((S^\diamond)^{\not\geq w}\right) \leq 
\ot\bigl((S^{\not\geq s_0})^\diamond\bigr)\oplus
\left(\alpha\otimes \ot\bigl((S^\diamond)^{\not\geq w'}\bigr)\right).$$
Put $\gamma=\ot(S^{\not\geq s_0})$, so $\gamma<\alpha$ and hence
$\ot\bigl((S^{\not\geq s_0})^\diamond\bigr)\leq\omega^{\gamma'}<
\omega^{\alpha'}$ by inductive hypothesis on $\alpha$. 
By inductive hypothesis on $w$ we have $\delta := \ot\bigl((S^\diamond)^{\not\geq w'}\bigr)<\omega^{\alpha'}$. Hence it suffice to show that
$\alpha\otimes\delta<\omega^{\alpha'}$. Write $\delta$ in Cantor normal
form as $\delta=\omega^{\delta_1}n_1+\cdots+\omega^{\delta_k}n_k$ with
ordinals $\delta_1>\cdots>\delta_k$ and positive integers $n_1,\dots,n_k$.
Then the Cantor normal form of $\alpha\otimes\delta$ has leading term
$\omega^{\beta\oplus\delta_1}n_1$. If $\alpha<\varepsilon_0$
then $\beta<\omega^\beta=\alpha$, and $\delta_1<\alpha=\alpha'$, hence
$\beta\oplus\delta_1<\alpha$, since $\alpha$ is additively indecomposable.
If $\alpha\geq\varepsilon_0$, then $\beta \leq \omega^\beta<\omega^{\beta+1}$ 
and $\delta_1<\omega\otimes\alpha=\omega^{\beta+1}$, hence
$\beta\oplus\delta_1<\omega^{\beta+1}=\omega\otimes\alpha$, since
$\omega^{\beta+1}$ is additively indecomposable. In both cases
we have $\beta\oplus\delta_1<\alpha'$, hence 
$\alpha\otimes\delta<\omega^{\alpha'}$ as desired.
\end{proof}

\begin{remarkunnumbered}
By Remark~\ref{multiset} and Proposition~\ref{oProp}, the lemma we just proved
also implies that
$\ot\bigl(S^\diamond,\preceq\!\!\!\preceq\bigr)\leq
\omega^{\alpha'}$ for any Noetherian ordered set $(S,\leq)$ of type $\alpha$.
See \cite{Weiermann} for a proof of the slightly better bound
$\ot\bigl(S^\diamond,\preceq\!\!\!\preceq\bigr)\leq
\omega^{\alpha}$.
\end{remarkunnumbered}

\subsection*{Some possibilities for totally ordering monomial ideals}
By Lemma~\ref{Decr-Lemma}, we have $\cF(\Nn^{m})\cong
\Decr\bigl(\Nn,\cF(\Nn^{m-1})\bigr)$ for $m>1$:
Every final segment $F$ of $\Nn^m$
can be written as the disjoint union
\begin{equation}\label{Decomp}
F=(F_0\times\{0\}) \cup (F_1\times\{1\}) \cup \cdots\cup (F_j\times\{j\}) \cup\cdots,
\end{equation}
where 
$$F_j:=\bigl\{(e_1,\dots,e_{m-1})\in\Nn^{m-1}:(e_1,\dots,e_{m-1},j)\in F\bigr\},$$
a final segment of $\Nn^{m-1}$ (possibly empty). 
In the notation introduced in the proof of Lemma~\ref{Decr-Lemma},
$F_j=\varphi_F(j)$ for all $j\in\Nn$.
We have $(F_0,F_1,\dots)\in \cF(\Nn^{m-1})^{(\geq)}$, that is,
$F_0\subseteq F_1\subseteq\cdots$ is an ascending chain of 
final segments of $\Nn^{m-1}$ (and hence becomes eventually stationary).
Moreover, $F\supseteq G$ if and only if $F_j\supseteq G_j$ for all $j$, 
that is, if and only if $(F_0,F_1,\dots) \leq (G_0,G_1,\dots)$ in the 
ordering of $\cF(\Nn^{m-1})^{(\geq)}$.
The decomposition \eqref{Decomp} for final segments of ${\cal F}(\Nn^m)$ can
be used to explicitly construct a {\it total}\/ ordering $\trianglelefteq$
on ${\cal F}(\Nn^m)$ which extends $\supseteq$.
By Corollary~\ref{Noetherian-Cor}
and Proposition~\ref{Folk}, this
ordering will then be a well-ordering. 
For the construction, we proceed as follows,
by induction on $m$: 
\begin{enumerate}
\item If $m=1$, then $F\trianglelefteq G :\Longleftrightarrow F\supseteq G$.
\item Let $m>1$, and suppose we have already constructed a total
ordering $\trianglelefteq$ on ${\cal F}(\Nn^{m-1})$. We then put
$F\trianglelefteq G$ if and only if $$(F_0,F_1,\dots)\trianglelefteq_{\text{lex}}
(G_0,G_1\dots)$$  in the lexicographic ordering on
${\cal F}(\Nn^{m-1})^{(\geq)}$ induced by $\trianglelefteq$.
(That is, $F\trianglelefteq G$ if and only if either $F=G$, or 
there is $j\in\Nn$ with
$F_0=G_0,\dots,F_{j-1}=G_{j-1}, F_{j} \vartriangleleft G_{j}$.)
\end{enumerate}
By induction on $m$ it follows easily that 
$F\supseteq G \Rightarrow F\trianglelefteq G$, for all $F,G\in{\cal F}(\Nn^m)$.
The empty final segment is the largest and 
the final segment $\Nn^m$ the smallest element of ${\cal F}(\Nn^m)$.

We shall not try to compute here the order type of $\bigl(\cF(\Nn^m),
\trianglelefteq\bigr)$ for general $m$. Let us just point out:

\begin{lemma}
$\ot\bigl(\cF(\Nn^2),\trianglelefteq\bigr)=\omega^{\omega+1}+1$.
\end{lemma}
\noindent
In order to see this, suppose 
that $S$ is a well-ordering. Then the restriction of the
lexicographic ordering on $S^\omega$ to
$S^{(\geq)}$ is a well-ordering which extends the product ordering. In the
next proposition we
compute the order type of $S^{(\geq)}$ in terms of the order type
of $S$. For $S=\cF(\Nn)$ this yields the lemma above.
We may assume that $S=\alpha$ is an ordinal.

\begin{prop}\label{Ordertype}
Let $\alpha$ be an ordinal. Then:
$$\ot\bigl(\alpha^{(\geq)}\bigr)=\begin{cases}
\alpha & \text{if $\alpha=0$ or $\alpha=1$} \\
\omega^{\alpha-1}+1 & \text{if $2\leq\alpha<\omega$} \\
\omega^\alpha & \text{if $\alpha\geq\omega$ is a limit} \\
\omega^\alpha + 1 & \text{if $\alpha\geq\omega$ is a successor.}
\end{cases}$$
\end{prop}
\begin{proof}
Clearly $\ot\bigl(0^{(\geq)}\bigr)=0$.
Observe that $\beta^{(\geq)}$ is an initial segment of $\gamma^{(\geq)}$
for any $\beta<\gamma$, so
\begin{equation}\label{Ord1}
\ot\bigl(\alpha^{(\geq)}\bigr)=\bigcup_{\beta<\alpha} 
\ot\bigl(\beta^{(\geq)}\bigr)\qquad\text{if $\alpha$ is a limit ordinal.}
\end{equation} 
Moreover,
$$(\alpha+1)^{(\geq)}=\bigcup_{i<\omega} B_i \cup \bigl\{(\alpha,\alpha,\dots)\bigr\},$$
where $B_i$ is the set of decreasing sequences in $\alpha+1$ that begin with
exactly $i$ many $\alpha$'s. Hence each $B_i$ is isomorphic to $\alpha^{(\geq)}$
(as ordered set), thus
\begin{equation}\label{Ord2}
\ot\bigl((\alpha+1)^{(\geq)}\bigr)=\underbrace{\ot\bigl(\alpha^{(\geq)}\bigr)+
\ot\bigl(\alpha^{(\geq)}\bigr)+\cdots}_{\text{$\omega$ many times}}+1=
\ot\bigl(\alpha^{(\geq)}\bigr)\omega+1.
\end{equation}
The formula for $\ot\bigl(\alpha^{(\geq)}\bigr)$ follows
by transfinite induction, using the relations
\eqref{Ord1} and \eqref{Ord2}.
\end{proof}

\noindent
The well-ordering $\trianglelefteq$ of $\cF(\Nn^m)$
introduced above has several disadvantages.
Most severely, from a practical point of view, suppose  $F$ and $G$ are two
final segments of $\Nn^m$,
given in terms of finite sets of generators, and we want to compare 
$F$ and $G$ with respect to $\trianglelefteq$. So we need to compute
representations \eqref{Decomp} for $F$ and $G$, and lexicographically compare
the resulting sequences of final segments of $\Nn^{m-1}$. This
gives rise to a computationally demanding recursion on $m$.
Sometimes, however, we have access to the Hilbert-Samuel polynomials of
monomial ideals (since they are needed for an auxiliary computation, say).
In this case, we may  use a variant of the ordering $\trianglelefteq$ 
for which comparing $F,G\in\cF(\Nn^m)$ can be done in a more efficient way: 
By Section~\ref{InvariantsSection} 
we obtain a well-ordering $\leq$ of $\cF(\Nn^m)$ extending $\supseteq$ with
minimal possible order type $\omega^m+1$ by defining
$$F \leq G \quad :\Longleftrightarrow\quad\text{$p_F \prec p_G$, or
$p_F=p_G$ and $F\trianglelefteq G$.}$$
This makes it necessary to decide $F\trianglelefteq G$ only to break ties, that is,
in case $F$ and $G$ have the same Hilbert-Samuel polynomial.

\subsection*{The Kleene-Brouwer ordering.}
In the rest of the paper we
study another ordering of monomial ideals which has the advantage that
comparison of monomial ideals specified by sets of generators is extremely
easy.

\begin{definition}
Let $(U,{\leq})$ be a totally ordered set. 
We define the {\it Kleene-Brouwer ordering}\/ $\leq_\KB$
of the tree $U^\ast$ as follows:
If $s=(s_1,\dots,s_m)$, $t=(t_1,\dots,t_n)$, then $s\leq_\KB t$ if and only if
either
\begin{enumerate}
\item $s\sqsupseteq t$, or
\item $(s_1,\dots,s_k) \leq_{\text{lex}} (t_1,\dots,t_k)$, where
$k=\min\{m,n\}$ and $\leq_{\text{lex}}$ denotes the lexicographic ordering
on $U^k$.
\end{enumerate}
\end{definition}
\noindent
It it easy to check that $\leq_\KB$ is a total ordering on $U^\ast$ extending
the initial segment relation $\sqsupseteq$.
We refer, e.g., to \cite{Kechris}, (2.12), for a proof of the following
fundamental fact:

\begin{lemma}
Let $(U,{\leq})$ be a well-ordered set and
$T$ a tree on $U$. Then $T$ is well-founded if and only if
the Kleene-Brouwer ordering restricted to $T$ is a well-ordering. 
\rom{(}In this case, we write $\ot_{\KB}(T)$ for the order type of
$\leq_\KB$.\rom{)} \qed
\end{lemma}

\noindent
Let now $(S,\leq)$ be a Noetherian ordered set, and
fix a total ordering $\leq'$ on $S$ extending $\leq$. As in 
Section~\ref{NoetherianOrderedSetsSection}
we let $\Ant_{\leq'}(S)$ be the well-founded tree
$$\Ant_{\leq'}(S) := \bigl\{(s_1,\dots,s_n)\in S^\ast:s_i\,\,||\,\, s_j
\text{ and } s_i<' s_j
\text{ for all $1\leq i<j\leq n$}\bigr\}$$
on $S$. We consider $\Ant_{\leq'}(S)$ as an
ordered set via the restriction of $\leq_\KB$.
We define a bijection $\varphi\colon\cF(S)\to \Ant_{\leq'}(S)$ by
$\varphi(F)=(a_1,\dots,a_n)$, where $a_1,\dots,a_n$ are the minimal
generators of the final segment $F$ of $S$, ordered in increasing
order with respect to $\leq'$.

\begin{lemma}
The map $\varphi\colon\cF(S)\to \Ant_{\leq'}(S)$ is strictly increasing.
\end{lemma}
\begin{proof}
Let $F\supset G$ be final segments of $S$, and let $a_1,\dots,a_n$ and
$b_1,\dots,b_m$ be the minimal generators of $F$ and $G$, respectively,
with $a_1 <' \cdots <' a_n$ and $b_1 <' \cdots <' b_m$. 
We need to show $(a_1,\dots,a_n)<_{\KB} (b_1,\dots,b_m)$.
If $(b_1,\dots,b_m) \sqsubset (a_1,\dots,a_n)$ we are done. 
Otherwise, there exists $r \le 
\min\{m,n\}$ such that $a_1=b_1,\dots,a_{r-1}=b_{r-1}$ and 
$a_r \neq b_r$. Since 
$F \supset G$, we have $a_i \le b_r$ for some $i$. 
Since $a_j = b_j$ for $j < r$ and 
$\{b_1,\dots,b_{r}\}$ is an antichain, we have $i \ge r$, hence
$a_r \le' a_i \le b_r$.
Since $a_r \neq b_r$ we have $a_r <' b_r$, and therefore
$(a_1,\dots,a_n)<_{\KB} (b_1,\dots,b_m)$ as required.
\end{proof}
\noindent
By means of the last lemma, we obtain a well-ordering on
$\cF(\Nn^m)$ extending $\supseteq$ as follows:
Fix a term ordering $\leq'$ on $\Nn^m$.
Given final segments $F$ and $G$ of $\Nn^m$, 
with minimal generators $a_1<'\cdots <'a_r$
and $b_1<'\cdots <'b_s$ (where $r,s\in\Nn$), define
$$F \leq_{\KB} G \qquad :\Longleftrightarrow\qquad (a_1,\dots,a_r) \leq_{\KB}
(b_1,\dots,b_s) \ \text{(in $\Ant_{\leq'}(\Nn^m)$).}$$
We shall call the well-ordering $\leq_{\KB}$ of $\cF(\Nn^m)$
the {\it Kleene-Brouwer ordering}\/ of $\cF(\Nn^m)$ (induced by $\leq'$),
and we put $\ot_{\KB}\bigl(\cF(\Nn^m)\bigr):=
\ot_\KB\bigl(\Ant_{\leq'}(\Nn^m)\bigr)$. If $I$ and $J$ are monomial
ideals in $K[X_1,\dots,X_m]$ (where $K$ is a field)
corresponding to final segments
$F$ and $G$ of $\Nn^m$, respectively, we put $I\leq_{\KB} J$ if
$F\leq_{\KB} G$. This yields a well-ordering on the set of monomial
ideals of $K[X]$ which extends $\supseteq$.

For lex-segments, the
Kleene-Brouwer ordering induced by the degree-lexicogra\-phic ordering
has an alternative description:

\begin{example}
Suppose that $\leq'$ is the degree-lexicographic ordering of $\Nn^m$,
and let $E=(a_1,\dots,a_r)$ with $a_1<'\cdots<'a_r$ be a lex-segment of
$\Nn^m$, and let $F=(b_1,\dots,b_s)$ with $b_1<'\cdots<'b_s$
be any final segment of $\Nn^m$. Then
$$E <_\KB F \quad\Rightarrow\quad E_0=F_0, \dots, E_{d-1}=F_{d-1}, E_d
\supset F_d\text{ for some $d\in\Nn$.}$$
\end{example}
\begin{proof}
Suppose that $E<_\KB F$. If $E\supset F$, we are done. Otherwise, there
is $t<\min\{r,s\}$ such that $a_1=b_1,\dots,a_t=b_t,a_{t+1}<'b_{t+1}$.
Put $d=\abs{a_{t+1}}$. Then $a_{t+1}\notin F_d$: otherwise we have
$a_{t+1}\geq b_i$ for some $i$. Since the $a_1,\dots,a_s$ form an
antichain (with respect to $\leq$) and $a_i=b_i$ for $1\leq i\leq t$,
we have $i>t$ and so $a_{t+1}\geq' b_{t+1}$, a contradiction. 
Moreover $F_j\subseteq E_j$ for all $0\leq j\leq d$: Let $x\in F_j$,
so $x\geq b_i$ for some $i\in\{1,\dots,s\}$. If $1\leq i\leq t$ we are
done (since $a_i=b_i$), so suppose $i>t$. Then $x\geq' b_{t+1}>'a_{t+1}$,
so $\abs{x}\geq a_{t+1}=d$, hence $j=d$. Since $E$ is a lex-segment and
$x>' a_{t+1}$ we get $x\in E_d$ as claimed. Finally, we have $E_j\subseteq
F_j$ for $0\leq j<d$: If $y\in E_j$ for $j\in\{1,\dots,d-1\}$, then $y\geq
a_i$ for some $i\in\{1,\dots,t\}$, hence $y\in F_j$.
\end{proof}

\begin{historicalremark}
The Kleene-Brouwer ordering of a tree plays an important role 
in descriptive set theory and recursion theory.
It appears for the first time in the work of Brouwer \cite{Brouwer}
(in his proof that intuitionistically, every real function is uniformly
continuous on closed intervals) and Lusin-Sierpinski \cite{Lusin-Sierpinski},
and was later used by Kleene \cite{Kleene}. (See the remarks in
\cite{Moschovakis}, p.~270.) A variant of the
Kleene-Brouwer ordering was independently discovered by Ritt in his seminal 
work on differential algebra, in his definition of the rank of
characteristic sets.
(See \cite{Kolchin}, p.~81.)
\end{historicalremark}

\subsection*{An upper bound for $\ot_\KB$.}
We want to investigate the order-theoretic complexity of $\ot_\KB$. We first
establish an upper bound on $\ot_\KB$.

\begin{notation}
Given an ordinal $\alpha$ and a sequence $(\alpha_n)_{n\in\Nn}$ of
ordinals, we write $$\alpha=\limsup_n \alpha_n$$ if $\alpha=
\sup\{\alpha_n : n\in\Nn\}$, and for every $n_0$ and $\beta<\alpha$
there exists $n\geq n_0$ with $\beta<\alpha_n$. (Equivalently,
$\alpha=\limsup_n \alpha_n$ if and only if $\alpha=\sup\bigl\{\alpha_{n_i}:
i\in\Nn\bigr\}$ for some increasing subsequence $\alpha_{i_0}\leq\alpha_{i_1}
\leq\cdots$ of $(\alpha_n)$.)
\end{notation}

\noindent
In the following, $U$ and $V$ will denote {\it countable}\/ sets.
For the purpose of this section, let us call a tree $T$ on 
$U$ {\it universal}\/ if $T$ is well-founded, and
every node $a$ of $T$ which is not a leaf has infinitely
many successors $a_0,a_1,\dots$, and $\hgt(a)=\limsup_n
\hgt(a_n)+1$.  Note that the property of being universal 
is preserved under passing to subtrees.

By Lemma~\ref{Kechris-Lemma}, if $S$ and $T$ are well-founded trees with
$\rk(S)\leq\rk(T)$, then there exists an increasing length-preserving
map $S\to T$. If $T$ is universal, we have the following result
(justifying our choice of terminology):

\begin{lemma}\label{Universal}
Let $S$ and $T$ be trees on $U$ and $V$, respectively. 
If $T$ is universal, then the following are equivalent:
\begin{enumerate}
\item $S$ is well-founded with $\rk(S)\leq\rk(T)$.
\item There exists a length-preserving embedding $S\to T$.
\item There exists a strictly increasing map $S\to T$.
\end{enumerate}
\end{lemma}
\begin{proof}
We prove (1)~$\Rightarrow$~(2) by induction on the rank of $T$,
the case $\rk(T)=0$ being trivial. Suppose that $\rk(T)=\alpha+1$ is
a successor. 
Since $T$ is universal, there exists a sequence
$(b_n)_{n\in\Nn}$ of pairwise distinct elements of $V$
such that $(b_0),(b_1),\dots$ are successors of $\varepsilon$ in $T$ of
height $\alpha$. 
Suppose that $S$ is well-founded with $\rk(S)\leq\rk(T)$, and
let $\bigl((a_n)\bigr)_{n<\lambda}$ (where $\lambda\leq\omega$
and $a_i\neq a_j$ for all $0\leq i<j<\lambda$) be
the successors of the root $\varepsilon$ in $S$. For every $n$, 
the subtree $S_{(a_n)}$  has rank $\leq\alpha$ and hence can be embedded 
into $T_{(b_n)}$
by a length-preserving embedding, by
induction hypothesis. 
Hence there exists a length-preserving 
embedding $a_n\,  T_{(a_n)} \to
b_n\, S_{(b_n)}$. 
Extending the union of these embeddings to a map
$S\to T$ by mapping the root of $S$ to the one of $T$ gives a length-preserving
embedding of $S$ into $T$. Finally, suppose that $\rk(T)$ is a 
limit ordinal.
Let $(b_n)_{n\in\Nn}$ be a sequence of elements of $V$
such that $\sup_n \rk(T_{(b_n)}) = \rk(T)$.
Suppose that $S$ is well-founded with $\rk(S)\leq\rk(T)$, and as before
let $\bigl((a_n)\bigr)_{n<\lambda}$ (where $\lambda\leq\omega$
and $a_i\neq a_j$ for all $0\leq i<j\leq\lambda$) be
the successors of the root $\varepsilon$ in $S$. For each 
$0\leq i <\lambda$
there exists $n_i\in\Nn$ such that $\rk(S_{(a_i)}) \leq \rk(T_{(b_{n_i})})<\rk(T)$.
Using the induction hypothesis we find a length-preserving embedding
$a_i\,  T_{(a_i)} \to
b_{n_i}\, S_{(b_{n_i})}$. Again it is not difficult to combine these
to obtain a length-preserving embedding $S\to T$ as required.
The implications (2)~$\Rightarrow$~(3) and (3)~$\Rightarrow$~(1) are clear.
\end{proof}

\begin{cor}
Suppose $U$ is infinite and  $T\neq\{\varepsilon\}$ a well-founded
tree on $U$. There exists a universal tree $T'$ on $U$ with $T\subseteq T'$
of the same rank as $T$.
\end{cor}
\begin{proof}
By Lemma~\ref{Universal},
it is enough to construct {\it some}\/ universal tree on
$U$ of rank $\rk(T)$. This is easy to accomplish by induction on $\rk(T)$.
\end{proof}

\noindent
The order type of $\leq_{\KB}$ is easy to compute for universal trees:

\begin{lemma}\label{Smooth-KB}
Let $T\neq\{\varepsilon\}$ be a universal tree on a well-ordered set $U$ of
order type $\omega$. Then $\ot_\KB(T)=\omega^{\rk(T)}+1$.
\end{lemma}
\begin{proof}
Let $a_0<a_1<\cdots$ be the successors of the root $\varepsilon$ of $T$,
listed according to their order in $U$. Note that
\begin{equation}\label{KB}
\ot_{\KB}(T) = \ot_{\KB}(T_0) + \ot_{\KB}(T_1) + \cdots + 1,
\end{equation}
where $T_n:=T_{a_n}$ is the subtree of $T$ with root at $a_n$.
We prove the lemma by induction on $\rk(T)>0$. 
The result is clear if $\rk(T)=1$. Suppose that $\rk(T)=\alpha+1$ where
$\alpha>0$.
Each $T_n$ is universal of rank $\leq\alpha$, so by induction hypothesis,
if $T_n\neq\{\varepsilon\}$, then
$\ot_{\KB}(T_n)=\omega^{\rk(T_n)}+1$ for all $n$.
Because $T$ is universal, there are infinitely many $n\in\Nn$
such that $\rk(T_n)=\alpha$. Since $\omega^\beta+\omega^\gamma=
\omega^\gamma$ whenever $\gamma>\beta$, it follows readily from
\eqref{KB} that $\ot_{\KB}(T)=\omega^{\rk(T)}+1$.
Suppose now that $\rk(T)$ is a limit ordinal. Hence
$\rk(T)=\sup\bigl\{\rk(T_n)+1:n\in\Nn\bigr\}$
is the limit of a strictly increasing subsequence of
$\bigl(\rk(T_n)\bigr)$. So by \eqref{KB} and induction hypothesis
we have
$$\ot_{\KB}(T) = \left(\bigcup_n \omega^{\rk(T_{n})}\right) + 1 = \omega^{\bigcup_n\rk(T_{n})} + 1 =
\omega^{\rk(T)} + 1,$$
as desired.
\end{proof}

%\begin{cor}
%Let $U$ and $V$ be equipped with well-orderings of
%order type $\omega$. If $S$ and $T$ are trees on $U$ and $V$, with  
%$S$ well-founded and $T$ smooth, such that $\rk(S)\leq\rk(T)$, then
%there exists an embedding $S\to T$ with respect to the Kleene-Brouwer 
%orderings on $S$ and $T$. \qed
%\end{cor}

\begin{cor}\label{KB-Upperbound}
For every well-founded
tree $T\neq\{\varepsilon\}$ on a well-ordered set $U$ of order
type $\omega$, we have $\ot_\KB(T) \leq\omega^{\rk(T)}+1$. \qed
\end{cor}
\noindent
Clearly it may happen that $\ot_\KB(T)<\omega^{\rk(T)}+1$, for example
if  $T\neq\{\varepsilon\}$ is finite. Another (infinite) example is given by
the tree
$T=\bigl\{ (i,i,\dots,i)\in \Nn^i : i\in\Nn\bigr\}$ of rank $\omega$ 
on $\Nn$.

\subsection*{The order type of the Kleene-Brouwer ordering.}
We now investigate the order type of the Kleene-Brouwer ordering on
$\cF(\Nn^m)$ in the case where $\leq'$ has order type $\omega$.
By Lemma~\ref{Ant-Prime} and
Corollary~\ref{KB-Upperbound} we obtain the upper bound
$\ot_\KB\bigl(\cF(\Nn^m)\bigr) \leq \omega^{\omega^{m-1}}+1$. We will show:

\begin{prop}\label{KB-Prop}
The tree $\Ant_{\leq'}(\Nn^m)$ contains a
universal tree on $\Nn^m$ with the same rank.
\end{prop}

\noindent
Using Lemma~\ref{Smooth-KB}, this immediately yields:

\begin{cor}\label{KB-Cor}
$\ot_\KB\bigl(\cF(\Nn^m)\bigr) = \omega^{\omega^{m-1}}+1$. \qed
\end{cor}
\noindent
Before we begin the proof, let us introduce some notations:
Given an element $\nu\in\Nn^m$ we will denote by $\tau_\nu$
the translation $$x\mapsto x+\nu \colon\Nn^m\to\Nn^m,$$ and
given a natural number $n$ we denote by $\iota_n$
the map $$\nu\mapsto\nu\, n\colon\Nn^m\to\Nn^{m+1}.$$
By component-wise application,
the map $\tau_\nu$ gives rise to a map 
$(\Nn^m)^\ast\to (\Nn^m)^\ast$
and $\iota_n$ gives rise to a map $(\Nn^m)^\ast\to (\Nn^{m+1})^\ast$,
denoted by the same symbols. 
We have $$\tau_\nu\bigl(\Ant_{\leq'}(\Nn^m)\bigr)
\subseteq\Ant_{\leq'}(\Nn^m)\quad\text{ and }\quad
\iota_n\bigl(\Ant_{\leq'}(\Nn^m)\bigr)
\subseteq\Ant_{\leq'}(\Nn^{m+1}).$$
For a sequence $a=(a_1,\dots,a_n)\in
(\Nn^m)^\ast$ we put $|a|=|a_1|+\cdots+|a_n|$.

\begin{proof}[Proof \rom{(Proposition~\ref{KB-Prop})}]
We proceed by induction on $m=1,2,\dots$. The case $m=1$ is trivial, since
$\Ant_{\leq'}(\Nn)$ itself is universal. Suppose that $m>1$
and let $U\subseteq\Ant_{\leq'}(\Nn^{m-1})$ be a universal tree of 
rank $\omega^{m-2}$. Put $T:= \Ant_{\leq'}(\Nn^{m})$. 
For any $k\geq 1$, we have $$T_{(0,\dots,0,k)}=
\Ant_{\leq'}\bigl(\Nn^{m-1}\times\{0,\dots,k-1\}\bigr),$$ hence
$\rk(T_{(0,\dots,0,k)})=\omega^{m-2}k$ by Lemma~\ref{Ant-Prime}.
Therefore it suffices to show that $T_{(0,\dots,0,k)}$ contains a universal
tree of rank $\omega^{m-2}k$. Starting with $V_0=\{\varepsilon\}$, we 
construct this tree in $k$ steps. 
Suppose that $V_i$ ($0\leq i<k$) is a universal
tree which is contained in $$T_{(0,\dots,0,k)}\cap \bigl(\Nn^{m-1}\times
\{k-i,\dots,k-1\}\bigr)^\ast$$ and has rank $\omega^{m-2}i$.
For each leaf $a$ of $V_i$ choose an element $v_a$ of
$\Nn^{m-1}$ with $|v_a|>|a|+i+1-k$. It is easy to check that
$\abs{v_a}$ is large enough to guarantee that every node of the universal tree
$a\ \iota_{k-i-1}\bigl(\tau_{v_a}(U)\bigr)$
is an antichain in $\Nn^m\times\{k-i-1,\dots,k-1\}$ arranged in $\leq'$-increasing order. Hence
the tree $$V_{i+1}:=V_i\cup \bigcup_a 
a\ \iota_{k-i-1}\bigl(\tau_{v_a}(U)\bigr)$$
(where the union runs over all leafs $a$ of $V_i$) is contained in
$$T_{(0,\dots,0,k)}\cap \bigl(\Nn^{m-1}\times
\{k-i-1,\dots,k-1\}\bigr)^\ast.$$
So $V_{i+1}$ is simply the tree obtained by ``implanting'' a copy of
$U$ (that is, the tree $a\ \iota_{k-i-1}(\tau_{v_a}(U))$) at $a$, for each leaf
$a$ of $V_i$. It is immediate that each non-leaf node of $V_{i+1}$
has $\omega$ many successors. Note also that the heights of nodes of $V_{i+1}$
that are coming from $V_i$ will increase by
$\rk(U)$ while the heights of nodes coming from $U$ will remain
unchanged.  More precisely we have: $\hgt_{V_{i+1}}(v) = \hgt_{V_{i}}(v)
+ \rk(U)$ if $v \in V_{i}$ and $\hgt_{V_{i+1}}(v) = \hgt_{U}(u)$ for
$v=a\ \iota_{k-i-1}(\tau_{v_a}(u))$ where $a$ is a leaf of $V_{i}$ and
$u \in U$. This observation clearly implies the $\limsup$ condition
hence universality for $V_{i+1}$, and
$$
  \rk(V_{i+1}) = \rk(V_{i}) + \rk(U) = \omega^{m-2}(i+1).
$$
Now one sees that the last tree $V_k$ constructed in this way has the
desired properties.
%If $v\in V_i$, then $\hgt_{V_{i+1}}(v) = \hgt_{V_i}+\hgt(U)$, and if
%$u\in U$, $a\in V_i$, then 
%$\hgt_{V_{i+1}}(v)=\hgt_{V_i}(v)+\rk(U)$ for
%$v=a \iota_{k-i-1}\bigl(\tau_{v_a}(u)\bigr)$. 
%It follows that $V_{i+1}$ is universal and has rank
%$$\rk(V_{i+1})=\rk(V_i)+\rk(U)=\omega^{m-2}(i+1)$$ as required.
%The last tree $V_k$ constructed in this way hence has the desired properties.
\end{proof}

\noindent
Combining Corollary~\ref{KB-Cor} and \eqref{UpperBound} we obtain as promised our 
estimates
on the type of the ordered set of monomial ideals:

\begin{cor}
$\omega^{\omega^{m-1}}+1 \leq \ot\bigl(\cF(\Nn^m)\bigr) \leq 
\omega^{(\omega+1)^{\otimes m}}$. \qed
\end{cor}

%\begin{remark}
%Fix a Hilbert polynomial $P$ of a non-zero final segment of $\Nn^m$, and
%consider $S=\bigl\{E\in\cF(\Nn^m): P_E=P\bigr\}$.
%What can one say about the Hilbert function of the smallest element
%of $S$ in the Kleene-Brouwer ordering?
%\end{remark}

\bibliographystyle{amsplain}
\bibliography{rd.bib}

\providecommand{\bysame}{\leavevmode\hbox to3em{\hrulefill}\thinspace}
\providecommand{\MR}{\relax\ifhmode\unskip\space\fi MR }
\providecommand{\MRhref}[2]{%
  \href{http://www.ams.org/mathscinet-getitem?mr=#1}{#2}
}
\providecommand{\href}[2]{#2}
\begin{thebibliography}{10}

\bibitem{Abraham}
U.~Abraham, \emph{A note on {D}ilworth's theorem in the infinite case}, Order
  \textbf{4} (1987), no.~2, 107--125.

\bibitem{Abraham-Bonnet}
U.~Abraham and R.~Bonnet, \emph{Hausdorff's theorem for posets that satisfy the
  finite antichain property}, Fund. Math. \textbf{159} (1999), no.~1, 51--69.

\bibitem{maschenb-hemmecke}
M.~Aschenbrenner and R.~Hemmecke, \emph{A finiteness theorem in stochastic
  integer programming}, in preparation, 2003.

\bibitem{maschenb-pong-rd}
M.~Aschenbrenner and W.~Y. Pong, \emph{The differential order}, in preparation,
  2003.

\bibitem{Bachmann}
H.~Bachmann, \emph{Transfinite {Z}ahlen}, Ergebnisse der Mathematik und ihrer
  Grenzgebiete, vol.~1, Springer-Verlag, Berlin, 1955.

\bibitem{BeckerWeispfenning}
T.~Becker and V.~Weispfenning, \emph{Gr\"obner {B}ases}, Graduate Texts in
  Mathematics, vol. 141, Springer-Verlag, New York, 1993.

\bibitem{BCR}
J.~Bochnak, M.~Coste, and M.-F. Roy, \emph{G\'eom\'etrie {A}lg\'ebrique
  {R}\'eelle}, Ergebnisse der Mathematik und ihrer Grenzgebiete (3), vol.~12,
  Springer-Verlag, Berlin, 1987.

\bibitem{Bonnet-Pouzet}
R.~Bonnet and M.~Pouzet, \emph{Extension et stratification d'ensembles
  dispers\'es}, C. R. Acad. Sci. Paris S\'er. A-B \textbf{268} (1969),
  A1512--A1515.

\bibitem{Brookfield}
G.~Brookfield, \emph{The length of {N}oetherian modules}, Comm. Algebra
  \textbf{30} (2002), no.~7, 3177--3204.

\bibitem{Brouwer}
L.~E.~J. Brouwer, \emph{Beweis, da\ss{} jede volle {F}unktion gleichm\"assig
  stetig ist}, Nederl. Akad. Wetensch. Proc. Ser. B \textbf{27} (1924),
  189--193.

\bibitem{Bruns-Herzog}
W.~Bruns and J.~Herzog, \emph{Cohen-{M}acaulay {R}ings}, Cambridge Studies in
  Advanced Mathematics, vol.~39, Cambridge University Press, Cambridge, 1993.

\bibitem{dJP}
D.~de~Jongh and R.~Parikh, \emph{Well-partial orderings and hierarchies},
  Indag. Math. \textbf{39} (1977), no.~3, 195--207.

\bibitem{Dershowitz}
N.~Dershowitz, \emph{Termination of rewriting}, J. Symbolic Comput. \textbf{3}
  (1987), no.~1--2, 69--115.

\bibitem{Dilworth}
R.~P. Dilworth, \emph{A decomposition theorem for partially ordered sets}, Ann.
  of Math. (2) \textbf{51} (1950), 161--166.

\bibitem{Dube}
T.~Dub{\'e}, \emph{The structure of polynomial ideals and {G}r\"obner bases},
  SIAM J. Comput. \textbf{19} (1990), no.~4, 750--775.

\bibitem{Eisenbud}
D.~Eisenbud, \emph{Commutative {A}lgebra with a {V}iew {T}oward {A}lgebraic
  {G}eometry}, Graduate Texts in Mathematics, vol. 150, Springer-Verlag, New
  York, 1995.

\bibitem{Fraisse}
R.~Fra\"\i{}ss\'e, \emph{Theory of {R}elations}, revised ed., Studies in Logic
  and the Foundations of Mathematics, vol. 145, North-Holland Publishing Co.,
  Amsterdam, 2000.

\bibitem{Friedman}
H.~Friedman, \emph{The {A}ckermann function in elementary algebraic geometry},
  manuscript, 1999.

\bibitem{Gallier}
J.~H. Gallier, \emph{What's so special about {K}ruskal's theorem and the
  ordinal {$\Gamma\sb 0$}? {A} survey of some results in proof theory}, Ann.
  Pure Appl. Logic \textbf{53} (1991), no.~3, 199--260.

\bibitem{Higman}
G.~Higman, \emph{Ordering by divisibility in abstract algebras}, Proc. London
  Math. Soc. (3) \textbf{2} (1952), 326--336.

\bibitem{Kechris}
A.~S. Kechris, \emph{Classical {D}escriptive {S}et {T}heory}, Graduate Texts in
  Mathematics, vol. 156, Springer-Verlag, New York, 1995.

\bibitem{Kleene}
S.~C. Kleene, \emph{On the forms of the predicates in the theory of
  constructive ordinals. \rom{II}}, Amer. J. Math. \textbf{77} (1955),
  405--428.

\bibitem{Kolchin}
E.~Kolchin, \emph{Differential {A}lgebra and {A}lgebraic {G}roups}, Pure and
  Applied Mathematics, vol.~54, Academic Press, New York-London, 1973.

\bibitem{ddp}
M.~V. Kondratieva, A.~B. Levin, A.~V. Mikhalev, and E.~V. Pankratiev,
  \emph{Differential and {D}ifference {D}imension {P}olynomials}, Mathematics
  and its Applications, vol. 461, Kluwer Academic Publishers, Dordrecht, 1999.

\bibitem{Kruskal}
J.~B. Kruskal, \emph{The theory of well-quasi-ordering: {A} frequently
  discovered concept}, J. Combinatorial Theory Ser. A \textbf{13} (1972),
  297--305.

\bibitem{Kriz-Thomas}
I.~K\v{r}\'{\i}\v{z} and R.~Thomas, \emph{Ordinal types in {R}amsey theory and
  well-partial-ordering theory}, in \cite{Nesetril-Roedl}, pp.~57--95.

\bibitem{Lusin-Sierpinski}
N.~Lusin and W.~Sierpinski, \emph{Sur un ensemble non mesurable {B}}, J. Math.
  Pures Appl. \rom{(}9\rom{)} \textbf{2} (1923), 53--72.

\bibitem{Macaulay}
F.~S. Macaulay, \emph{Some properties of enumeration in the theory of modular
  systems}, Proc. London Math. Soc. \textbf{26} (1927), no.~2, 531--555.

\bibitem{Maclagan}
D.~Maclagan, \emph{Antichains of monomial ideals are finite}, Proc. Amer. Math.
  Soc. \textbf{129} (2001), no.~6, 1609--1615 (electronic).

\bibitem{Martin-Scott}
U.~Martin and E.~Scott, \emph{The order types of termination orderings on
  monadic terms, strings and multisets}, J. Symbolic Logic \textbf{62} (1997),
  no.~2, 624--635.

\bibitem{Moreno-Socias}
G.~Moreno~Soc{\'\i}as, \emph{Length of polynomial ascending chains and
  primitive recursiveness}, Math. Scand. \textbf{71} (1992), 181--205.

\bibitem{Moschovakis}
Y.~N. Moschovakis, \emph{Descriptive {S}et {T}heory}, Studies in Logic and the
  Foundations of Mathematics, vol. 100, North-Holland Publishing Co.,
  Amsterdam, 1980.

\bibitem{NW2}
C.~St. J.~A. Nash-Williams, \emph{On well-quasi-ordering finite trees}, Proc.
  Cambridge Philos. Soc. \textbf{59} (1963), 833--835.

\bibitem{NW}
\bysame, \emph{On well-quasi-ordering infinite trees}, Proc. Cambridge Philos.
  Soc. \textbf{61} (1965), 697--720.

\bibitem{Nesetril-Roedl}
R.~Ne\v{s}et\v{r}il and V.~R\"odl (eds.), \emph{Mathematics of {R}amsey
  {T}heory}, Algorithms and Combinatorics, vol.~5, Springer-Verlag,
  Berlin-Heidelberg-New York, 1990.

\bibitem{Perles}
M.~Perles, \emph{On {D}ilworth's theorem in the infinite case}, Israel J. Math.
  \textbf{1} (1963), 108--109.

\bibitem{Rado}
R.~Rado, \emph{Partial well-ordering of sets of vectors}, Mathematika
  \textbf{1} (1954), 89--95.

\bibitem{Robbiano}
L.~Robbiano, \emph{Term orderings on the polynomial ring}, EUROCAL '85, Vol.\ 2
  (Linz, 1985), Lecture Notes in Comput. Sci., vol. 204, Springer, Berlin,
  1985, pp.~513--517.

\bibitem{Rust}
C.~J. Rust, \emph{Rankings of {D}erivatives for {E}limination {A}lgorithms and
  {F}ormal {S}olvability of {A}nalytic {P}artial {D}ifferential {E}quations},
  Ph.D. thesis, University of Chicago, 1998.

\bibitem{Schmidt}
D.~Schmidt, \emph{Well-{P}artial {O}rderings and {T}heir {M}aximal {O}rder
  {T}ypes}, Ha\-bi\-li\-tations\-schrift, Mathematisches Institut der
  Universit\"at Heidelberg, 1979.

\bibitem{Simpson-Book}
S.~Simpson, \emph{Subsystems of {S}econd {O}rder {A}rithmetic}, Perspectives in
  Mathematical Logic, Springer-Verlag, Berlin, 1999.

\bibitem{Sit}
W.~Sit, \emph{Well ordering of certain numerical polynomials}, Trans. Amer.
  Math. Soc. \textbf{212} (1975), no.~1, 37--45.

\bibitem{Sturmfels}
B.~Sturmfels, \emph{Gr\"obner {B}ases and {C}onvex {P}olytopes}, University
  Lecture Series, vol.~8, American Mathematical Society, Providence, RI, 1996.

\bibitem{vdDries-Ehrlich-Erratum}
L.~van~den Dries and P.~Ehrlich, \emph{Erratum to
  \rom{{\cite{vdDries-Ehrlich}}}}, Fund. Math. \textbf{168} (2001), no.~3,
  295--297.

\bibitem{vdDries-Ehrlich}
\bysame, \emph{Fields of surreal numbers and exponentiation}, Fund. Math.
  \textbf{167} (2001), no.~2, 173--188.

\bibitem{vdH:phd}
J.~van~der Hoeven, \emph{Asymptotique automatique}, Ph.D. thesis, \'Ecole
  Polytechnique, Paris, 1997.

\bibitem{Weiermann}
A.~Weiermann, \emph{Proving termination for term rewriting systems}, Computer
  {S}cience {L}ogic \rom{(}{B}erne, 1991\rom{)}, Lecture Notes in Comput. Sci.,
  vol. 626, Springer, Berlin, 1992, pp.~419--428.

\bibitem{Weispfenning-Orders}
V.~Weispfenning, \emph{Admissible orders and linear forms}, ACM SIGSAM Bulletin
  \textbf{21} (1987), no.~2, 16--18.

\bibitem{Wolk}
E.~S. Wolk, \emph{Partially well ordered sets and partial ordinals}, Fund.
  Math. \textbf{60} (1967), 175--186.

\end{thebibliography}

\end{document}